\documentclass[12pt]{amsart}
\usepackage{epsf,righttag,latexsym,amssymb,amscd,xspace}
\newcommand{\titlestring}{Resolutions of Orbifold Singularities and
  the Transportation Problem on the {M}c{K}ay Quiver}
\newcommand{\runningheadstring}{Orbifold Singularities and McKay Flows}
\newcommand{\msc}{(1991): 14M25 (Primary) 05C35, 90C08, 32S45 (Secondary)}
\newcommand{\name}{Alexander V.\ Sardo Infirri}
\newtheorem{thm}{Theorem}[section]
\newtheorem{nonumberthm}{Theorem}  
\newtheorem{nonumbercor}{Corollary}  
\newtheorem{prop}[thm]{Proposition}
\newtheorem{cor}[thm]{Corollary}
\newtheorem{lemma}[thm]{Lemma}

 \theoremstyle{definition}

\newtheorem{dfn}[thm]{Definition}
\newtheorem{problem}[thm]{Problem}
\newtheorem{conj}{Conjecture}
  
\newtheorem{example}[thm]{Example}

\newtheorem{notation}[thm]{Notation}

\theoremstyle{remark}

\newtheorem{rmk}[thm]{Remark}      %   \renewcommand{\thermk}{}
\newtheorem{question}{Question} 

\newtheorem{fact}[thm]{Fact}

\newcommand{\frk}{\mathfrak{k}}

\newcommand{\Adm}{{{\Lambda^{0,0}_{\R}}}}
\newcommand{\dsh}{\mbox{---}}
\newcommand{\C}{{\mathbb C}}
\newcommand{\CX}{{\C^*}}
\newcommand{\R}{{\mathbb R}}
\newcommand{\Z}{{\mathbb Z}} 
\newcommand{\Q}{{\mathbb Q}} 
\newcommand{\N}{{\mathbb N}}

\newcommand{\PP}{{\mathbb P}}
\newcommand{\nneg}{_{\ge 0}}

\newcommand{\T}{\protect\overline{T}{}}
\newcommand{\Mtw}{\widetilde{M}}

\newcommand{\Ptw}{\widetilde{P}}
\newcommand{\cC}{{\mathcal C}}
\newcommand{\cD}{{\mathcal D}}
\newcommand{\cE}{{\mathcal E}}

\newcommand{\cK}{{\mathcal K}}

\newcommand{\cM}{{\mathcal M}}
\newcommand{\cN}{{\mathcal N}}
\newcommand{\cO}{{\mathcal O}}

\newcommand{\cQ}{{\mathcal Q}}

\newcommand{\cS}{{\mathcal S}}
\newcommand{\cT}{{\mathcal T}}

\newcommand{\cW}{{\mathcal W}}

\newcommand{\bgamma}{\boldsymbol{\gamma}}

\newcommand{\bdR}{{\mathbf R}}
\newcommand{\bdV}{{\mathbf V}}

\newcommand{\Tn}{T^\Pi}
\newcommand{\clos}{\overline}
\newcommand{\SC}{\clos{S}}
\newcommand{\TC}{\clos{T}}

\newcommand{\alphatw}{\tilde \alpha}

\newcommand{\ga}{\Gamma}

\newcommand{\gahat}{\widehat\Gamma}

\newcommand{\gitquot}[1]{/\!\!/_{\!#1}\,}

\newcommand{\ie}{i.e.\xspace}
\newcommand{\cf}{c.f.\xspace}
\newcommand{\resp}{resp.\xspace}

\newcommand{\git}{geometric invariant theory}

\newcommand{\kah}{K\"ahler}

\DeclareMathOperator{\Map}{Map}

\DeclareMathOperator{\pr}{pr}

\DeclareMathOperator{\Image}{Im}
\DeclareMathOperator{\rk}{rank}
\DeclareMathOperator{\spn}{span}

\DeclareMathOperator{\trace}{trace}
\DeclareMathOperator{\Proj}{Proj}
\DeclareMathOperator{\Spec}{Spec}

\DeclareMathOperator{\Hom}{Hom}
\DeclareMathOperator{\End}{End}

\DeclareMathOperator{\Rep}{Rep}

\DeclareMathOperator*{\In}{In}
\DeclareMathOperator*{\Out}{Out}
\DeclareMathOperator{\supp}{supp}
\DeclareMathOperator{\interior}{int}

\DeclareMathOperator{\ext}{ext}

\DeclareMathOperator{\Lie}{Lie}

\DeclareMathOperator{\GL}{GL}

\DeclareMathOperator{\PGL}{PGL}

\DeclareMathOperator{\SL}{SL}

\DeclareMathOperator{\U}{U}

\DeclareMathOperator{\SU}{SU}

\DeclareMathOperator{\PU}{PU}

\def\qsing#1/#2(#3){\frac{#1}{#2}(#3)} % to get 1/r(a,b,c,...,d) notation
\def\<#1>{\langle#1\rangle}

\newcommand{\trans}{{}^t\kern -0.2em}

\newcommand{\aftersub}{\nopagebreak}
\newcommand{\abs}[1]{\lvert#1\rvert}
\newcommand{\norm}[1]{\lVert#1\rVert}
\newcommand{\card}[1]{\sharp\,#1}

\newcommand{\map}[5]{$$\begin{array}{rccc} 
#1\colon & #2 & \longrightarrow & #3 \\ 
         & #4 & \longmapsto & #5
\end{array}$$}

\newcommand{\corresp}[4]{$$\begin{array}{ccc} 
#1 & \to     & #2 \\ 
#3 & \mapsto & #4
\end{array}$$}

\numberwithin{equation}{section}
\def\*#1:#2*{\S\ref{sec:#1}.\ref{sec:#1:#2}}

        {\begin{list}{}{%
                \settowidth{\labelwidth}{{\emph{#1:}}}%
                \setlength{\leftmargin}{\labelwidth+\labelsep}}}%
        {\end{list}}
\newenvironment{entry}
        {\begin{list}{}%
                {%
                  \setlength{\labelwidth}{12mm}%
                  \setlength{\leftmargin}{14mm}%
                }%
        }%
        {\end{list}}
\newlength{\Mylen}
\newcommand{\Lentrylabel}[1]{%
        \settowidth{\Mylen}{\emph{#1}}%
        \ifthenelse{\lengthtest{\Mylen > \labelwidth}}%
                {\parbox[b]{\labelwidth}%     term > labelwidth 
                        {\makebox[0pt][l]{\emph{#1}}\\}}%
                {\emph{#1}}%                 term < labelwidth 
        \hfil\relax}

\newenvironment{Pentry}        
        {%
                \begin{entry}}
        {\end{entry}}

\newcommand{\Mentrylabel}[1]%
        {\raisebox{0pt}[1ex][0pt]{\makebox[\labelwidth][l]%
                {\parbox[t]{\labelwidth}{\hspace{0pt}\emph{#1:}}}}}
        {\begin{entry}}%
        {\end{entry}}
\begin{document}
\title[\runningheadstring]{\titlestring\footnote{Maths Subject Classification \msc}}
\author{\name}
\email{sacha@kurims.kyoto-u.ac.jp}
\address{Research Institute for Mathematical Sciences\\ Ky\=oto University\\ 
  Oiwake-ch\protect\=o\\ Kitashirakawa\\ Saky\protect\=o-ku\\ Ky\=oto
  606-01\\ Japan}

\date{2 October 1996} 

\begin{abstract}

Let~$\Gamma$ be a finite group acting linearly on~$\C^n$, freely
outside the origin, and let $N$ be the number of conjugacy classes of
$\Gamma$ minus one.

In~\cite{sacha:thesis,sacha:ale} a generalisation of Kronheimer's
construction \cite{kron:ale} of moduli of Hermitian-Yang-Mills bundles
with certain invariance properties was given.  This produced varieties
$X_\zeta$ (parameterised by $\zeta\in\Q^N$) which are partial
resolutions of $\C^n/\Gamma$.

In this article, it is shown that $X_\zeta$ can be described as moduli
spaces of representations of the McKay quiver associated to the action
of $\ga$, subject to certain natural commutation relations.
  
This allows a complete description of these varieties in the case when
$\ga$ is abelian.  They are shown to be toric varieties corresponding
to convex polyhedra which are the solution sets for a generalisation
of the transportation problem on the McKay quiver.
  
The generalised transportation problem is solved for a general quiver
to give a description of the extreme points, faces, and tangent cones
to the solution polyhedra in terms of certain distinguished trees in
the underlying graph to the quiver.  Applied to the case of the McKay
quiver, this gives an explicit computational procedure for calculating
$X_\zeta $, its Euler number, and giving a complete list of the
singularities which can occur for all $\zeta$.  The
$\zeta$-parameter-space $\Q^N$ is thus partitioned into a finite
disjoint union of cones inside which the biregular type of $X_\zeta$
remains constant.  Passing from one cone to the other correponds to a
birational transformation.
  
The example $\C^3/\Z_5$ (weights $1,2,3$) is worked out in detail:
there are two types of $\zeta$-cones: ones for which $X_\zeta $ is a
smooth resolution, and others where it has a singularity isomorphic to
a cone over a quadric in $\C^4$.  This gives some evidence for the
conjecture expressed in~\cite{sacha:ale} acording to which the
singularities of $X_\zeta$ are at most quadratic for a generic
$\zeta$.  Computer calculations also show that the cases where
$\Gamma\subset SU(3)$, and $|\Gamma|\leq 11$ yield crepant $X_\zeta$.
For generic $\zeta$, the Euler number of $X_\zeta$ is also equal to
$|\Gamma|$, and the $X_\zeta$ give smooth crepant resolutions fo the
singularity.
  
A further example of a picture of the polyhedron corresponding to
$X_\zeta$ is drawn for the singularity $\frac{1}{11}(1,4,6)$.

\end{abstract}
\maketitle
\tableofcontents

\setcounter{section}{-1}
\section{Introduction}
\label{sec:intro}

This paper is concerned with affine \emph{orbifold singularities},
namely with singularities of the type $X=\C^n/\Gamma$ for $\ga$ a
finite group acting linearly on $\C^n$.

In~\cite{sacha:thesis,sacha:ale} a method using moduli of invariant
Hermitian-Yang-Mills bundles was given for constructing partial
resolutions of $X$ carrying natural asymptotically locally Euclidean
(ALE) metrics.  In this article, the same construction is described in
terms of representation moduli of quivers.  This allows a complete
description of the moduli to be given for the case of abelian groups
using the language of toric varieties.  It turns out that the convex
polyhedra describing the moduli are the solutions of a generalisation
of a well-known network optimization problem (the transportation
problem) on the McKay quiver, in which the quiver plays the role of a
network where commodities are transported, and the parameter $\zeta$
specifies the supplies and demands at each vertex.

\subsection{Background}
\label{sec:intro:back}

The resolution of the Kleinian singularities $\C^2/\Gamma$ for $\ga$ a
finite subgroup of $\SL(2)$ is a classical subject; their minimal
resolution $\widetilde X$ was first constructed by Du Val, and
Brieskorn~\cite{briesk} showed that the components of the exceptional
divisors form graphs which are dual to the homogeneous Dynkin diagrams
for the Lie algebras $A,D,E$.  In 1980, McKay~\cite{mckay:graphs}
remarked that this establishes a correspondence between the extended
homogeneous Dynkin diagrams $\clos A,\clos D,\clos E$ and the
irreducible representations of $\ga$.  More recently~\cite{reid_ito}
another correspondence was constructed between conjugacy classes ``of
weight 1'' and crepant divisors for any quotient singularity generated
by a finite subgroup of $\SL(n)$.

In the case where $\Gamma\subset\SL(3)$, Dixon, Harvey, Vafa and
Witten proposed a definition of the orbifold Euler
characteristic~\cite{dhvw:i,dhvw:ii} and conjectured the existence of
smooth resolutions with trivial canonical bundle (\ie which are
crepant) and whose Euler number is given by the DHVW orbifold Euler
number.  The final cases in the proof of this conjecture were only
recently completed
\cite{mar_ols_per,roan:mirror_cy,mark:res_168,roan:res_a5,ito:trihedral,roan:crepant}.

For the case $\Gamma\subset\SL(4)$, the author has obtained some
interesting analogous results if one considers terminalisations rather
than resolutions~\cite{sacha:sl4}.

A promising link between the representation theory of $\Gamma$ and the
construction of resolutions was established in 1986;
Kronheimer~\cite{kron:thesis,kron:ale} constructed a family of
hyper-\kah\ quotients $X_\zeta $ by looking at a dimensional reduction
of the Hermitian Yang-Mills (HYM) equations on a $\ga$-equivariant
bundle over $\C^2$. He used the representation theory of $\ga$ and
McKay's observation to prove that, for generic $\zeta$, these
quotients are isomorphic to the minimal resolution $\widetilde X$.
The author's thesis~\cite{sacha:thesis} (from which the present paper
and its companion~\cite{sacha:ale} are in large part extracted) grew
out of the ambition to generalise the results as far as possible to
dimensions higher than two.

Despite the fact that there is little hope for the results
for $\GL(n)$ and general $n$ to be as strong as those for $\SL(2)$
groups, some of the nice properties of the $\SL(2)$ case generalise to
the $\GL(n)$ case.\footnote{Some interesting new aspects are also
  present for the $\SL(3)$ case~\cite{sacha:thesis}.}

In~\cite{sacha:thesis} several descriptions were given of the
generalised construction.  Probably the most concise is the one given
in~\cite{sacha:ale}: construct moduli spaces $X_\zeta$ of instantons
on the trivial bundle $\C^n\times R\to\C^n$.  Here $R$ denotes the
regular representation space for the group $\Gamma$, and $\zeta$ is a
linearisation of the bundle action.  The instantons are required to
satisfy Hermitian-Yang-Mills-type equations, as well as additional
$\Gamma$-equivariance and translation-invariance properties.

In the present paper, an entirely different view point is adopted,
namely that of considering $X_\zeta$ are representation moduli of the
McKay quiver.  More precisely, one considers representations
$\Rep_{\cQ,\cK}(\bdR)$ of the McKay quiver $\cQ$ into the multiplicity
space $\bdR$ of the regular representation of $\Gamma$, subject to
certain commutation relations $\cK$.  A reductive group $\PGL(\bdR)$ acts on the above space, and one obtains GIT quotients  $X_\zeta$
which depend on a rational parameter $\zeta$.

\subsection{Main Results}
\label{sec:intro:main}

This paper assumes that $\Gamma\subset\GL(n)$ acts on $\C^n$ freely
outside the origin for any $n\geq 2$, which means that $X=\C^n/\Gamma$
has an isolated singularity.\footnote{This is for the purpose of
  simplicity --- the method would seem to be applicable to the general
  case with some modifications.} The main results are as follows.

Firstly, $X_\zeta$ are identified with representation moduli of the
McKay quiver.  This is fairly straight-forward and involves mostly
translating quiver concepts over to the language used
in~\cite{sacha:ale}.  There are two further sets of results.

\subsubsection{Toric description of $X_\zeta$}

The first is a toric description of the representation variety
$\Rep_{\cQ,\cK}(\bdR)$ and its moduli $X_\zeta$.  The convex
polyhedron $C_\zeta$ corresponding to $X_\zeta$ is identified with the
projection to $\R^n$ of the solution polyhedron for the
transportation problem on the McKay quiver.  The transportation
problem is a well-known linear network optimization
problem~\cite{kenn_helg,gond_mino:graphs}: given an assignment of real
numbers to the vertices of a quiver\footnote{The term
  \emph{network} is usually used in this context instead of quiver.}  $\cQ$
(thought of as representing the demand and supply of certain
commodities), the aim is to find an assignment of non-negative real
numbers to the arrows (a \emph{flow} on $\cQ$), in such a way that the
demands and supplies at each vertex are satisfied (\ie the equation
$\partial f=\zeta$ holds).  For any given assignment of weights
$\zeta$, the solution set is a convex polyhedron inside $\R^{\cQ_1}$
denoted by $F_\zeta$.

The notation in use is as follows:

\begin{notation}
    Let $\ga$ be the cyclic group of order $r$ acting freely outside the
  origin in $Q$ with weights $w_1,\dots,w_n\in\Z_r=\Z/r\Z$ \ie via
  \map{\rho}{\mu_r\subset\CX}{\CX^n}{\lambda}{\begin{pmatrix}
\lambda^{w^1} \\
 & \lambda^{w^2}  \\
 &\\
 &  &  & \lambda^{w^n}
\end{pmatrix}.}
  The McKay
  quiver $\cQ$ of $\Gamma$ has vertices 
  $\cQ_0=\Z_r$, and arrows $a_v^i:=v\to v-w_i \pmod r$ for each $v\in\Z_r$ and
  $i=1,\dots,n$.
  
  Let $\pi\colon\cQ_1\to\{1,\dots,n\}$ be the map defined by
  $\pi(a_v^i):=i$ and let $\pi\colon\R^{\cQ_1}\to\R^n$ be the induced
  linear map defined by mapping the standard bases.  (Here and
  elsewhere we use the exponential notation
  $\R^{\cQ_1}:=\Map(\cQ_1,\R)$.  The standard basis of $\R^{\cQ_1}$ is
  given by the ``indicator'' functions $\chi_a$ defined by
  $\chi_a(a')=1$ if $a=a'$ and $\chi_a(a')=0$ otherwise).  Let $\Pi$
  be the sub-lattice of $\Z^n$ of index $r$ defined by
  $$\Pi := \ker\hat\rho=\{ x\in \Z^n | \sum x_iw_i\equiv 0\pmod r\}.$$
  For each
  element $\zeta\in\Z^{\cQ_0}$ such that $\sum_{v\in\cQ_0}
  \zeta(v)=0$, consider the corresponding solution polyhedron
  $F_\zeta$ to the transportation problem on $\cQ$, and let
  $C_\zeta=\pi F_\zeta$ be its projection to $\R^n$.
\end{notation}

Adopting the above notation, one can state the first main theorem.

\begin{nonumberthm}[\cf Thms.\ \ref{thm:toric-d-moduli} and \ref{thm:slices}] 
  The moduli $X_\zeta$ are isomorphic to the toric varieties
  $$X_\zeta\cong\T^{\Pi,C_\zeta},$$
  where $C_\zeta$ correspond to
  the projection to $\R^n$ of the solution polyhedra to the
  transportation problem on $\cQ$.
\end{nonumberthm}

\subsubsection{Generalised Transportation Problem}
\label{sec:intro:ab:trans}

The problem of determining $C_\zeta$ can be considered as a
generalised transportation problem.  The solution to the classical
transportation problem is well known:

\begin{nonumberthm}[Classical Transportation Problem (\cf Theorem~\ref{thm:classical})]
  The extreme points of $F_\zeta$ are precisely those flows in
  $F_\zeta$ whose supports\footnote{The support of a flow is the set
    of arrows on which it is non-zero.} are \emph{trees}, \ie contain
  no cycles.
\end{nonumberthm}

The basic idea behind the proof of this theorem is to associate to
each cycle $c$ a corresponding flow $\tilde\chi_c$ as follows.  A
\emph{cycle} in $\cQ$ is a sequence $c_1,\dots,c_m$ of arrows such
that they form a cycle when their orientation is disregarded.  The
disjoint union of the arrows $c_i$ whose direction agrees (\resp disagrees)
with the ordering $c_1,\dots,c_k$ is called the \emph{positive\/}
(\resp \emph{negative\/}) part of $c$ and is denoted it by $c^+$
(\resp $c^-$).

For each arrow $a\in\cQ_1$, let $\chi_a$ denote the flow which takes
the value $1$ on the arrow $a$ and zero elsewhere.  To each cycle $c$,
one can define the \emph{basic flow} associated to $c$ by
$$\tilde\chi_c:=\sum_{c_i\in c^+}\chi_{c_i} - \sum_{c_i\in c^-}\chi_{c_i}\in\Z^n.$$

Note that for any cycle $c$, $\partial \tilde\chi_c=0$, \ie the
associated flow does not contribute anything at any vertex.  If $f\in
F_\zeta$ contains $c$ in its support, then one can add
$\pm\epsilon\tilde\chi_c$ to $f$ for some small $\epsilon$ and still
remain in $F_\zeta$.  Hence $f$ cannot be an extreme point.
Vice-versa, if $f$ is not extreme, then one can reconstruct a cycle in
its support.

The generalisation of this theorem to describe the projection
$C_\zeta=\pi(F_\zeta)$ is quite straight-forward.  One begins by
defining the \emph{type} of any given cycle as the element
$\pi(\tilde\chi_c)\in\Z^n$.  Cycles of type $0\in\Z^n$ play a special
role, much as ordinary cycles do in the classical case.

\begin{dfn}
  The \emph{closure} of $S\subset\cQ_1$ is defined to be the
  smallest over-set $\clos{S}\supseteq S$ such that $$c^-\subseteq
  \clos{S}\iff c^+\subseteq\clos{S},$$
  for all cycles $c$ of type
  zero.  Two configurations $S,S'$ will be called \emph{equivalent}
  (written $S\sim S'$) if $\clos{S}=\clos{S'}$.
\end{dfn}
  
Th same methods as in the classical case allows one to prove the
following generalisation.

\begin{nonumberthm}[Extreme Points of $C_\zeta$]
  The extreme points of $C_\zeta$ are the images under $\pi$ of
  precisely those flows in $F_\zeta$ whose support contains no cycles
  of non-zero type.
\end{nonumberthm}

In fact, by taking into account all cycles in the support of a given
flow, it is possible to determine the dimension of the face which
$\pi(f)$ belongs to.  Some additional terminology will be convenient
to state the results.

Let $\cC$ denote the set of all \emph{configurations}, \ie the set of
non-empty subsets of $\cQ_1$. For a configuration $S\subset\cQ_1$, let
$F_0(S)$ be the cone generated by the flows $\tilde\chi_c$ for the
cycles $c$ whose negative part is included in $S$ and let $Z_0(S)$ be
its maximal vector subspace (\ie the subspace generated by the flows
$\tilde\chi_c$ for the cycles $c\subseteq S$).  The \emph{rank} of $S$
is defined to be the dimension of $\pi Z_0({\clos S})$.  The set of
configurations (\resp trees) of rank $k$ is denoted $\cC^k$ (\resp
$\cT^k$).  Also, write $\cC_\zeta$ (\resp $\cT_\zeta$) for the subset of
configurations (\resp trees) which are \emph{admissible for
  $\zeta$},\/ namely configurations (\resp trees) $S\in\cC$ which
arise as the support of some element in $F_\zeta$.

\begin{nonumberthm}[\cf Theorem \ref{thm:faces}]
  \label{thm:intro:faces}
  For all $\zeta$, the map
  \map{\text{Face}_\zeta}{\cC_\zeta}{\text{Faces of }\pi
    F_\zeta}{S}{\pi F_\zeta \cap(\pi f+\pi Z_0(\clos{S})))} is
  independent of the choice of $f\in F_\zeta\cap\supp^{-1}(S)$, and
  induces a bijection
  $$\text{Face}_\zeta\colon \cC^k_\zeta/\!\!\sim
  \xrightarrow{\;\cong\;} \text{$k$-faces of }\pi F_\zeta.$$
  Furthermore, for all $[S]\in\cC^k_\zeta/\!\!\sim$,
  $$T_{\text{Face}_\zeta(S)}\pi F_\zeta =\pi F_0(\clos{S}),$$
  where
  the left-hand side denotes the tangent cone to the polyhedron $\pi
  F_\zeta$ at the face $\text{Face}_\zeta(S)$.

  In other words, $\pi F_0(\clos{S})$ gives the tangent cone
  corresponding to the configuration $S$ (which is independent of the
  value of $\zeta$) and $\text{Face}_\zeta(S)$ gives the corresponding
  face of $C_\zeta$ (whose direction is also independent of $\zeta$).
\end{nonumberthm}

This theorem has a number of important
corollaries.

\begin{nonumbercor}[\cf Cor.\ \ref{cor:iso_C_zeta}]
  \label{cor:intro:iso_C_zeta}
  If $\zeta$ and $\zeta'$ have the same admissible configurations
  ($\cC_\zeta=\cC_{\zeta'}$) or even just the same admissible trees
  ($\cT_\zeta=\cT_{\zeta'}$) then the corresponding polyhedra
  $C_{\zeta}$ and $C_{\zeta'}$ are geometrically isomorphic.  Two
  polyhedra are said to be \emph{geometrically isomorphic} if they are
  combinatorially isomorphic and their tangent cones at the
  corresponding faces are identical.  In particular, their associated
  fans and toric varieties must be identical.
\end{nonumbercor}
\begin{nonumbercor}[\cf Cor. \ref{cor:number_extreme_pts}]
  \label{cor:intro:number_extreme_pts}
  Let $\cT\subset\cC$ denote the configurations which are {\em
    trees,\/}  and let $\ext
  C_\zeta$ denote the extreme points of $C_\zeta$.  Then
  $$\card{\ext C_\zeta}=\card{\cC^0_\zeta/\!\!\sim}=\card{\cT^0_\zeta/\!\!\sim}.$$
\end{nonumbercor}
\begin{rmk}
  Since any tree admits a unique flow such that $\partial f = \zeta$,
  it is very easy to determine $\cT_\zeta$.  Furthermore, if $\zeta$
  is generic, then $\cC_\zeta\subset\cC_{\text{span}}$, the subset of
  \emph{spanning configurations},\/ \ie configurations whose arrows
  join any two vertices of the quiver (not necessarily in an oriented
  way).  In this case, $T\in\cT^0_\zeta$ if and only if the tree $T$
  admits an assignment of elements of $\Z^n$ to the vertices
  which  satisfies particular properties (a so-called \emph{$n$-weighting}
  in the terminology of Section~\ref{sec:2:examples:weight:fixed}).  This
  condition can again be checked by a very simple algorithm.  
  
  To sum up: the above corollary allows one to determine the extreme
  points of $C_\zeta$ for any $\zeta$ very easily.
\end{rmk}
\begin{nonumbercor}[\cf Cor.\ \ref{cor:fan} and Lemma \ref{lemma:generic_flow}]
  \label{cor:intro:fan}
  The extreme points of the polyhedron $\pi F_\zeta$ correspond to the
  trees $T$ in $\cT^0_\zeta$ and the tangent cone to $\pi F_\zeta$ at
  the point corresponding to $T$ is $\pi F_0(\TC)=\pi F_0(T)$. Thus
  the fan associated to the polyhedron $\pi F_\zeta $ is given by the
  dual cones $\pi F_0(T)^\vee$ for the trees $T\in\cT^0_\zeta$ and all
  their faces. 
%These cones correspond (with respect to the lattice
%  $\Pi$) to the possible singularities of $C_\zeta$ for generic
%  $\zeta$.
\end{nonumbercor}

The theorem and corollaries above allow a
complete understanding of the extreme points, faces, and tangent cones
to the solution polyhedra of a generalised transportation problem.
Applied to the case of the McKay quiver, they allow one to determine
the Euler number and the singularities of $X_\zeta$ as well as giving
a complete list of the singularities which can occur for all $\zeta$.

Several interesting questions regarding these moduli nevertheless
remain, such as whether they produce smooth (\resp terminal)
resolutions in the $\SL(3)$ (\resp $\SL(4)$) case.  A conjecture is
given in the companion paper~\cite{sacha:ale} and several examples are
computed in Section~\ref{sec:2:examples} of this paper.

Depending on the point of view, isomorphic objects are denoted by
different notations.  For ease of reference, Table~\ref{tab:notation}
gives a correspondence table between the point of view
in~\cite{sacha:ale} and in Parts~\ref{sec:1} and~\ref{sec:2} of this
paper.

\renewcommand{\arraystretch}{1.4}
\begin{table}[htbp]
  \begin{center}
    \leavevmode
\setlength{\fboxrule}{0pt}
\begin{tabular}{|l|l|l|}
  \hline
  Moduli of Bundles & Moduli of Representations&
  Toric Varieties\\
   \cite{sacha:ale} & Part 1 of this paper & Part 2 of this paper\\
  \hline
  $M^\Gamma  $&$ \Rep_\cQ(\bdR) $&$ \T^{\Lambda^1_\R,\Lambda^1_{\R+}}=\C^{\cQ_1}$\\
  \hline
  $\cN^\Gamma $&$ \Rep_{\cQ,\cK}(\bdR)  $&$ \T^{\Lambda,C}$\\
  \hline
  $G^\Gamma=\PGL^\Gamma(R) $&$ \PGL(\bdR) $&$ \T^{\Lambda^{0,0}}=\T^{0,0}$\\
  \hline
  $X_\zeta $&$ \cM_\zeta $& $\T^{\Pi,C_\zeta}$\\
  \hline
\end{tabular}
\vspace{1em}
\caption{Correspondence between the notations in use in this paper and in the 
  companion paper~\cite{sacha:ale}.}
    \label{tab:notation}
  \end{center}
\end{table}

\subsection{Methods and Outline}
\label{sec:intro:outline}

The methods used to prove the main results are notationally cumbersome
due to the quiver notation, but on the whole straightforward.  A
familiarity with toric geometry will make reading easier, although all
the required facts and notation are explained when needed.

The paper is divided into two parts; the first dealing with general
groups, and the second specializing to the case of abelian ones.  Each
part begins with a section summarizing the notation in use.

\subsubsection{Part \ref{sec:1}}

Section~\ref{sec:1:notation} summarises the notation.

Section~\ref{sec:1:quivers} outlines the basic facts and notation
concerning quivers, their representations and their moduli.  The
necessary details regarding \git (GIT) quotients are also given.

Section~\ref{sec:1:mckay} explains the case of the McKay
quiver associated to the representation of a finite group.  The McKay
correspondence is briefly mentioned, and the ``canonical'' commutation 
relations $\cK$ are defined.

\subsubsection{Part \ref{sec:2}}

Section~\ref{sec:2:notation} summarises the notation.

Section~\ref{sec:2:abel} specialises further the discussion from
Part~\ref{sec:1} to the case of abelian groups.  The McKay quiver and
commutation relations for these can be described quite simply and some
specialised notation is introduced for this purpose.  The
representation variety and its moduli are proved to be toric varieties
corresponding to an $(n+\card{\Gamma}-1)$-dimensional convex cone $C$
and to its $n$-dimensional convex polyhedral slices $C_\zeta$
respectively.

From Section~\ref{sec:2:flow} onwards, the focus is on the
transportation problem; its solution polyhedra coincide with $C_\zeta$
in the case of the McKay quiver.

Section~\ref{sec:2:flow} explains the ordinary transportation problem 
on a network.  

Section~\ref{sec:2:gtp} explains the generalised  transportation
problem and states the theorems describing its solution polyhedra.

An example of an application to the McKay quiver and several important
corollaries are given in Section~\ref{sec:2:mckay}.

The proofs of the theorems are given in
Section~\ref{sec:2:flow:proofs}, except for the proof of some
technical lemmas regarding flows which are left until
Section~\ref{sec:2:exact}.

Section~\ref{sec:2:sing} contains a discussion of the singularities of
$X_\zeta$.

The paper concludes in Section~\ref{sec:2:examples} with some
practical examples and computations.

\subsection{Acknowledgments}
\label{sec:intro:ack}

The present paper and its companion paper~\cite{sacha:ale} consist
mostly\footnote{Minor portions have been rewritten to include
  references to advances in the field made since then
  (notably~\cite{ito:trihedral,roan:crepant,reid_ito}).} of excerpts
of my D.Phil.\ thesis~\cite{sacha:thesis}, and I wish to acknowledge
the University of Oxford and Wolfson College for their hospitality
during its preparation.  I am grateful to the Rhodes Trust for
financial support during my first three years, and to Wolfson College
for a loan in my final year.  The conversion from thesis to article
format was done while I was a Research Assistant in RIMS, Kyoto.

I also take the opportunity to thank my supervisors Peter Kronheimer
and Sir~Michael Atiyah who provided me with constant advice,
encouragement and support and whose mathematical insight has been an
inspiration.  I also wish to thank William Crawley-Boevey, Michel
Brion, Gavin Brown, Jack Evans, Partha Guha, Katrina Hicks, Frances
Kirwan, Alistair Mees, Alvise Munari, Martyn Quick, David Reed, Miles
Reid, Michael Thaddeus, and, last but not least, my parents and
family.

\section{Summary of Notation for Part \ref{sec:1}}
\label{sec:1:notation}

\subsection{General}

\begin{Pentry}
\item[$k$] Algebraically closed field.
\item[$\Z_+$] Non-negative integers.
\item[$\R_+$] Non-negative reals.
\end{Pentry}

\subsection{Quivers, Representations}

\begin{Pentry}
\item[$Q$] $n$-dimensional complex vector-space.
\item[$\Gamma$] Finite sub-group of $\SL(Q)$
\item[$\cQ$]     Generic quiver $\cQ=(\cQ_0,\cQ_1)$.\\
  From Section~\ref{sec:1:mckay} onwards, $\cQ=\cQ_{\Gamma,\cQ}$, the
  McKay quiver associated to $(\Gamma,Q)$.
\item[$\cQ_0$] Vertices of $\cQ$; in the case of the McKay quiver,
  elements of $\cQ_0$ index the irreducible representations $R_v$ of
  $\Gamma$.
\item[$\cQ_1$] Arrows of $\cQ$.
\item[$R_v$] Irreducible representation of $\Gamma$ corresponding to
  vertex $v$ of the McKay quiver.
\item[$R$] Regular representation of $\Gamma$; $R=\oplus_{v\in\cQ_0}
  {\bdR_v}\otimes R_v$.
\item[$\bdR_v$] The trivial $\Gamma$-module giving the multiplicity of 
  $R_v$ in $R$.
\item[$\bdR$] Multiplicity space for $R$;  $\bdR=\oplus_{v\in\cQ_0}{\bdR_v}$
\item[$\bdV$] Generic representation space of $\cQ$; $\bdV=\oplus_{v\in\cQ_0}\bdV_v$.
\item[$R_\bdV$] The $\Gamma$-module $\oplus_{v\in\cQ_0} \bdV_v\otimes R_v$
  corresponding to $\bdV$.
\item[$\alpha$] $\Gamma$-invariant element of $(Q\otimes\End R)$.
\item[$q_i$] Basis of $Q$ ($i=1,\dots,n$).
\item[$\alpha_i$] Component of $\alpha$ with respect to $\{q_i\}$:
  $\alpha=\sum_{i=1}^n q_i\otimes\alpha_i$.
\item[$\tilde\alpha$] Corresponding representation of the McKay quiver
  into $\bdR$.
\item[$a$] Arrow of $\cQ$; also written $t(a)\to h(a)$, where $t(a)$
  denotes the tail and $h(a)$ the head of $a$.
\item[$\tilde\alpha_a$] Value of $\tilde\alpha$ on the arrow $a$;
  $\tilde\alpha_a\colon \bdR_{t(a)}\to \bdR_{h(a)}$.
\item[$\cK$] Relations on $\cQ$.
\item[$\Rep_{\cQ,\cK}(\bdR)$] Space of all representations of $\cQ$ into
  $\bdR$ satisfying the relations $\cK$.
\item[$\PGL(\bdR)$] The isomorphism group for representations of $\cQ$
  into $\bdR$
  $$\PGL(\bdR):=\times_{v\in\cQ_0}\GL({\bdR_v})/\CX.$$
\item[$\frk$] The Lie algebra of $\PU(\bdR)$ (a real form of $\PGL(\bdR)$).
\item[$\cM_{\cQ,\cK,\zeta}$] Representation moduli of $(\cQ,\cK)$
  $$\cM_{\cQ,\cK,\zeta}:= \Rep_{\cQ,\cK}(\bdR)\gitquot{\chi}\PGL(\bdR).$$
  The linearisation $\chi$ is related to $\zeta$ by
  $\zeta=d\chi(1)_{|\frk}$.
\item[$X_{\zeta}$] The moduli above in the case when $\cQ$ is the
  McKay quiver, and $\cK$ are the commutation relations defined in
  Section~\ref{sec:mckay:mckay:ccr:relns}.
\item[$\rho_{\zeta}$] The natural projective morphism $X_{\zeta}\to
  X_{0}$ (a partial resolution).
\end{Pentry}

\part{General Groups}
\label{sec:1}

\section{Quivers, Relations and Representation Moduli}
\label{sec:1:quivers}

\subsection{Quivers}
\label{sec:1:quivers:quivers}
%\label{sec:mckay:quiver:reps}

A \emph{quiver} is an oriented graph, possibly with multiple arrows
between the same vertices and with loops (arrows which begin and end
at the same vertex).  Formally, a quiver $\cQ$ consists of a pair of
finite sets $\cQ=(\cQ_0,\cQ_1)$ with two maps $\cQ_1
\overset{h}{\underset{t}\rightrightarrows} \cQ_0$; the elements of
$\cQ_0$ are called \emph{vertices} and those of $\cQ_1$ are called
\emph{arrows}.\/ The elements $t(a)$ and $h(a)$ are called the
\emph{tail\/} and \emph{head\/} of the arrow $a\in\cQ_1$ respectively.
The arrow $a$ is also sometimes denoted $t(a)\to h(a)$.

\subsection{Representations}
\label{sec:1:quivers:reps}

Let $k$ be an algebraically closed field.

A \emph{representation of a quiver} is a realization of its diagram of
vertices and arrows in some category: it corresponds to replacing the
vertices by objects and the arrows by morphisms between the objects.
From now on, only the category of $k$-vector-spaces is considered; a
representation $\bdV$ of a quiver $\cQ$ is thus taken to be a collection
of finite dimensional vector-spaces $\bdV_v$, indexed by the vertices
$v\in \cQ_0$, and of linear maps $\bdV_{v\to v'}\colon \bdV_{v}\to \bdV_{v'},$
indexed by the arrows $v\to v' \,\in\cQ_1$.

The set of
all representations of $\cQ$ into a fixed $\cQ_0$-graded vector-space
$\bdV=\oplus_{v\in \cQ_0} \bdV_v$ forms a vector-space, denoted $\Rep_\cQ
\bdV$. 
%% An alternative notation often encountered is to set $\dim \bdV
%%:=(\dim \bdV_v)_{v\in\cQ_0}$ be the \emph{dimension
%%  vector}  and to write $\Rep_\cQ(\dim \bdV)$ for
%%$\Rep_\cQ(\oplus_{v\in\cQ_0}\C^{\dim \bdV_v})\cong\Rep_\cQ(\bdV)$. 

\begin{example}
\label{ex:endo}
Representations of the quiver $\cQ$ with one vertex and one loop are just 
 endomorphisms of a vector-space. 
\end{example}

\begin{example}
  A \emph{primitive} representation of
  $\cQ$ is an element of
  $$\Rep_{\cQ}(k^{\cQ_0}),$$ where $k^{\cQ_0}$ denotes
  the free $k$-vector-space on the vertices.  A primitive representation
  thus corresponds to an assignment of an element of $k$ to each arrow
  in $\cQ$, \ie to an element of the vector-space $k^{\cQ_1}$.
\label{ex:primitive}
\end{example}

\subsection{Relations}
\label{sec:1:quivers:relations}
%\label{sec:mckay:quiver:moduli}

If  the morphisms $\bdV_{v\to v'}$ are required to satisfy   relations
between them, then one talks about a \emph{representation of a quiver with
  relations}.

More formally, a \emph{relation} is defined as a formal sum of paths in
a quiver.  A (non-trivial) \emph{path} $p$ is a sequence
$a_1\dots a_n$ of arrows which \emph{compose}, \ie such that
$t(a_i)=h(a_{i+1})$ for $1\leq i < n$:
$$\stackrel{h(p)}\bullet \xleftarrow[{a_1}]{} \bullet \xleftarrow[{a_2}]{} \dots
\xleftarrow[{a_n}]{}\stackrel{t(p)}\bullet.$$ The vertex $h(a_1)$ ($t(a_n)$) is 
called the
\emph{head} (\emph{tail}) of the path $p$ and denote it by $t(p)$ ($h(p)$). If $\bdV$ is a representation
of $\cQ$, then there is an induced morphism
$$\bdV(p)=\bdV_{a_1}\bdV_{a_{2}}\dots \bdV_{a_n}\colon \bdV_{t(p)} \to \bdV_{h(p)}$$
corresponding to each path. The \emph{trivial path} $e_v$ 
consists of a single vertex $v$ and no arrows; $\bdV_{e_v}$ is of
course the identity endomorphism of $\bdV_v$.

A \emph{representation \/} is said to \emph{satisfy}
$r=\sum\lambda_ip^i$ if $\sum\lambda_i\bdV(p^i)=0$.  If $\cK$ denotes
a set of relations, then the set of representations of $\cQ$ into
$\bdV$ satisfying the relations $\cK$ will be denoted
$\Rep_{\cQ,\cK}(\bdV)$.  It is an affine variety inside
$\Rep_\cQ(\bdV)$.

\begin{rmk}
  There are in general many classes of relations one can consider, but
  in this paper, only \emph{commutation relations}
  will be considered, namely relations generated by differences of two
  paths, both paths having the same tail and head.
\end{rmk}

\subsubsection{Morphisms of Representations}
\label{sec:1:quivers:morph}
%\label{sec:mckay:quiver:moduli:morph}
\aftersub

Let $\bdV,\bdV'$ be two representations of the same quiver $\cQ$ (possibly
with relations).  A \emph{morphism} of
representations $\theta\colon \bdV\to \bdV'$ is given by morphisms
$\theta_v\colon \bdV_v\to \bdV'_v$ for each vertex $v$, which satisfy
$\bdV'_{a}\theta_{t(a)}=\theta_{h(a)}\bdV_a$ for each arrow $a$.  If
$\bdV_v=\bdV_v'$ and the linear maps $\theta_v$ are all isomorphisms, then
$\theta$ is called an \emph{isomorphism of representations}.

\subsubsection{Representation Moduli}
\label{sec:1:quivers:moduli}
\aftersub

Given a quiver $\cQ$, some relations $\cK$ and a representation space
$\bdV$, one is naturally interested in the \emph{moduli space} of representations.  Loosely speaking, this is the
set of isomorphism classes of representations of $\cQ$ into $\bdV$.  More
precisely, the group 
$$\GL(\bdV):=\times_{v\in\cQ_0}\GL(\bdV_v)$$
acts on
$\Rep_{\cQ,\cK}(\bdV)$ and its orbits consist of equivalence classes
of representations.  The scalar subgroup $\CX\subset\GL(\bdV)$ acts
trivially, and one is left with a free action of the quotient
$\PGL(\bdV):=\GL(\bdV)/\CX$.

Since  $\PGL(\bdV)$ is not compact, one must resort to \git\ (GIT) in order
to obtain quotients which are well-defined as quasi-projective
varieties.  The ``canonical'' quotient is the affine GIT quotient
\begin{equation}
  \label{eq:zero-quot}
  \cM_{\cQ,\cK,0}:= \Rep_{\cQ,\cK}(\bdV)\gitquot{} \PGL(\bdV).
\end{equation}
There are other possible quotients however and their existence is in
fact the basis of this paper.  In general, given a complex
affine variety $X$ acted upon by a reductive group $G$, and given a
character $\chi\colon G\to\CX$, one defines
$$X\gitquot\chi G := \Proj \bigoplus_{r\in\N} \cO_X(rL_\chi)^G,$$
where $L_\chi$ is the trivial bundle over $X$, on which $G$ acts via
$\chi$. The reader is refered to~\cite{sacha:ale,sacha:thesis} for a
detailed treatment.  

When the complex reductive group $G$ has a real form $K$, rather than
specifying the character $\chi$ of $G$ one can specify instead its
derivative at the identity $\Lie G\to \R$, or even the restriction of
the later to the Lie algebra of the real form, giving an element of
$\zeta=d\chi(1)_{|\Lie K}\in(\Lie K)^*$.  

For the moduli of representations of $\cQ$ into $\bdV$, a real form of
$\PGL(\bdV)$ is given by $\PU(\bdV):=\times_{v\in\cQ_0}\U(\bdV_v)/\U(1)$, 
and the dual of its Lie algebra is the subspace 
$${\frk}:=\{\zeta\in \times_{v\in\cQ_0} \End_\R(\bdV_v)^*|
\sum_{v\in\cQ_0}\trace \zeta_v=0\},$$
where $\zeta_v$ denotes the restriction
of $\zeta$ to $\End_\R(\bdV_v)$.
For $\zeta$ an integral element of $\frk$, define
\begin{equation}
  \label{eq:zeta-quot}
  \cM_{\cQ,\cK,\zeta_\chi} := \Rep_{\cQ,\cK}(\bdV)\gitquot\chi \PGL(\bdV),\quad\text{where }\zeta_\chi=d\chi(1)_{|\Lie K}
\end{equation}
 This is seen to coincide with Definition~\ref{eq:zero-quot} in the case 
$\zeta_\chi=0$ (\ie $\chi=1$).

\begin{rmk}
  In the above, ``$\zeta$ integral'' means integral with respect to
  the natural lattice in $\frk$ (\ie the kernel of the exponential
  map).  Strictly speaking, one is not restricted to integral values,
  any rational value will do, provided one makes sense of ``fractional
  linearisations'' in the obvious manner.  We do not bother, as
  nothing substantially new is gained from this approach (the moduli
  for $k\zeta$ are isomorphic to those for $\zeta$).
\end{rmk}

\section{The McKay Quiver}
\label{sec:1:mckay}

The McKay quiver $\cQ_{\ga,Q}$ is a quiver which is naturally
associated to a representation $Q$ of a finite group $\ga$. As
explained below, its vertices are the irreducible representations of
$\ga$ and the arrows describe how the tensor product of $Q$ with each
irreducible decomposes into a sum of irreducibles.

It seems natural to expect a relation between $\cQ_{\ga,Q}$ and the
quotient singularity $X=Q/\ga$. In fact, for the case of a finite
subgroup $\ga\subset \SU(2)$ there is a remarkable relation between
$\cQ_{\ga,Q}$ and the minimal resolution
$\widetilde{X}\to X$. McKay remarked~\cite{mckay:graphs} that the
quivers $\cQ_{\ga,\C^2}$ are precisely the extended Dynkin diagrams of
type $\overline A,\overline D$, and $\overline E$. The ordinary Dynkin
diagrams of type $A,D$ and $E$ had previously been shown by
Brieskorn~\cite{briesk} to be the dual graphs to the graphs of rational
curves in the exceptional fibre of $\widetilde{X}\to X$.

Kronheimer~\cite{kron:thesis,kron:crendus,kron:ale} showed that one
can construct the minimal resolution $\widetilde{X}$ by considering
what turn out to be, upon closer inspection, moduli spaces of
representations of $\cQ_{\ga,\C^2}$ into the regular representation
space of $\Gamma$.  The representations are also required to
satisfy some commutation relations, denoted $\cK$.

In this section, the generalisation of these commutation relations to
any McKay quiver is described.

\subsection{The McKay Quiver}
\label{sec:mckay:mckay:dfn}
\aftersub 
Let $\ga$ be a finite group and let $\{R_i, i\in I\}$ be the 
set of irreducible representations of $\ga$. Any 
$\ga$-module $R_\bdV$ decomposes into a sum of irreducibles: 
\begin{equation}
R_\bdV= \bigoplus_{i\in I}\bdV_i\otimes R_i,
\label{eq:irred_decomp}  
\end{equation}
and gives an $I$-graded vector-space (trivial as a $\ga$-module)
$\bdV=\oplus_{i\in I}\bdV_i$, called the \emph{multiplicity space} for $R_\bdV$.
The dimension of $\bdV_i$ is called the \emph{multiplicity} of $R_i$ in
$R_\bdV$ and the vector $\dim \bdV=(\dim \bdV_i)_{i\in I}$ is called the
\emph{dimension vector} of the $\Gamma$-module $R_\bdV$. Conversely, given
any $I$-graded vector-space $\bdV$, one can construct a corresponding $\ga$-module
$R_\bdV$ by formula~\eqref{eq:irred_decomp}.

The \emph{McKay quiver} $\cQ=\cQ_{\ga,Q}$ of a
representation $Q$ is constructed as follows.  The  vertices $\cQ_0=I$
are the irreducible representations of $\ga$, and there are  $a_{ij}$ (possibly
zero) arrows from vertex $i$ to vertex $j$.  The non-negative
integers $a_{ij}$ (which form what is called the \emph{adjacency
  matrix} of $\cQ$) are defined by the following irreducible
decompositions (for each $i\in I$)
\begin{equation}
\label{eq:aijdecomp}
Q\otimes R_i=\bigoplus_j a_{ji}R_j.
\end{equation}
More invariantly, one may write 
\begin{align}
  Q\otimes R_i &= \bigoplus_j A^*_{ji}\otimes R_j\\
\intertext{where}
A_{ji}:&=\Hom_\ga(Q\otimes R_i,R_j),
\end{align}
and think of the arrows from $i$ to $j$ as giving a basis for the
$a_{ij}$-dimensional vector-space $A_{ij}$.

If $\bdV$ is a $\cQ_0$-graded vector-space, then the representations of
$\cQ$ into $\bdV$ correspond to the $\ga$-module endomorphisms $R_\bdV\to
Q\otimes R_\bdV$. In fact, 
\begin{equation}
\begin{split}
   \Hom_\ga(R_\bdV,Q\otimes R_\bdV) 
&= \Hom_\ga(\oplus_i \bdV_i\otimes R_i, Q\otimes (\oplus_j\bdV_j\otimes R_j))\\ 
&= \bigoplus_{i,j} A_{ij}^*  \otimes \Hom(\bdV_i,\bdV_j),
\end{split}
\end{equation}
so given $\alpha\in\Hom_\ga(R_\bdV,Q\otimes R_\bdV)$, one 
can pair it with an arrow $a=i\to j$ (considered as a basis element of
$A_{ij}$) and obtain a map $\bdV_i\to \bdV_j$.  Thus
\begin{equation}
  \label{eq:hom-rep}
  \Hom_\ga(R_\bdV,Q\otimes R_\bdV) = \Rep_\cQ(\bdV).
\end{equation}
The representation of $\cQ$ which corresponds to
$\alpha\in\Hom_\ga(R_\bdV,Q\otimes R_\bdV)$ will be denoted $\alphatw$.

\begin{rmk}
  Note that the orientation of $\cQ$ is reversed when the
representation $Q$ is replaced by its dual $Q^*$, so when $Q$ is
self-dual all the arrows $i\to j$ have opposite arrows $j\to i$. The
pair $i\rightleftarrows j$ is usually denoted $i\dsh j$.
\end{rmk}

%\subsection{The McKay Correspondence}
%\label{sec:mckay:mckay:ex}

\begin{example}[The McKay Correspondence]
\label{ex:dynkin}
For a subgroup $\ga<\SU(2)$, the standard 2-dimensional
representation $Q$ is always self-dual, since $\Lambda ^2 Q\cong \C$
as $\ga$-modules.   McKay's observation~\cite{mckay:graphs} is that
the quivers $\cQ_{\ga,Q}$ coincide with the extended homogeneous\footnote{A
  Dynkin diagram is called \emph{homogeneous} if it has
  no multiple bonds} Dynkin diagrams $\overline A_k,\overline D_k,
\overline E_6,\overline E_7,\overline E_8$ represented in
Figure~\ref{fig:dynkin}.
\begin{figure}[htbp]
\begin{center}
  \leavevmode
\epsfysize= 5cm \epsfbox{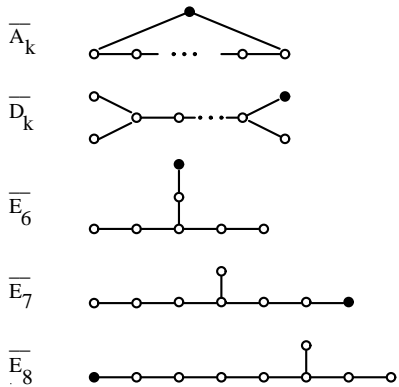}
  \end{center}
\caption[The extended homogeneous Dynkin diagrams $\overline 
A_k,\overline D_k, \overline E_6,\overline E_7,\overline E_8$]{The 
extended homogeneous Dynkin diagrams $\overline A_k,\overline D_k, 
\overline E_6,\overline E_7,\overline E_8$  (The extra vertex is 
indicated marked $\bullet$.)}
  \label{fig:dynkin}
\end{figure}
\end{example}

\subsection{The Commutation Relations}
\label{sec:mckay:mckay:ccr:relns}
\aftersub

Recall from Section~\ref{sec:mckay:mckay:dfn} that if $R_\bdV$ is any 
$\ga$-module, the $\ga$-module homomorphisms $R_\bdV \to Q\otimes R_\bdV$ 
correspond to representations of the McKay quiver into the $I$-graded 
vector-space $\bdV$ which is the multiplicity space for $R_\bdV$. 

One natural $\Gamma$-module to consider as a candidate for $R_\bdV$ is
the ``canonical one'' given by the regular representation space $R$ of 
$\Gamma$. Its  multiplicity space will be denoted $\bdR=\oplus_{v\in \cQ_0}
\bdR_v$.  Its components satisfy the well-known equalities $\dim \bdR_v=\dim R_v$.

The simplest way to define the commutation relations $\cK$ is to use
the isomorphism~\eqref{eq:hom-rep}.  Let $q_i$ be a basis of $Q$ which
is orthonormal with respect to some positive definite $\ga$-invariant
hermitian inner product. Let $\alphatw\in\Rep_\cQ(\bdR)$ be a
representation of $\cQ$ into $\bdR$, and let $\alpha\in Q\otimes\End_\C
R$ be the $\ga$-invariant element which corresponds via the
isomorphism~\eqref{eq:hom-rep}.  The group $\ga$ acts on $\End_\C R$
in the natural way, \ie by conjugation, and the element $\alpha $
decomposes with respect to the basis $q_i$ into endomorphisms $\alpha
_i\in\End_\C R$ for $i=1,\dots, n$ which satisfy the following
equivariance condition
\begin{equation}\label{eq:alpha_equivariance_condition}
  \sum_l \gamma_{kl} \alpha_{l} = 
  \varphi(\gamma)\alpha_{k}\varphi(\gamma)^{-1},\rlap{$\qquad \forall 
  k,\gamma$,}
\end{equation}
where $\bgamma=(\gamma_{kl})$ is the matrix corresponding to 
the
action of the element $\gamma$ on $Q$ with respect to the basis
$\{q_l\}_{l=1}^n$.

The commutation relations are defined by the condition
\begin{equation}
\label{eq:alpha_commute}
[\alpha _i,\alpha _j]=0.
\end{equation}
Thus the variety $\cM_{\cQ,\cK}(\bdR)$ of representations of $\cQ$ into
$\bdR$ satisfying these relations coincides with the variety $\cN^\ga$
defined in~\cite{sacha:ale,sacha:thesis}.

\subsection{Representation Moduli}
\label{sec:1:mckay:moduli}

Given $\zeta\in\frk$, the  representation moduli of the McKay
quiver of $(\Gamma,Q)$ in the multiplicity space $\bdR$ for the regular
representation $R$ of $\Gamma$ are defined by the GIT quotients
$$\cM_\zeta= \Rep_{\cQ,\cK}(\bdR) \gitquot\zeta \PGL(\bdR).$$

Since the group $\PGL(\bdR)$ coincides with the group $\PGL^\ga(R)$ (the 
projectivization of the group of $\Gamma$-invariant linear 
endomorphisms of $R$), one sees that the moduli $\cM_\zeta$ coincide 
with the moduli $X_\zeta$ of Hermitian-Yang-Mills type bundles defined 
and studied in~\cite{sacha:ale}.

\part{Abelian Groups}
\label{sec:2}

In Part~\ref{sec:1} the moduli $X_\zeta$ were shown to coincide with
the moduli of Hermitian-Yang-Mills bundles defined in~\cite{sacha:ale}
and hence ($\Gamma$ acts freely outside the origin) to provide partial
resolutions $\rho_\zeta\colon X_\zeta \to X_0=Q/\Gamma$ of the
isolated quotient singularity.

The focus from now on will be the moduli $X_\zeta$ in the case where
$\Gamma$ is an abelian group (of order $r$) acting linearly on
$Q\cong\C^n$.  (Later we shall specialise to the cyclic group of order
$r$.)

When $\ga$ is abelian, the singularity $X$ is toric and can be
resolved within the toric category, in general in many ways. In fact,
as shown in the first section below, the space of representations
$\Rep_{\cQ,\cK}(\bdR)$ and its quotients $X_\zeta$ are also
toric varieties; in  toric notation (reviewed in~\ref{sec:2:abel:toric})
$$\Rep_{\cQ,\cK}(\bdR) = \T^{\Lambda,C}\text{ and }X_\zeta =
\T^{\Pi,C_\zeta},$$
where $\Lambda,\Pi$ are certain
lattices, and where the $n$-dimensional convex polyhedra $C_\zeta$ are
obtained by `slicing' the $n+r-1$-dimensional cone $C$ by affine
$n$-planes.

The rest of this paper is devoted to describing the slices $C_\zeta $,
and in particular their extreme points and their tangent cones. The
number of extreme points is the Euler characteristic of $X_\zeta $,
and the geometry of the tangent cones of $C_\zeta $ at these points
describes the singularities of $X_\zeta $.

There are two main steps in the study of the slices $C_\zeta $. The
first step is to reduce their description to a network flow problem on
the McKay quiver.  The particular problem turns out to be a
generalisation of the classical transportation
problem~\cite{kenn_helg,gond_mino:graphs}. This is done
in Section~\ref{sec:2:abel:toric}.

The second step is to describe the polyhedra which solve the
generalised transportation problem.  A solution for a general quiver
is given in Section~\ref{sec:2:flow}.  The main result is
Theorem~\ref{thm:general}, which says that the flows which correspond
to the extreme points of $C_\zeta$ are precisely the ones whose
support contains no cycles of non-zero type in its ``closure'' (see
\*2:flow:cycle* for the definition).  This is then applied to the case
of the McKay quiver in Sections~\ref{sec:2:exact} and~\ref{sec:2:sing} to
describe the slices $C_\zeta $ for all $\zeta$.

The last section in this part chapter concerns various conjectures, in
particular for singularities of dimension~3. One hope is that, in
general, the quotients $X_\zeta $ provide a class of partial
desingularisations of the isolated quotient singularity $X$, which
will be in some sense ``natural,'' and, in good cases,
non-singular and minimal with respect to this
property.  This is motivated by the fact that
Kronheimer~\cite{kron:thesis} has shown that for the case $\Gamma
\subset\SU(2)$, the $X_\zeta $ coincide with the minimal smooth
resolution of the singularity for generic
values of $\zeta$.  One candidate for a ``good case'' is the case when
$\Gamma\subset\SU(3)$.

\begin{conj}
  If $\Gamma\subset\SU(3)$ is abelian and acts on $\C^3$ freely
  outside the origin then $\rho_\zeta\colon X_\zeta \to X_0=\C^3/\Gamma$
  is crepant and is a smooth resolution for generic values of $\zeta$.
\end{conj}

Using the description of $C_\zeta$, this conjecture is reduced in
Section~\ref{sec:2:crep} to to a statement concerning the existence
and combinatorial properties of certain trees in the McKay quiver.
Brute-force computer calculations show the conjecture to be true for
the singularities $\qsing 1/6(1,2,3)$, $\qsing 1/7(1,2,4)$, $\qsing
1/8(1,2,5)$, $\qsing 1/9(1,2,6)$, $\qsing 1/10(1,2,7)$ and $\qsing
1/11(1,2,8)$, but there is as yet no general proof.  Other conjectures
regarding, for example, the Euler number of the resolutions also
translate to graph theoretical and combinatorial statements on the
McKay quiver.

\section{Summary of Notation for Part \ref{sec:2}}
\label{sec:2:notation}

\begin{Pentry}
\item[$\ga$] Finite abelian group $\ga$ of order $r$ (usually
  $\mu_r$).
\item[$\mu_{r}$] Group of $r$-th roots of unity in $\CX$.
\item[$\gahat$] Character group $\gahat:=\Hom(\ga,\CX)$.
\item[$\rho$] The morphism giving the action of $\Gamma$ on $\CX^n$:
  $\rho\colon\ga\to\CX^{n}$.
\item[$w_{i}$] Weights ($w_i\in \Z_{r}; i=1,\dots,n$) for the action of
  $\ga$ on $\C^{n}$.
\item[$\hat{\rho}$] Dual morphism $\hat\rho\colon\Z^{n}\to\gahat$.
\item[$\Pi$] A sub-lattice of $\Z^{n}$ of index
  $r$) given by $$\Pi := \ker\hat\rho= \{ x\in \Z^n | \sum x_iw_i\equiv 0\pmod
  r\}.$$
\item[$\qsing1/r(w_1,\dots,w_n)$] The quotient singularity
  $\C^n/\mu_r$, where $\mu_r$ acts on $\C^n$ with weights $(w_1,\dots,w_n)$.
\item[$R_{v}$] One-dimensional representation on which
  $\ga$ acts by $\lambda\to\lambda^{v}$.
\item[$\nu$] Natural morphism $\nu\colon\ga\to\PGL(\bdR)$.
\item[$\hat{\nu}$] Dual morphism $\hat\nu\colon\Lambda^{0,0}\to\gahat$.
\item[$a_{v}^{i}$] Arrow of type $i$: $a_{v}^{i}:=v\to v-w_{i}$ ($v\in\cQ_{0}$,
  $i\in\{1,\cdots,n\}$).
\item[$\tilde{\alpha}(a_{v}^{i})$] Value of $\tilde{\alpha}$ on
  $a_{v}^{i}$; Equal to the $(v-w_{i},v)$-th entry of the matrix
  $\alpha$.
\end{Pentry}

\subsection{Toric Geometry}
\begin{Pentry}
\item[$M$] Generic lattice.
\item[$M_{\Q}$] Associated rational vector-space.
\item[$T^{M}$] Algebraic torus $T^{M}:=\Spec\C[M]$.
\item[$P$] Generic polyhedron in $M_{\Q}$.
\item[$\T^{M,P}$] Quasi-projective toric variety
  $\Proj\C[\Ptw\cap\Mtw]$ defined by $M$ and $P$.
\item[$\Mtw$] Product lattice $\Mtw:=\Z\times M$
\item[$\Ptw$] (Closure of) the cone over $P$:  $\Ptw :=
  \overline{\Q\nneg(\{1\}\times P)}\subset M_{\Q}$.
\item[$T_{F}P$] Tangent cone of $P$ at face $F$.
\end{Pentry}

\subsection{Flows}
\begin{Pentry}
\item[$k^A$] Set (lattice, vectorspace, cone) of maps $A\to k$.  Used for
  $k=\Z,\Z_+, \R, \R_+, \C$, and $A=\cQ_0, \cQ_1, S$.
\item[$f$] Generic flow (element of $\R^{\cQ_{1}}$).
\item[$f^{\pm}$] Positive/negative part of $f$.
\item[$\Lambda^{1}$] Lattice $\Z^{\cQ_{1}}$ of integer-valued flows on
  $\cQ$.
\item[$\Lambda^{1}_{\R}$] Vectorspace of real-valued flows on $\cQ$.
\item[$\Lambda^{1}_{\R_{+}}$] First quadrant in $\Lambda^{1}_{\R}$
  (non-negative real-valued flows).
\item[$\chi_{a}$] Flow taking the value 1 on $a$ and zero elsewhere.
  Basis element of $\Lambda^{1}$.
\item[$\chi_{v}^{i}$] Alternative notation for $\chi_{a_{v}^{i}}$.
\item[$\Lambda^{0}$] Lattice $\Z^{\cQ_{0}}$ of assignments of integers
  to the vertices of $\cQ$; $\Lambda^0_\R$ coincides with $\Lie \GL(\bdR)$.
\item[$\Lambda^{0,0}$] Sub-lattice $\{\zeta\in\Z^{\cQ_{0}}|
  \sum_{v\in\cQ_{0}}\zeta(v)=0\}$ of $\Lambda^0$; $\Lambda^{0,0}_\R$
  coincides with $\Lie \PGL(\bdR)$.
\item[$\zeta$] Element of $\Lambda^{0,0}$.
\item[$\In(v)$] Set of arrows $a$ such that $h(a)=v$.
\item[$\Out(v)$] Set of arrows $a$ such that $t(a)=v$.
\item[$\partial$] Natural map
  $\partial\colon\Lambda^{1}\to\Lambda^{0,0}$ which calculates the net
  contribution of a flow at each vertex:
  $$\partial f(v)=\sum_{a\in\In(v)} f(a) - \sum_{a\in\Out(v)} f(a).$$
\item[$\chi_{v}$] Function taking the value 1 on $v$ and zero
  elsewhere.  Basis element of $\Lambda^{0}$.
\item[$\Lambda^{2}$] Sub-lattice of $\Lambda^{1}$ generated by the
  elements corresponding to the commutation relations
  $$\chi_v^i +\chi_{v- w_i}^j-\chi_v^j-\chi_{v- w_j}^i,\qquad\text{for
    } i,j=1,\dots, n$$
\item[$\Lambda$] Quotient lattice $\Lambda^{1}/\Lambda^{2}$.
\item[$C$] Image of cone $\Lambda^{1}_{\R_{+}}$ inside the quotient
  $\Lambda$.  Also written $\Lambda_{\R_{+}}$.
\item[$C_{\zeta}$] Convex polyhedron obtained by slicing $C$ with the
  hyper-plane $\partial f=\zeta$.  Equal to the image of $F_{\zeta}$
  under $\pi$.
\item[$\pi$] Map $\pi\colon\cQ_1\to\{1,\dots,n\}$ which assigns to each
  arrow of the McKay quiver its type: $\pi(a^i_v) := i$.  Induces a
  map $\pi\colon\Lambda^{1}_{\R_{+}}\to\R^{n}$ on the flows:
  $$\pi(f)=\left(\sum_{v\in\cQ_{0}} f(a_{v}^{1}),\dots,\sum_{v\in\cQ_{0}}
  f(a_{v}^{n})\right).$$
\end{Pentry}

\subsection{Configurations}
\begin{Pentry}
\item[$S$] Configuration of arrows (non-empty subset of $\cQ_{1}$).
\item[$\supp f$] Support of the flow $f$ (arrows on which $f$ takes a
  non-zero value).
\item[$\cC$] Set of all configurations of arrows.
\item[$\cT$] Set of all configurations of arrows which are trees (\ie
  which contain no cycles).
\item[$F_{\zeta}$] Admissible (\ie non-negative) flows which satisfy
  $\partial f=\zeta$.
\item[$F_{\zeta}(S)$] Flows satisfying $\partial f=\zeta$ which are
  non-negative outside $S$.  (A convex cone if $\zeta=0$).
\item[$Z_{\zeta}(S)$] Flows satisfying $\partial f=\zeta$ which are
  zero outside $S$. (Equal to the maximal subspace contained in
  $F_{0}(S)$ if $\zeta=0$).
\item[$f_{\zeta}$] Map $f_{\zeta}\colon\cT\to\Lambda^{1}_{\R}$
  assigning to each tree $T$ the unique flow $f=f_{\zeta}(T)$ such
  that $\partial f=\zeta$.
\item[$\Adm(T)$] Open convex cone in $\Lambda^{0,0}_{\R}$ of values
  of $\zeta$ for which the tree $T$ is the support of some flow in
  $F_{\zeta}$.
\item[$\cD$] Generic subset of $\cC$.
\item[$\cD_{\zeta}$] Elements of $\cD$ which are the support of some
  $f\in F_{zeta}$.
\item[$\cD_{\spn}$] Elements of $\cD$ which span all vertices of
  $\cQ$.
\item[$\cD^{k}$] Elements of $\cD$ of rank $k$ (\ie elements $S$ such that $\rk \pi F_{0}(\clos{S})=k$.
\end{Pentry}

\subsection{Cycles and Paths}
\begin{Pentry}
\item[$p$] Generic path in $\cQ$.
\item[$p^{+}$] Positive part of $p$.
\item[$p^{-}$] Negative part of $p$.
\item[$\tilde{\chi}_{p}$] Basic flow associated to the path $p$.
\item[$c$] Generic cycle in $\cQ$.
\item[$\{v\}(j_{0},\dots,j_{k})$] Notation for basic flows.
\end{Pentry}

\subsection{Miscellaneous}
\begin{Pentry}
\item[$\cW_S$] Weighting of the vertices of $\cQ_{0}$ (a map
  $\cQ_{0}\to\Z^{n}$ with special properties).
\end{Pentry}

\section{Representations and Relations in the Abelian Case}
\label{sec:2:abel}
\subsection{The McKay Quivers}
\label{sec:2:abel:mckay}
%\label{sec:mckay:ab:quiv}

In the abelian case, all the irreducible representations are
one-dimensional and are determined by an element of the dual group
$\gahat$, which is a finite abelian group isomorphic to $\ga$. Thus
the vertex set of $\cQ$ is nothing but $\gahat$. The group $\ga$ will
always be identified with a subgroup of $\C^{*n}\subset\GL(n)$, and thus
$\gahat$ will be thought of as a product of finite groups $\Z_{r_i}$ with the group
operation denoted additively.
\begin{example}
\label{ex:5123}
Let $\ga=\mu_5\subset\CX$, the group of $5$-th roots of unity, acting 
on $\C^3$ with \emph{weights} $1,2$ and $3$, \ie via 
$$\lambda\to 
\begin{pmatrix}
\lambda^1 & 0 & 0 \\
0 & \lambda^2 & 0 \\
0 & 0 & \lambda^3
\end{pmatrix}.
$$   This action is denoted symbolically by $\frac{1}{ 5}(1,2,3)$.  The
character group of $\mu_5$ is $\widehat{\mu_5}=\Z_5:=\Z/5\Z $, so
$\cQ$ has $5$ vertices. For each $i\in\Z_5$, let $\chi_i\colon
\lambda\mapsto\lambda^i$ be the corresponding character, so that 
$\chi_i\chi_j=\chi_{i+j}$. Writing $R_i$ for the irreducible 
representations, one has $R_i\otimes R_j = R_{i+j}$. The representation 
$Q$ is just $R_1\oplus R_2\oplus R_3$, so that $$Q\otimes 
R_i=R_{i+1}\oplus R_{i+2}\oplus R_{i+3}.$$ Thus, the quiver has arrows 
from the vertex $i$ to the vertices $i-1, i-2$ and $i-3$ for each 
$i$ (see Figure~\ref{fig:all5123} below).\footnote{The additions are modulo 5.}
%%%
% QUIVER 5 123
%%%%
\begin{figure}[htbp]
  \begin{center} 
\leavevmode 
\epsfysize= 3cm \epsfbox{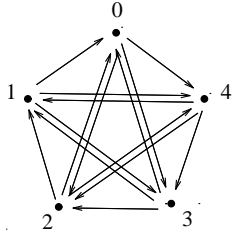}
  \end{center}
  \caption{The McKay quiver for the action $\frac{1}{ 5}(1,2,3)$}
  \label{fig:all5123}
\end{figure}
\end{example}

In fact, it is easy to see that the McKay quiver for the action of any
cyclic group has a similar appearance.  To see what the arrows are, decompose $Q$ into a sum of 
one-dimensional irreducibles 
\begin{equation}
\label{eq:qdecomp}
Q=\bigoplus_{i=1}^n R_{w_i},
\end{equation}
where $w_i\in\gahat$ for $i=1,\dots, n$ are  the \emph{weights} of the action of $\Gamma$ on $Q$. Using 
~\eqref{eq:irred_decomp}, one sees easily that 
\begin{equation} 
Q\otimes R_v = R_{v+w_1}\oplus \dots\oplus R_{v+w_n}, 
\label{eq:qtensorrv}
\end{equation}
so there is one arrow from $v$ to $v-w_i$ for each vertex $v$ and 
each weight $w_i$, giving a total of $nr$ arrows.  The arrows
corresponding to the weight $w_i$ are written $$a_v^i :=
v\to v -w_i,\quad \text{for } v\in\cQ_0$$  and are said to be of
\emph{type}~$i$. 

The McKay quiver for a general abelian group $\ga$ can be obtained by 
decomposing $\ga$ into products of cyclic groups. Define the 
\emph{product} of $\cQ$ and ${\cQ'}$ to be the quiver 
with vertices $\cQ_0\times{\cQ'}_0$ and with arrows 
\begin{multline}
  \{(v,t(a'))\to (v,h(a')) : v\in\cQ_0,a'\in{\cQ'}_1\}\\
\cup \{( t(a),v')\to (h(a),v') : v'\in{\cQ'}_0,a\in\cQ_1\}.
\end{multline}
Then the McKay quiver for $\ga$ is given by taking the product of the 
quivers for the cyclic factors.

\subsection{Commutation Relations}
\label{sec:2:abel:comm}

Pick a basis $q_i$ of $Q$ such that the action of $\ga$ on $Q$ 
is diagonal, with $q_i$ corresponding to the irreducible component 
$R_{w_i}$. Then  decomposing $R$ into irreducibles and using Shur's Lemma
\begin{equation}
\begin{split}
\Hom_\ga(R,R_{w_i}\otimes R) 
&= \bigoplus_{v,v'\in\cQ_0}\Hom_\ga({\bdR_v}\otimes R_v,\bdR_{v'}\otimes 
R_{w_i}\otimes R_{v'})\\
&= \bigoplus_{v\in\cQ_0\phantom{v'}} \Hom({\bdR_v},\bdR_{v-w_i })
\end{split}
  \label{eq:homga}
\end{equation}
This means that the scalar map $\alphatw(a_v^i)\colon \bdR_{v}\to \bdR_{v-w_i 
}$ is multiplication by the $(v-w_i,v)$-th entry of the $i$-th 
component matrix $\alpha_i\in\End R$. The condition 
$[\alpha_i,\alpha_j]=0$ therefore translates to the commutation relation 
\begin{equation}
\label{eq:comrelations2}
\alphatw(a_v^i)\alphatw(a_{v-w_i }^j)-\alphatw(a_v^j)\alphatw(a_{v-w_j }^i)=0
\end{equation}
for all $v\in\cQ_0$ and $i,j\in\{1,\dots, n\}$.
\begin{figure}[htbp]
\begin{center} 
\leavevmode 
%\special{epsf: comm7_123.eps}
\epsfysize= 3cm \epsfbox{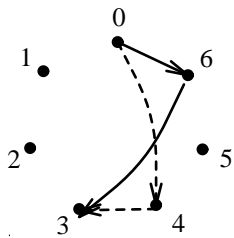} 
\end{center} 
\caption[A picture of the commutation 
relation~\eqref{eq:comrelations2}]{A picture of the commutation 
  relation~\eqref{eq:comrelations2} for $v=0$, $i=1$ and $j=3$ in the
   McKay quiver for $\frac{1}{ 7}(1,2,3)$. (The product of the representation along the continuous arrows must equal the product along the dotted arrows.)}
\label{fig:comm7_123} 
\end{figure}

\subsection{Toric Descriptions of the Representation Moduli}
\label{sec:2:abel:toric}
%\label{sec:mckay:ab:toric-d-moduli}

\subsubsection{Toric Varieties}

Let $M$ be an integral lattice isomorphic to $\Z^n$,  let
$M_\Q:=M\otimes\Q$ be the associated rational vector
space, and let $C\subset M_\Q$ be a convex $n$-dimensional cone.

Let $T^M:=\Spec\C[M]$ be the algebraic torus whose
character group is $M$.  Then 
$$\T^{M,C}:=\Spec\C[C\cap M]$$
is a compactification of $T^M$, called
the \emph{toric variety associated to the cone} $C$ with respect to the
lattice $M$.  The covariant notation $\T_{M^\vee, C^\vee}$ which uses
the dual objects\footnote{The dual of a cone $C$ is the cone $C^\vee
  := \{n\in M^\vee | n(c) \geq 0, \forall c\in C\}$} is widely in use in the literature, but the
contravariant version is more convenient for the purposes of this
paper.

This definition can be extended from cones to general convex polyhedra
as follows. If $P$ is any convex polyhedron in $M_\Q$, denote by
$\Ptw$ the closure of the cone over $P$, namely
\begin{equation}
  \label{eq:Ptw}
  \Ptw :=  \overline{\Q\nneg(\{1\}\times P)} \subset \Mtw_\Q=\Q\times
M.
\end{equation}
Then one can define:
\begin{equation}
  \label{eq:toric-defn}
  \T^{M,P}:=\Proj\C[\Ptw\cap\Mtw],
\end{equation}
and this is easily seen to extend the previous definition for cones.

\subsubsection{GIT quotients of Toric Varieties}

Note that there is a natural line bundle $L^{M,P}:=\cO(1)$ on $\T^{M,P}$
on which $T^M$ acts and whose sections are given simply by evaluation
on the lattice points of $P$.  If one translates $P$ by an element
$\zeta\in M$, one gets the same toric variety
($\T^{M,P+\zeta}=\T^{M,P}$) and the same line bundle $L$, but with a
\emph{different \/} linearisation of the action of $T^M$, differing
from the old by the character $\zeta$.  

The consequence of this for GIT quotients is that if $M'$ is any
sub-lattice of $M$, and if $\zeta$ an element of $M'$ (and hence a
character of $T^{M'}$) one can consider either the quotient
$$\T^{M,P}\gitquot{\zeta}T^{M'}$$
with  $T^{M'}$ taken to act on $L^{M,P}$ via multiplication by the
section $\zeta$ or, equivalently, the quotient
$$T^{M,\zeta+P}\gitquot{} T^{M'}$$ where $T^{M'}$ now acts on
$L^{M,\zeta+P}$ in the ordinary way.

Using this notation, one can state the following proposition which
describes how toric varieties behave under GIT quotients
(\cf~\cite{thaddeus:flips}).
\begin{prop}
\label{prop:toric_quotients}
  Let $0\to M''\to M\to M'\to 0$ be an exact sequence of lattices and 
$P\subset M_\Q$ be a convex polyhedron. Write 
$$(M,P)\gitquot{}M' := 
(M'',P\cap M_\Q'').$$ Then  
  $$\T^{ M,P}\gitquot{} T^{ M'} = \T^{(M,P)\gitquot{}M'}.$$
\end{prop}

\begin{proof}
  An element $x^m\in\C[ M\cap P]$ is invariant under $T^{ M'}$ if $m$
  belongs to the stabiliser of $T^{ M'}$ in $T^M$, which is nothing
  but $M''$. Thus 
  $\C[\Ptw\cap \Mtw]^{T^ {M'}} =\C[\Ptw\cap\widetilde{M''}]= \C[\widetilde{(P\cap M_\Q'')}\cap\widetilde{M''}]$ and the result follows by taking
  $\Proj$.
\end{proof}

\subsubsection{The Representation Variety}

The description of the previous sections shows that
$\Rep_{\cQ,\cK}(\bdR)$ is a subvariety of $\C^{\cQ_1}$ defined by the
equations~\eqref{eq:comrelations2}.  The goal is now to describe the quotients
by the group $\PGL(\bdR)$.  It turns out that $\Rep_{\cQ,\cK}(\bdR)$ is 
a toric variety and that $\PGL(\bdR)$ is an algebraic sub-torus, so that one is
able to use Proposition~\ref{prop:toric_quotients}.

Consider the lattice $\Lambda ^1:= \Z^{\cQ_1}$ and the
first quadrant $\Lambda ^1_{\R+}:=\R_+^{\cQ_1}$ inside its associated
real vector-space $\Lambda^1_\R$. Recall that the toric variety which
corresponds to $\C^{\cQ_1}$ is
$$\T^{\Lambda ^1,\Lambda ^1_{\R+}} := \Spec \C[\Lambda ^1_+],$$
where $\Lambda ^1_+$ is the semi-group $\Z_+^{\cQ_1}$ and $\C[\Lambda 
^1_+]$ denotes its group algebra.  More precisely, the $\C$-points of 
the scheme $\Spec \C[\Lambda ^1_+]$ are into one-one correspondence 
with the points of $\C^{\cQ_1}$ as follows. If $x\colon \C[\Lambda 
^1_+]\to\C$ is an algebra homomorphism representing a point of $\Spec 
\C[\Lambda ^1_+]$ then, evaluating it on the generators of $\Lambda 
^1_+$ gives a map $\cQ_1\to \C$, \ie a point of $\C^{\cQ_1}$. 
Conversely, if $\alphatw \in\C^{\cQ_1}$,  there is an induced morphism 
of semi-groups 
\corresp{(\Lambda ^1_+,+)}{(\C,\cdot)}{f}{\prod_a(\alphatw (a))^{f(a)},}
and this induces an algebra morphism $\C[\Lambda ^1_+]\to\C$.

Under this identification, the points of $\Rep_{\cQ,\cK}(\bdR)$ correspond to the
semi-group morphisms $ \Lambda ^1_+ \to\C$ which are the identity when
restricted to the 
sub-semi-group $\Lambda ^2\cap\Lambda ^1_+$ of $\Lambda ^1_+$.  Here 
$\Lambda ^2$ is the sub-lattice of $\Lambda ^1$ generated by the 
elements 
$$\chi_v^i +\chi_{v- w_i}^j-\chi_v^j-\chi_{v- w_j}^i,\qquad\text{for }
i,j=1,\dots, n$$
(and $\chi_v^i$ denotes the indicator function
$\chi_{a_v^i}$ for the element $a_v^i$, \ie the basis element of
$\Lambda^1$ which takes the value $1$ on the arrow $a_v^i=v\to
v-w_i\in \cQ_1$ and zero elsewhere). 
\begin{prop}
\label{prop:rep-variety-toric}
  The quotient $\Lambda:=\Lambda ^1/\Lambda ^2$ is a lattice (\ie 
  is torsion free).  If one denotes by $\Lambda _+$ (\resp 
  $C=\Lambda_{\R+}$) the lattice (\resp the cone) generated by the
  image of the elements of $\Lambda ^1_+$ in $\Lambda$, then the
  variety of representations of the McKay quiver with
  relations  into $\bdR$ is given by
 $$\Rep_{\cQ,\cK}(\bdR) = \T^{\Lambda,C} = \Spec \C[\Lambda _+].$$
\end{prop}
\begin{proof}
  The only fact which has to be proved is that $\Lambda $ is a
  lattice.  This will be postponed till
  Lemmas~\ref{lemma:kerc_lambda2}--\ref{lemma:imagec}.
\end{proof}

The group $\GL(\bdR)$ on the other hand coincides with
$\CX^{\cQ_0}$.  If $\Lambda^0$ denotes the lattice
$\Z^{\cQ_0}$, then one can write
$\GL(\bdR)=T^{\Lambda^0}$ in the notation of Section~\ref{sec:2:abel:toric}.
Similarly, writing  $\Lambda^{0,0}$ for the 
sub-lattice of $\Lambda^0$ whose elements $\zeta$ satisfy
$\sum_{v\in\cQ_0}(v)=0$, one has  $\PGL(\bdR)= T^{\Lambda^{0,0}}$.

 Applying Proposition~\ref{prop:toric_quotients}, one obtains the
 following toric description of the moduli $X_\zeta$.
\begin{thm}
\label{thm:toric-d-moduli}
The moduli spaces $X_\zeta$ are quasi-projective toric varieties
$$X_\zeta = \T^{\Pi ,C_\zeta },$$
where
$\Pi:=\ker\hat\rho=\pi\times\partial(\Lambda)/\Lambda^{0,0}$ is a
sub-lattice of $\Z^n$ of index $\card{\Gamma}$ and $C_\zeta$ are convex
polyhedra in $\Pi_\R$ in obtained by slicing the cone $C$ with the
affine planes $\zeta +\Lambda ^{0,0}_\R$.
\end{thm}
\begin{proof}
  By definition
  $$X_\zeta=\T^{\Lambda,C}\gitquot{\zeta}T^{\Lambda^{0,0}},$$
  and this
  coincides with $$\T^{\Lambda,\zeta+C}\gitquot{}T^{\Lambda^{0,0}}.$$
  If $\Lambda':=\pi\times\partial(\Lambda)/\Lambda^{0,0}$ was a
  lattice then applying Proposition~\ref{prop:toric_quotients} would
  give $X_\zeta = \T^{\Lambda' ,C_\zeta }$.
  
  It therefore remains to prove that $\Lambda'=\Pi$.  This is done in
  Section~\ref{sec:2:exact}:
  Lemmas~\ref{lemma:kerc_lambda2}--\ref{lemma:imagec} show that there
  is an exact sequence of abelian groups
  \begin{equation}
    \label{eq:exact1}
    0\to\Lambda\xrightarrow{\pi\times\partial} \Z^n\times\Lambda^{0,0}
  \xrightarrow{\hat\rho-\hat\nu} \Z^n/\Pi \cong\gahat \to 0.
  \end{equation}
\end{proof}

Note that the morphisms appearing in the exact sequence in the above
proof are:
\begin{itemize}
\item The projection $\pi\colon\Lambda^1\to\Z^n$ (or better, its descent to
  $\Lambda\to\Z^n$).
\item The (descent to $\Lambda$) of the morphism of lattices
  $\partial\colon\Lambda^1\to\Lambda^{0,0}$ dual to the action of
  $\PGL(\bdR)=T^{0,0}$ on $\Rep_{\cQ}(\bdR)=\C^{\cQ_1}$ given by
  \map{\Hat\partial}{T^{0,0}}{T^1}{\lambda}{a\mapsto\lambda_{t(a)}^{-1}\lambda_{h(a)}.}
  Here $\lambda_v$ denotes the component of
  $\lambda\in T^{0,0}$ corresponding to the vertex $v\in\cQ_0$.
\item The morphism of lattices
  \map{\hat\nu}{\Lambda^{0,0}}{\gahat}{\zeta}{\sum_{v\in\cQ_0}\zeta(v)v}
  dual to the action of \/$\Gamma$ on $\End R$ by conjugation.
\item The morphism of
  lattices\map{\hat\rho}{\Z^n}{\gahat}{\zeta}{\sum_{v\in\cQ_0}x_iw_i}
  dual to the action $\rho\colon \Gamma\to \CX^n$.
\end{itemize}

The fact that $\pi\times\partial$ descends to $\Lambda$ corresponds to
the fact that $\GL^\ga(Q)=\CX^n$ and $\PGL(\bdR)=T^{0,0}$ act on
$\Rep_{\cQ,\cK}(\bdR)$ rather than simply on
$\Rep_{\cQ}(\bdR)=\C^{\cQ_1}$.  On the other hand, the fact that their
image maps into $\ker (\hat\rho-\hat\nu)$ corresponds to the fact that
$\Rep_{\cQ,\cK}(\bdR)$ can be identified with the $\Gamma$-invariant
part of $\Hom(R,Q\otimes R)$.

The action of $\Tn$ arises from the existence of the map
$$\pi\colon\cQ_1\to \{1,\dots, n\},$$ which assigns to each arrow its 
type. 
This induces a $\Z$-linear map $\pi 
\colon {\Lambda^1}\to \Z^n$, which is denoted by 
the same letter, and one obtains an action of $\CX^n$ on 
$\C^{\cQ_1}$ via the corresponding morphism of algebraic tori: 
$\tau\cdot\alphatw:=\Hat\pi (\tau)\alphatw$. 
Explicitly, for $\tau\in\CX^n$, $\alphatw \in\C^{\cQ_1}$, and 
$a\in\cQ_1$,
$$  (\tau\cdot \alphatw)(a) 
%%%%%%= \Hat\pi\tau(\chi_a)\alphatw(a)
%%%%%%=\tau(\pi\chi_a)\alphatw(a)=\tau(e_{\pi (a)})\alphatw(a)
=\tau_{\pi(a)}\alphatw(a),
$$
where $\tau\in \CX^n$ has components $\tau_i:=\tau(e_i)$ with
respect to the standard basis $e_i$ of $\Z^n$, and $\chi _a$ denotes
the basis element of $\Lambda ^1$ which is the indicator function of
the singleton $\{a\}\subset\cQ_1$. The action of $\CX^n$ corresponds
to multiplying the value of $\alphatw $ on all arrows $a$ with the
same $\pi(a)$ by the same factor $\tau_{\pi(a)}$.  Since $\Gamma$ is
abelian, $\rho(\ga)\subset\CX^n$ acts trivially on the
representations, and one has an action of $\Tn:=\CX^n/\rho(\ga)$.
According to the toric formalism, one has
$$\Tn = \CX^n/\Image\rho = T^{\ker\hat\rho},$$
and thus $\Pi=\ker\hat\rho$.

Note also that the exact sequence implies that the cone
$C:=\Lambda_{\R+}$ is isomorphic to
$({\pi\times\partial})(\Lambda^1_{\R+})$.  Therefore, one has the
following theorem.

\begin{thm}
\label{thm:slices}  The slices $C_\zeta$ (which describe
the toric moduli $X_\zeta$ when $\zeta$ is integral) are given by 
\begin{align}
\label{eq:czeta}
C_\zeta &\phantom{:}= \pi_\R(F_\zeta ),\\
\intertext{ where }
\label{eq:fzeta}
F_\zeta &:= \R_+^{\cQ_1}\cap \partial^{-1}_\R(\zeta ).
\end{align}
\end{thm}

The next section is devoted to the study of the polyhedra $F_\zeta $ 
and $C_\zeta $.

\section{Flow Polyhedra}
\label{sec:2:flow}

The question of determining the polyhedra $C_\zeta$ is in fact a
generalisation of a basic problem, well-known to network optimization
specialists:
\begin{problem}[The Transportation Problem ({\protect\cite{kenn_helg,gond_mino:graphs}})]
\label{pb:transp}
Given a quiver --- a ``network" in optimization parlance --- $\cQ$
and an element $\zeta\in\R^{\cQ_0}$ specifying given demands and
supplies of some commodity at the vertices, find
all\footnote{Actually, the interest in the transportation problem
  usually centers around finding the extreme flow which minimizes a
  certain cost function, not in finding \emph{all\/} extreme flows.
  Furthermore, there are usually capacity constraints for the maximal
  flow allowed along each arrow.} the non-negative elements
$f\in\R_+^{\cQ_1}$, representing flows of commodities along each
arrow, whose net contribution at each vertex  balances the demands and supplies specified by $\zeta$.
\end{problem}
Writing $\In(v)$ ($\Out(v)$) for the arrows
which have their head (tail) at a given vertex $v$, then
$\partial\colon\R^{\cQ_1}\to\R^{\cQ_0}$ is given by
$$\partial f(v)=\sum_{a\in\In(v)} f(a) - \sum_{a\in\Out(v)} f(a).$$
The element $\partial f$ is called the \emph{net contribution} of $f$.
The flow $f\in\R^{\cQ_1}$ is said to be a \emph{$\zeta$-flow} if
\begin{equation}
\partial f = \zeta.
\label{eq:demand}
\end{equation}
 If $f$ is
non-negative,  it is called an
\emph{admissible $\zeta$-flow}.

\begin{figure}[htbp]
\begin{center} 
\leavevmode 
\epsfysize= 3cm \epsfbox{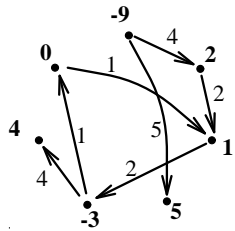} 
\end{center} 
\caption[An admissible flow on the McKay quiver for $\frac{1}{ 
7}(1,2,3)$ and its net contribution at each vertex]{An admissible 
flow on the McKay quiver for $\frac{1}{ 7}(1,2,3)$. (The numbers along 
the edges indicate the values of
  $f$ and the numbers at the vertices
  indicate the values of $\partial f=\zeta$.)}
\label{fig:flow7123} 
\end{figure}

The problem is therefore to determine the set $F_\zeta$ of all
admissible flows for a given $\zeta$. Note that the solution set is
empty unless $\zeta$ satisfies $\sum_{v\in\cQ_0}\zeta(v)=0$ (\ie 
supply = demand).  Recall that the sub-space for which this holds was
denoted $\Lambda^{0,0}_\R$.  The set $F_\zeta$ is the
intersection of the cone $\R_+^{\cQ_1}$ (the first quadrant) with an
affine translate of the vector-space $\ker\partial$, so $F_\zeta$ is a convex
polyhedron; the interest lies in determining its extreme points. For
instance, if one is trying to minimize a convex cost function of the
flows, then the minimum will be attained at one of these extreme
points. A solution to this problem is given by the easy:
\begin{thm}\footnote{This theorem is the basis for the ``Simplex on a graph'' algorithm~\cite{luen}. See Proposition~3.16 in \cite{kenn_helg}, noting that trees correspond to what is called a ``basic solution'' or a ``basis'' for the linear programming problem~\cite[\S 2.3, Defn., p.17]{luen}.  For the case of the permutation polytope, see~\cite[\S 5, Th.~1.1]{yem_kov_kra}}
\label{thm:classical}
The extreme points of $F_\zeta$ are precisely the admissible
$\zeta$-flows whose support contains no cycles.
\end{thm}
Here, the \emph{support} of a flow is the set of
arrows on which it takes non-zero values. A \emph{cycle}  means a
sequence of arrows in $\cQ_1$ which form a cycle in the underlying
graph to $\cQ$, when their orientation is disregarded.
\begin{figure}[htbp]
  \begin{center}
    \leavevmode
\epsfysize= 3cm \epsfbox{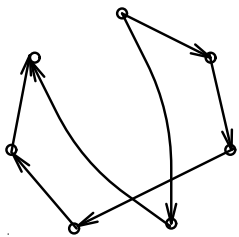}
  \end{center}
  \caption{A cycle in the McKay 
quiver for $\frac{1}{7}(1,2,3)$}
  \label{fig:uncycle7123}
\end{figure}

\section{Generalised Transportation Problem}
\label{sec:2:gtp}

In the present paper, 
the interest is in  the image $C_\zeta$ of
$F_\zeta$  under the projection $\pi\colon \R^{\cQ_1}\to\R^n$,  and
one is therefore let to the  following slightly more general problem.
\begin{problem}
\label{pb:labeledq}
Given a quiver $\cQ$, a natural integer $n$ and  a projection
$\pi\colon\R^{\cQ_1}\to\R^{n}$, characterize the (supports of the)
flows which, under $\pi$, map to an extreme point of $\pi F_\zeta $
for a given $\zeta$.
\end{problem}
\begin{rmk}
  This problem can be thought of as a transportation problem with some
  extra structure: instead of associating a cost with \emph{each\/}
  arrow, one supposes that there are linear relations between the arrow
  costs which are parametrised by $n$ independent variables --- in
  other words, suppose that the cost function really is  a (convex) function of these $n$
  variables, so that the  cost-minimizing flows will correspond
  to extreme points of $\pi F_\zeta$.
\end{rmk}

If $f$ is any flow with $\pi(f)=x\in \pi F_\zeta $, then the support 
$S$ of $f$ will be called a $\zeta$-\emph{configuration (of arrows)} for 
$x$, or simply a \emph{configuration} of $x$ (there will in general 
be \emph{many \/} $\zeta $-configurations for a given $x$). If $x$ is 
an extreme point of $\pi F_\zeta$, then any configuration $S$ of $x$ 
will be called an \emph{extreme configuration}. The  problem is thus to 
determine all the possible extreme configurations. 

\subsection{Cycle Type and Closure}
\label{sec:2:flow:cycle}

In order to state the solution to the problem, one needs to introduce
some concepts relating to cycles. Suppose that $c=(c_1,\dots, c_k)$ is
a cycle, \ie a sequence of arrows $c_i\in\cQ_1$ such that they form
a circuit, when suitably oriented. Consecutive arrows $c_i$'s are not allowed 
to be the same (although they can join the same
vertices)\footnote{Note that it is necessary to state this explicitly
  for quivers, due to the possibility of several arrows joining the
  same two vertices.}. An arrow $c_i$ is said to \emph{agree} with the ordering
$c_1,\dots c_k$ if the tail of $c_i$ is an extremity of $c_{i-1}$ and
the head of $c_i$ is an extremity of $c_{i+1}$, and is said to
\emph{disagree} otherwise. Call the disjoint union\footnote{The
  disjoint union is
  used in case some cycles contain the same arrow twice.} of the
arrows $c_i$ whose direction agrees (disagrees) with the ordering
$c_1,\dots,c_k$ the \emph{positive\/} (\emph{negative\/}) part of $c$
and denote it by $c^+$ ($c^-$). 

For each arrow $a\in\cQ_1$, let $\chi_a$ denote the basis element of
$\Z^{\cQ_1}$ which takes the value $1$ on the arrow $a$ and zero
elsewhere.  To each cycle $c$, one can define the \emph{basic flow}
 associated to $c$ by
$$\tilde\chi_c:=\sum_{a\in c^+}\chi_a - \sum_{a\in c^-}\chi_a\in\Z^n.$$

The \emph{type} of $c$ is defined to be element $\pi(\tilde\chi_c)\in\Z^n$.
Cycles of
\emph{type zero} are of special interest.
They correspond to cycles having zero total number of arrows of any
given type, the number being counted algebraically, according to
whether the arrow agrees or disagrees with the
orientation of $c$.
\begin{figure}[htbp]
  \begin{center}
    \leavevmode
\epsfysize= 3cm \epsfbox{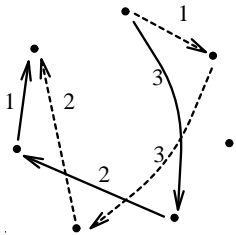}
  \end{center}
  \caption[A cycle of type zero in the McKay quiver for $\frac{1}{ 
7}(1,2,3)$]{A cycle of type zero in the McKay quiver for $\frac{1}{ 
7}(1,2,3)$. (The numbers next to the arrows indicate their type.)}
  \label{fig:nullcycle7123}
\end{figure}

The \emph{cycle closure} of a configuration $S$ is defined to be the 
smallest over-set $\SC$ of $S$ such that, for any cycle $c$ in 
$\cQ$ of type $0$,
\begin{equation}
c^+\subseteq \SC\iff c^-\subseteq \SC.
  \label{eq:cycle_closure}
\end{equation}
The closure of $S$ can be computed by searching for all the
cycles $c$ of type $0$ satisfying $c^+\subseteq S$, adjoining $c^-$ to
$S$, and repeating this procedure until~\eqref{eq:cycle_closure} is
satisfied. Note that, even if $S$ contains no cycles of non-zero type
(for instance, if $S$ is a tree), the arrows one adjoins may create
such cycles in $\SC$. 

Two configurations $S,S'$ will be called
\emph{equivalent} (written $S\sim S'$) if
$\clos{S}=\clos{S'}$.

\subsection{Statements of the Generalised Theorems}
\label{sec:2:flow:state}

The generalised version of
Theorem~\ref{thm:classical} can now be stated.
\begin{thm}[Generalised Extreme Flows]
\label{thm:general}
The extreme points of the polyhedron $\pi F_\zeta$ are the images under $\pi$ of the 
admissible $\zeta$-flows whose supports have no cycles of non-zero
type in their closures. If $S$ is a such a configuration 
and if there is a $\zeta$-flow with support $S$, then the image of 
that flow under $\pi$ is an extreme point of $\pi F_\zeta$. 
\end{thm}
Note that if $\pi$ is the identity map, then \emph{any\/} cycle is of
non-zero type and one recovers Theorem~\ref{thm:classical}.

The above theorem can be generalised to get a description of the faces
of $\pi F_\zeta $ of all dimensions. Recall that the \emph{tangent
  cone\/} of $P$ at one of its faces $F$ is the convex cone $$T_F
P:=\R_+(P-F):=\R_+(P-f), \text{ for any }f\in \interior F,$$
where $\interior F$ denotes the relative interior of the face $F$.

For any configuration of arrows  $S$ define 
$$
F_\zeta(S) := \{f\in\partial^{-1}(\zeta) : a\not\in S \implies
f(a) \geq 0\}$$
and
$$Z_\zeta(S) := \{f\in\partial^{-1}(\zeta) : a\not\in S
\implies f(a) = 0\}.$$
One has $F_\zeta(\emptyset)=F_\zeta$ and
$Z_\zeta(S)=F_\zeta\cap\supp^{-1}(S)$.  When $\zeta=0$, $F_0(S)$ is a
cone and $Z_0(S)$ is its maximal vector subspace.  The cone
$F_0(S)$ (\resp the vector-space $Z_0(S)$) is generated
over $\R_+$ by the flows $\tilde\chi_c$  for cycles $c$ such that
$c^-\subseteq S$ (\resp $c\subseteq S$).

Let $\cC$ denote the set of all \emph{configurations},\/ \ie
the set of non-empty subsets of $\cQ_1$.  The
\emph{rank} of $S$ is defined to be the rank
of $\pi Z_0(\clos S)$.  It is trivial to see that the rank
function determines a partition of $\cC$ into non-empty disjoint sets
$\cC=\cC^0\amalg\cC^1\amalg\dots\amalg\cC^n$ which respects the
equivalence relation $\sim$ induced by $S\mapsto\clos{S}$.  As a
further piece of notation, for any subset of $\cD\subseteq\cC$, denote
by  $\cD^k$ the subset of configurations in $\cD$ which
have rank $k$.  Also, write
$\cD_\zeta$ for the subset of
configurations which are
\emph{admissible for $\zeta$},\/ namely
configurations $S\in\cD$ which arise as the support of some element in
$F_\zeta$.

Theorem~\ref{thm:general} says that the configurations corresponding
to the extreme points of $F_\zeta$ are precisely those belonging to
$\cC^0_\zeta$. In general one has the following complete description
of the extreme faces and tangent cones of $C_\zeta$:

\begin{thm}[Extreme Faces and Tangent Cones]
  \label{thm:faces}
For all $\zeta$, the map
  \map{\text{Face}_\zeta}{\cC_\zeta}{\text{Faces of }\pi
    F_\zeta}{S}{\pi F_\zeta \cap(\pi f+\pi Z_0(\clos{S})))} is independent of
  the choice of   $f\in F_\zeta\cap\supp^{-1}(S)$, and induces a bijection
  $$\text{Face}_\zeta\colon \cC^k_\zeta/\!\!\sim \xrightarrow{\;\cong\;} \text{$k$-faces of
    }\pi F_\zeta.$$ Furthermore, for all $[S]\in\cC^k_\zeta/\!\!\sim$,
  $$T_{\text{Face}_\zeta(S)}\pi F_\zeta =\pi F_0(\clos{S}),$$
  where the left-hand side denotes
  the tangent cone to the polyhedron $\pi F_\zeta$ at the face
  $\text{Face}_\zeta(S)$.  In other words, $\pi F_0(\clos{S})$ gives the
  tangent cone corresponding to the configuration $S$ (which is
  independent of the value of $\zeta$) and $\text{Face}_\zeta$ gives the
  corresponding face of $C_\zeta$ (whose direction is also independent 
  of $\zeta$).
\end{thm}
\begin{rmk}
  Note that both the direction of the face corresponding to $S$ and 
the tangent cone of $\pi F_\zeta$ at this face are independent of the 
value of $\zeta$. 
\end{rmk}

In fact, (see Lemma~\ref{lemma:0cycles}) $Z_0(S)$ is generated
by the flows for the cycles supported in $S$, so that
Theorem~\ref{thm:general} corresponds to the  case when $S$ has rank zero.

The concepts needed for the proof of Theorem~\ref{thm:faces} are
developed in Section~\ref{sec:2:flow:basic}.  The theorem is then proved
in Section~\ref{sec:2:flow:classpf} for the classical case (where
$\R^{\cQ_1}=\R^n$ and $\pi\colon\cQ_1\to\cQ_1$ is the identity). The general case follows
easily from this and is treated in Section~\ref{sec:2:flow:genpf}.

Before this, an example of an application and several  important
corollaries are given.

\section{Application to the McKay Quiver}
\label{sec:2:mckay}

\subsection{Example}
\label{sec:2:flow:examples}

\begin{figure}[htbp]
  \begin{center}
    \leavevmode
\epsfysize= 1in \epsfbox{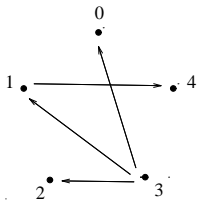}
  \end{center}
  \caption{A configuration $T\in\cC^0$ in the McKay quiver for $\frac{1}{ 5}(1,2,3)$.}
  \label{fig:tree5123}
\end{figure}

\begin{figure}[htbp]
  \begin{center}
    \leavevmode
\epsfysize= 7.5cm \epsfbox{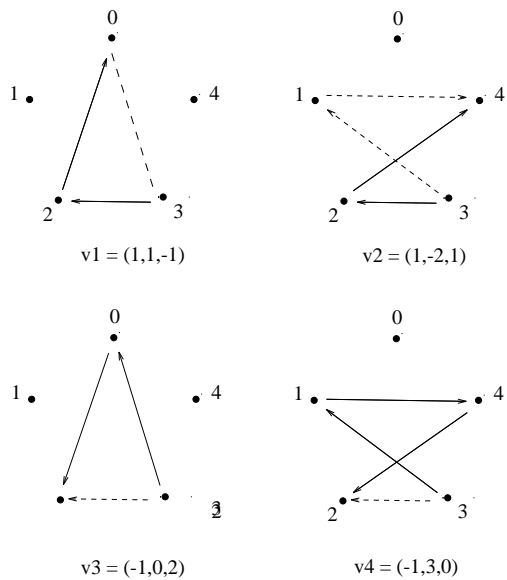}
 \end{center}
  \caption{The different types of cycles which generate the cone $\pi  F_0(T)$ for the configuration $T$ in Figure~\ref{fig:tree5123}.}
  \label{fig:cycles5123}
\end{figure}

Consider the McKay quiver
for the action $\frac{1}{ 5}(1,2,3)$. Recall that  identifying
$\cQ_0$ with $\Z_5$, the arrows are $a_v^i := v\to v-i$ for $v\in
\Z_5$ and $i\in\{1,2,3\}$. Let $\pi $ be the map $\cQ_1 \to\{1,2,3\}$
which assigns to each arrow $a_v^i$ its type $i$. This induces a
projection $\pi\colon \R^{\cQ_1}\to\R^3$ which assigns to the basis element
$\chi_v^i := \chi_{a_v^i}$ the basis element $e_i$ of $\R^3$. With
respect to this map, the cycles of type zero are those with total
number of arrows of any given type equal to zero, where the number of
arrows is counted algebraically according to the orientation of the cycle. Consider
the configuration  $T\subset\cQ_1$ represented in Figure~\ref{fig:tree5123}.
What is the closure of~$T$? A little thought 
shows\footnote{See Section~\ref{sec:2:examples:comm} for more details on this.} that 
$\TC$ is simply $T$ itself. Since $\TC$ has no cycles, $Z_0(\TC)=0$, so 
$T\in\cC^0$. If $T$ is admissible for $\zeta $, and  $f$ is any $\zeta $-flow with support $T$, then  the theorem says 
that $\pi(f)$ is a 
$0$-face of $\pi F_\zeta $, \ie an extreme point. Furthermore, the 
tangent cone to $\pi F_\zeta $ at $\pi(f)$ is the cone $\pi  F_0(T)$; this is 
generated by the types of the cycles whose negative part is contained 
in $T$. These are listed in Figure~\ref{fig:cycles5123}. 

There are four different types\nopagebreak
\begin{align*}
v_1 &= (1,1,-1)\\
v_2 &= (1,-2,1)\\
v_3 &= (-1,0,2)\\
v_4 &= (-1,3,0),
\end{align*}
and they generate $\pi  F_0(T)$ as a cone. Any three of these form a basis 
for the lattice $\Pi\subset\Z^3$ of index $5$ given by
$$\Pi:=\ker\hat\rho = \{(a,b,c): a+2b+3c \equiv 0\pmod 5\},$$
and
they satisfy the single relation $v_1+v_3 = v_2+v_4$.  This
corresponds via the usual toric formalism to a singularity of
$X_\zeta$ of the type $xw=yz \subset\C^4$.  The variety $X_\zeta$
therefore has such a singularity whenever the tree $T$ is admissible
for $\zeta$, \ie whenever the flow $f_\zeta(T)$ is 
positive. Writing this out explicitly, one sees that the admissible
cone for $T$, $\Adm(T)$ 
is given by the following inequalities\footnote{See
  \ref{sec:2:mckay:flow:trees} for the exact definition of $\Adm(T)$.}
\begin{align*}
-\zeta_0 &>0\\
-\zeta_2 &>0\\
-\zeta_4 &>0\\
-\zeta_4-\zeta_1 &>0
\end{align*}

\subsection{Corollaries}
\label{sec:2:flow:cor}

\subsubsection{$k$-Adjacency}

The ``simplex on a graph''
algorithm~\cite[Alg.~3.3]{kenn_helg} can be interpreted as moving from
one extreme point of $F_\zeta$ to an adjacent one by varying the flow
along a single cycle in the network\footnote{See~\cite[\S 5.4]{evans}, or, for
the special case of the ``permutation polytope,''~\cite[\S 5,
Th.~1.3]{yem_kov_kra}.}.  Geometrically this says that two extreme flows
are joined by an edge of $F_\zeta$ if and only if the union of their
supports contains only \emph{one\/} cycle (up to multiples obtained by
going around the cycle $k\in\Z$ times).

In order to generalize the notion of being joined by an edge,  the
following terminology is introduced:    two points of a
polyhedron are called \emph{$k$-adjacent} if they are contained in a
$k$-dimensional face, but in no face of smaller dimension.  Note that
by convexity, two points  are contained in the same face if
and only if their midpoint is also in that face.  Since  a
configuration of the midpoint is obtained by taking the union of configurations
of the two points, the answer to the question is deduced as a
corollary to Theorem~\ref{thm:faces}:
\begin{cor}[$k$-Adjacency]
  \label{cor:k_adjacent} 
  Two points $x,x'\in\pi F_\zeta$ are $k$-adjacent if and only if for
  some (and hence for any) configurations $S$ of $x$ and $S'$ of $x'$,
  $${S\cup S'}\in\cC^k,$$ \ie their union has rank $k$.
\end{cor}
Note that the case~ $k=0$ gives the condition for two configurations
to give rise to the same extreme point. The case $k=1$ says that two
configurations are joined by an edge of~ $\pi F_\zeta $ if and only if
their union only contains cycles whose types are multiples of a fixed
element $v\in\R^n$. For the case when $\pi$ is the identity, this
reduces to the classical statement that the union $S\cup S'$ contains
only one cycle (up to multiples obtained by going around the cycle
$k\in\Z$ times).

\subsubsection{Trees}
\label{sec:2:mckay:flow:trees}

It can be shown easily (see Corollary~\ref{cor:spanning_trees}) that
any extreme point of~ $\pi F_\zeta$ has at least one configuration~
$T$ which is a \emph{tree}, \ie it has no cycles. Thus, to know the
extreme points and the tangent cones  of $\pi
F_\zeta$ at these points, one need only answer the question of which
trees occur as extreme configurations.  Denoting the set of
trees in $\cC$ by $\cT$,  the generalised theorem
implies that the extreme points of $C_\zeta$ are given by the images
of the flows whose supports are members of~ $\cT^0_\zeta$.  This set
is easy to determine: for a given $\zeta$, there is, on each tree $T$, a
unique flow $f_\zeta(T)\in\Lambda^1$ which meets the demand $\zeta$
and is supported on $T$; $f_\zeta(T)$ is called the \emph{($\zeta$-)flow
 admitted by the tree\/} 
$T$. (See \ref{sec:2:examples} for an example). This defines a map
$$f_\zeta\colon\cT\to\Lambda^1.$$  The set $f_\zeta^{-1}(\Lambda^1_+)$ of 
 trees $T$ which admit non-negative $\zeta$-flows is  of course
the set $\cT_\zeta$ of \emph{admissible tree}\emph{s\/} for $\zeta $. 

The classical Theorem~\ref{thm:classical} says that 
$$\ext F_\zeta = f_\zeta(\cT_\zeta).$$  Similarly, the
generalised Theorem~\ref{thm:general} gives the following:

\begin{cor}
  \label{cor:number_extreme_pts}
  Let $f_\zeta$ be the map which assigns to each tree $T$ the unique
  $\zeta$-flow with support equal to $T$.  Then the extreme points of
  $C_\zeta$ are given by
  $$\ext \pi F_\zeta = \pi f_\zeta (\cT^0_\zeta ),$$ where $\cT^0 :=
  \cC^0\cap \cT$ is the set of trees whose closures have type zero.  Furthermore,  Theorem~\ref{thm:faces} implies that
  $$\card{\ext C_\zeta}=\card{\cC^0_\zeta/\!\!\sim}=\card{\cT^0_\zeta/\!\!\sim},$$
where $\sim$ is the equivalence relation induced by closure.
\end{cor}

Note (see Lemma~\ref{lemma:maximal_config} for a proof) that for the
sets $S\in \cC^0$ one has  $\pi  F_0(\SC)=\pi  F_0(S)$.  This implies the
following corollary:

\begin{cor}
  \label{cor:fan}
  The extreme points of the polyhedron $\pi F_\zeta$ correspond to the
  trees $T$ in $\cT^0_\zeta$ and the tangent cone to $\pi F_\zeta$ at
  the point corresponding to $T$ is $\pi  F_0(\TC)=\pi  F_0(T)$. Thus the
  fan\footnote{Recall that the \emph{fan \/}
    associated to a convex polyhedron $P$ is the collection of dual
    cones to the tangent cones of $P$ at all its faces.} $\Sigma
    _\zeta$ associated to the polyhedron $\pi F_\zeta $ is given by
  the dual cones $\pi F_0(T)^\vee$ for the trees $T\in\cT^0_\zeta$ and
  all their faces. In fact, its $k$-skeleton, \ie the set of its
  $k$-dimensional cones, is $\Sigma _\zeta^{(k)} :=\{\pi F_0(\SC)^\vee
  : S\in\cC^{n-k}_\zeta \}$.
\end{cor}
\begin{rmk}
  This corollary says that the singularities of $C_\zeta$ are
  precisely those given (with respect to the lattice $\Pi$) by the
  cones $\{\pi  F_0(T)| T\in\cT^0_\zeta\}$.
\end{rmk}

\subsubsection{Variation of the Flow Polyhedra with $\zeta $}
\label{sec:2:flow:cor:var}

The following corollary of Theorem~\ref{thm:faces} 
describes when two different values of $\zeta$ give isomorphic
polyhedra:

\begin{cor}
  \label{cor:iso_C_zeta}
  If $\zeta$ and $\zeta'$ have the same admissible configurations
  ($\cC_\zeta=\cC_{\zeta'}$) or even just the same admissible trees
  ($\cT_\zeta=\cT_{\zeta'}$) then the corresponding polyhedra
  $C_{\zeta}$ and $C_{\zeta'}$ are geometrically isomorphic.  Two 
  polyhedra are said to be   geometrically isomorphic is they are combinatorially
  isomorphic and their tangent cones at the corresponding faces are
  identical.  In particular, their associated fans are identical, and 
their corresponding toric varieties are isomorphic.
\end{cor}

Note that when one multiplies $\zeta$ by a non-zero number, the
polyhedron $C_\zeta$ is simply scaled-up.  Fixing a tree $T\in\cT^0$ determines
an open convex cone $\Adm(T)\subseteq \Lambda^{0,0}_\R$ of values of
$\zeta$ for which $T$ is $\zeta$-admissible.  The cone $\Adm(T)$ is 
called the
\emph{admissible cone of $T$}.\/ If $\cE$ is a fixed subset of
$\cT^0$, then the condition $\cT^0_\zeta=\cE$ defines an open cone in
the $\zeta$-parameter space $\Lambda^{0,0}_\R$.  As $\cE$ varies, one
obtains a partition of $\Lambda ^{0,0}_\R$ into a union of open cones
inside which the polyhedra $C_\zeta $ are geometrically isomorphic.

\subsubsection{Degeneracies}

Another fact which will be of interest to us is that, for {\em
  generic\/} values of $\zeta$, any $\zeta$-flow has a support which
is a spanning subgraph of the quiver, \ie which connects any two
vertices.  This is written $\cC_\zeta\subset\cC_{\spn
  }$.

For instance, the faces of $\Adm(T)$ consist of degenerate
values of $\zeta$
for which some extreme $\zeta$-flows have supports which are strict
subsets of $T$ and so cannot be spanning subsets.

In the next section, basic flows associated to paths and cycles are
studied in more detail.  They provide the key to the proofs of the
other results.

\section{Proofs}
\label{sec:2:flow:proofs}

\subsection{Basic Flows}
\label{sec:2:flow:basic}

Many proofs in the context of network flows use the basic technique of
decomposing a flow into certain basic components associated to paths
in $\cQ$. A \emph{path} in $\cQ$ means a sequence $p=(p_1,\dots,p_k)$
of arrows in $\cQ_1$ which form a connected path in the underlying
graph to $\cQ$, once their orientation has been disregarded.
Consecutive $p_i$'s are not allowed to be the same (although they can
join the same vertices). As for cycles, the disjoint union of the
arrows of the path $p=(p_1,\dots,p_k)$ which agree\footnote{See
  Section~\ref{sec:2:flow:cycle} for the definition of
  \emph{agree.\/}} (\resp disagree) with the sense of traversal specified by
the sequence $p_1,\dots,p_k$ are called the \emph{positive}
(\resp\emph{negative}) arrows of $p$ and denoted $p^+$ (\resp $p^-$).  A path
will sometimes be confused with its set of arrows, for instance in
statements such as ``a path $p$ is {\em in\/} a set $S\subseteq\cQ_1$
(written $p\subseteq S$)'' which means of course that all its arrows
belong to the set $S$.

To each path $p$ the \emph{basic}\footnote{Basic flows are also
  termed \emph{simple flows}~\cite{busacker_saaty:graphs}
  or \emph{elementary flows} by other authors.}
flow $\tilde\chi_p\colon\cQ_1\to\Z$ is
defined by
\begin{equation}
  \tilde\chi_p(a)=\sum_{a\in p^+}\chi_a - \sum_{a\in p^-}\chi_a\in\Z^n.
    \label{eq:basic_flow}
\end{equation}
Note that if $p$ is a path from $v$ to $v'$, then $\partial
\tilde{\chi}_p=\chi _{v'}-\chi _v$. Conversely, one has the following lemma.
\begin{lemma}
\label{lemma:basicf_path}
If $f$ is an integral flow with $\partial f=\chi _{v'}-\chi _v$ for 
two vertices $v,v'$ of $\cQ$, then there exists a path $p$ from $v$ 
to $v'$ such that $p^{\pm}\subseteq\supp f^\pm$. 
\end{lemma}
\begin{proof}
  By induction on the $1$-norm of $f$:
  $\norm{f}_1:=\sum_{a\in\cQ_1}\abs{f(a)}$.  If $\norm{f}_1=1$ then,
  obviously, $f=\chi _a$ or $-\chi_a$, for some arrow $a\in\cQ_1$
  which is either $v\to v'$ or $v'\to v$ respectively. Either way,
  $f=\tilde{\chi }_{p}$ for the corresponding one-arrow path $p$ from
  $v$ to $v'$. Now suppose that $\norm{f}_1>1$ and that $\partial
  f=\chi _{v'}-\chi _v$. There must be an arrow $a$ such that one of
  the following statements holds
\begin{enumerate}
\item[(1)] $t(a)=v$ and $f(a)>0$, or
\item[(2)] $h(a)=v$ and $f(a)<0$.
\end{enumerate}
If (1) holds then the flow $f'=f-\chi _a$ satisfies the induction
hypothesis with $\partial f'=\chi _{v'}-\chi _{h(a)}$, so there exists
a path $p'\subseteq\supp f'$ from $h(a)$ to $v'$. But then $p=ap'$ is
a path from $v$ to $v'$ with $p^+ = p^{\prime +}\cup
\{a\}\subseteq\supp f^{\prime+}\cup\{a\} = \supp f^+$. Case (2)
follows in a similar way setting $f'=f+\chi _a$ and $p=(-a)p$.
\end{proof}

\begin{figure}[htbp]
  \begin{center}
    \leavevmode
\epsfysize= 3cm \epsfbox{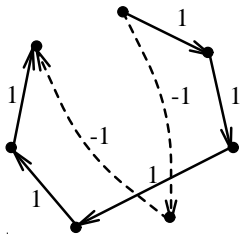}
  \end{center}
  \caption{A basic $0$-flow associated to the 
  cycle in Figure~\protect\ref{fig:uncycle7123}}
  \label{fig:basicf7123}
\end{figure}

Obviously, if $p$ is actually a cycle, then the basic flow associated
to $p$ is a $0$-flow (it satisfies $\partial \tilde\chi_p = 0$), and
is supported in $p$.   Figure~\ref{fig:basicf7123}
gives an example of such a flow.

The converse to this statement is given by the following lemma~\cite[Th.~7.2]{busacker_saaty:graphs}, \cite[{\S5.1.6,Th.~3}]{gond_mino:graphs}, .
\begin{lemma}[Decomposition into basic flows]
  Any\/ $0$-flow $f$ can be decomposed into a positive linear
  combination of basic flows for finitely many cycles $c_i$ such that
  $c^\pm_i\subseteq \supp f^\pm$.
\label{lemma:0cycles} 
\end{lemma} 
\begin{proof} 
  Let $f\in\ker\partial$ and $S$ be its support.
  If $S$ contains a cycle $c$ then let $v$ be a vertex such that $c^+$
  has an arrow $a$ with $t(a)=v$. Then $f':=
  f-f(a)\tilde\chi_c\in\ker\partial$ and $\supp f' = \supp f \setminus
  \{a\}$ is strictly smaller than $\supp f$. Repeating this reasoning
  finitely many times, one obtains a flow $f'\in\ker\partial$ such that
  $S'=\supp f'$ contains no cycles.  Suppose $S'$ is non-empty, and
  let $p$ be a path in $S'$ which is \emph{maximal, \/} \ie is not
  contained in any strictly larger path.  Let $p$  go from vertex $v$ to
  vertex $v'$ and $a$ denotes the first arrow of $p$ (the one joined
  to $v$).  Then this is the only arrow in the quiver which is joined to
  $v$ and on which $f'$ is non-zero, so $\partial f'(v)$ must be non-zero.
  Since this is impossible, $S'$ must be empty, and $f=\sum_i
  x_i \tilde\chi_{c_i}$ for some cycles $c_i\subseteq S$ and real
  numbers $x_i$. Replacing $c_i$ by $-c_i$ if necessary, one may assume
  $x_i>0$. (Note also that if $f$ is integral to start with, the $x_i$
  are also integral.)

  I claim that these cycles may be chosen in such a way that their
  orientations are \emph{conformal}, \ie such that any
  two cycles $c_i,c_j$ satisfy $c^+_i\cap c^-_j=\emptyset$. Indeed
  suppose that $c_i$ and $c_j$ are not conformal, and suppose, say,
  $x_i\geq x_j$. 
  On the subset $U\subseteq\cQ_1$ on which their orientations
  disagree, the sum of $x_ic_i$ and $x_jc_j$  will cancel each other out, and only $x_i-x_j$ units
  of flow will survive. The complement $c_i\cup c_j\setminus U$ will
  consist of one or more disjoint cycles $d_1,\dots,d_k$, all
  conformal to $c_i$.  Thus upon adding the two basic flows 
$$x_i\tilde\chi_{c_i}+x_j\tilde\chi_{c_j} = 
(x_i-x_j)\tilde\chi_{c_i} + x_j(\tilde{\chi }_{d_1}+\dots+\tilde{\chi 
}_{d_k}),$$ giving $k+1$ conformal cycles with all coefficients
positive or zero. Repeating this procedure for 
all the pairs of non-conformal cycles gives the desired conformal 
decomposition. 

  Now writing $$f=\sum_i\tilde\chi_{c_i}=\sum_i
  (\chi_{c^+_i}-\chi_{c^-_i}),$$ and since $c^+_i\cap
  c^-_j=\emptyset$, one has
  $$f^\pm = \sum_i\chi_{c^\pm_i},$$
  which is the desired result.
\end{proof} 
\begin{rmk}
 The decomposition lemma above can be viewed as a consequence of the 
following homological argument. Regard $\cQ$ as a CW-complex, and for 
each cycle $p$ in $\cQ$, adjoin a 2-cell whose boundary is the element 
$\tilde\chi_{p}$. Then the resulting CW-complex $\widetilde\cQ$ has 
$H^1(\widetilde\cQ)=0$. 
An  elementary proof was given because it  illustrates the basic
technique of decomposing a flow into basic flows. 
\end{rmk}

\subsection{Proof of Theorem~\protect\ref{thm:faces} --- Classical Case}
\label{sec:2:flow:classpf}

The classical analog of Theorem~\ref{thm:faces} (\ie the case when 
$\pi$ is the identity map) can now be proved.  The statement says that 
the $k$-dimensional faces of $F_\zeta$ correspond to the $\zeta 
$-configurations $S$ such that $\rk Z_0(S)=k$, \ie whose cycles 
span a lattice of rank $k$.  This theorem follows from the following 
four facts:
\begin{fact}
\label{fact:1}
If $f\in F_\zeta$ then the cone $\R_+(F_\zeta-f)$ contains a subspace
of dimension $k$ if and only if $f$ is contained in the (relative)
interior of a face of dimension $k$ (and hence in no lower dimensional
face).
\end{fact} 
\begin{fact}
  \label{fact:2}
For any $f\in F_\zeta$, one has $\R_+(F_\zeta-f)=  F_0(\supp f)$.
\end{fact}
\begin{fact}
  \label{fact:3}
The maximal vector subspace in the cone $ F_0(S)$ is $Z_0(S)$.
\end{fact}
\begin{fact}
  \label{fact:4}
The subspace $Z_0(S)$ is generated by the basic flows for the
  cycles supported in $S$.
\end{fact}

Fact~\ref{fact:1} follows from the definition of a $k$-dimensional
face, \ref{fact:3} is trivial and~\ref{fact:4} follows from Lemma~\ref{lemma:0cycles}. Fact~\ref{fact:2} is proved in the
following lemma.
\begin{lemma}
\label{lemma:tangent_fzeta}
  If $f\in F_\zeta$, then  $\R_+(F_\zeta -f) =  F_0(\supp f)$.
\end{lemma}
\begin{proof}
  Let $f\in F_\zeta$, so that
  $F_\zeta=(f+\ker\partial)\cap\R_+^{\cQ_1}$. The non-negative
  multiples of elements in $(f+\ker\partial)\cap\R_+^{\cQ_1} -f$
  belong to $\ker\partial$ and are non-negative outside $\supp f$, so
  belong to $ F_0(\supp f)$. Conversely, suppose $m\in  F_0(\supp f)$.
  Then $\partial(f+\epsilon m)=\zeta $ for all $\epsilon\in\R$, and it
  suffices to show that there exists $\epsilon >0$ such that
  $f+\epsilon m\in \R_+^{\cQ_1}$. But this follows because $f$ is
  bounded below by a positive number on $\supp f$, whereas $m\geq 0$
  outside $\supp f$.
\end{proof}

\subsection{Proof of Theorem~\protect\ref{thm:faces} --- General Case}
\label{sec:2:flow:genpf}

In order to prove the general case of Theorem~\ref{thm:faces},  a bit 
more has to be said about the various configurations which can occur for
a point $x\in\pi F_\zeta$.
\begin{lemma}
\label{lemma:maximal_config}
If $S$ is a $\zeta $-configuration for $x$ then so is $\SC$.
Furthermore, all $\zeta $-configurations for $x$ have the same
closure.
\end{lemma}
\begin{proof}  
  Let $f\in F_\zeta$ be a flow such that $\pi f=x$ and let $S=\supp
  f$.  The proof begins by constructing a flow $\clos f\in F_\zeta\cap
  \pi^{-1}(x)$ whose support is $\SC$. If $c$ is a cycle of type $0$
  in $\cQ$ such that $c^-\subseteq S$ and $c^+\not\subseteq S$, one
  can add a small positive multiple of $\tilde\chi_c$ to $f$ and
  obtain a non-negative flow $f'$.  Since $c$ is a cycle,
  $\partial\tilde\chi_c=0$, and since $c$  has type $0$, $\pi f=\pi
  f'$ and hence $f'\in F_\zeta$. Continuing in
  this way until all cycles of type $0$
  satisfy~\eqref{eq:cycle_closure}, one obtains the required flow
  $\clos f$.

For the second statement of the lemma, note that if $f,f'$ are two 
elements of $F_\zeta\cap\pi^{-1}(x) $, then $\partial(f-f')=0$, so by 
Lemma~\ref{lemma:0cycles}, has a decomposition into basic flows for 
cycles $c_i$: $$f-f' = \sum_i x_i\tilde\chi_{c_i},$$ with 
$c_i^+\subseteq f$, $c_i^-\subseteq f'$ and $x_i>0$. Now 
$c_i^+\subseteq S\implies c_i^-\subseteq \clos S$, so $$f'=f+\sum_i 
\tilde\chi_{c_i}$$ implies that $S'\subseteq \clos S$. By symmetry, one 
also has $S\subseteq \clos{S'}$, and so $\clos S = \clos{S'}$. 
\end{proof}
A little corollary  needed later is
\begin{cor}
\label{cor:spanning_trees}
Any $S\in\cC^0$  contains a tree $T$ such that $\TC=\SC$. 
\end{cor}
\begin{proof}
  Let $S=\supp f$ be a configuration for $x$. Eliminate any cycles in
  the support of $f$ by adding basic flows corresponding to those
  cycles, as in the proof of Lemma~\ref{lemma:maximal_config}.  The
  resulting flow $f'$ is supported in a tree $T$ and, since all the
  basic flows have type zero it satisfies $\pi(f')=x$.
  Lemma~\ref{lemma:maximal_config} shows that $\TC=\SC$.
\end{proof}

\begin{proof}[Proof of the general Theorem~\ref{thm:faces}]
  Consider a point $x\in\pi F_\zeta$. The cone $\R_+(\pi F_\zeta-x)$
  gives the tangent cone to $\pi F_\zeta$ at the minimal face of
  $F_\zeta$ containing $x$. It is obtained by taking the union over
  all $f\in F_\zeta\cap\pi^{-1}(x)$ of $\pi \R_+(F_\zeta-f)$. By
  Lemma~\ref{lemma:tangent_fzeta} this gives
  $$\bigcup\{\pi  F_0(S) : S \text{ a $\zeta $-configuration for }x\},$$
  which by Lemma~\ref{lemma:maximal_config} is $\sigma(S)=\pi
   F_0(\SC)$, for any $S=\supp f$ and $f\in F_\zeta\cap\pi^{-1}(x)$.
  The minimal face of $\pi F_\zeta$ containing $x$ is given by
  intersecting $\pi Z_0(\SC)$, the largest subspace in the cone $\pi
   F_0(\SC)$, with the tangent cone of $\pi F_\zeta $ at $x$ and then
  translating by $x$.  This gives
  $$x+(\pi F_\zeta )_x \cap Z_0(\SC),$$ which is precisely the face
  $\text{Face}_\zeta(S)$ mentioned in the theorem and all faces of
  $\pi F_\zeta$ are obtained in this way.  It is obvious that
  equivalent configurations have the same rank and give the same face,
  so the partition of $\cC$ has the stated properties and
  $\text{Face}_\zeta$ is a bijection $\cC^k_\zeta/\!\!\sim\to k\text{-faces}$.
\end{proof}

%\section{-}
\section{Exactness Results}
\label{sec:2:exact}

The exactness of the sequence~\eqref{eq:exact1} is proven in this
section.

\subsection{Basic Flows: Sequential Notation}
\label{sec:2:exact:sequ}
\aftersub

We introduce some notation which is convenient to describe basic
flows. This is used in~ Section~\ref{sec:2:exact:exact} to give an
elementary proof of the exactness of~\eqref{eq:exact1}.

For $v\in\cQ_0$ and $j\in\{1,\dots,n\}$, write $\chi_v^j$ for the 
basis element $\Lambda ^1$ which is the indicator function of the 
singleton $\{a_v^j\}\subset\cQ_1$. 
Define the following 
symbols:
\begin{align}
  \{v\}(j) &:= \chi_v^j\\
\{v\}(-j) &:= -\chi_{v+w_j}^j.
\end{align}
For $k>0$, and  $j_0,\dots,j_k\in\pm\{1,\dots,n\}$, define
$$ \{v\}(j_0,\dots,j_k):=\{v\}(j_0)+\{v-w_{j_0}\}(j_1,\dots,j_k).$$
Also define $(j)\{v\} := \{v+w_j\}(j)$ and
$$(j_k,\dots,j_0)\{v\}=(j_k,\dots,j_1)\{v+w_{j_0}\} + 
(j_0)\{v\}.$$

\begin{figure}[htbp]
  \begin{center}
    \leavevmode
\epsfysize= 3cm \epsfbox{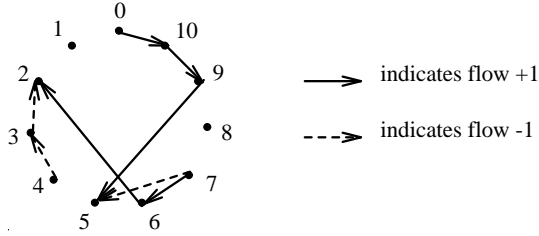}
  \end{center}
  \caption{The basic flow $\{0\}(1,1,3,-2,1,3,-1,-1)$ in the
McKay quiver for $\frac{1}{ 11}(1,2,4)$.}
  \label{fig:bfl11124}
\end{figure}

This sequential notation is designed especially for representing basic
flows and has the advantage of including only the relevant
information.  For instance, if $p=(p_1,\dots,p_k)$ is a 
path in $\cQ$, then the
basic flow associated to $p$ is, in this notation,
$$\tilde\chi_p=\{t(p_1)\}(j_1,\dots,j_k),$$ where 
$$j_i:=
\begin{cases}
\phantom{-}\pi(p_i),\qquad & p_i\in p^+\\
-\pi(p_i),\qquad & p_i\in p^-
\end{cases}
$$
For example, the basic flow represented in 
Figure~\ref{fig:basicf7123} 
can be written as $$\{0\}(1,1,2,1,1,-3,-3)$$ in this notation.  
Figure~\ref{fig:bfl11124} shows yet another example.

Often, when there is a need to specify both endpoints, the notation  
$\{v\}(j_0,\dots,j_k)=:\{v\}(j_0,\dots,j_k)\{v'\}$ with 
$v'=v-\sum w_{j_i}$ will be adopted. The 
following identities are easily checked for any $v\in\cQ_0$ 
and 
$j,j_i\in\pm\{1,\dots,n\}$: 
\begin{gather}
(j)\{v\}+\{v\}(-j) = 0 \label{eq:v-id1},\\
\{v\}(j_0,\dots,j_k,-j_k,\dots,-j_0)=0 \label{eq:v-id2},\\
\{v\}(j_0,\dots,j_k)\{v'\}+\{v'\}(j_{k+1},\dots,j_l)=
\{v\}(j_0,\dots,j_l), 
\label{eq:v-id3}
\end{gather}
where the last identity holds of course for $v'=v-\sum_{i=0}^k 
w_{j_i}$.

If $\sum_i w_{j_i}=0$, \ie the corresponding flow 
$\{v\}(j_1,\dots, j_k)$ is a $\partial$-closed flow associated to a cycle,  
there is another identity which allows us to cyclically permute 
the indices: 
\begin{equation}
\label{eq:v-id4}
\{v\}(j_0,\dots,j_k)=\{v-w_{j_0}\}(j_1,\dots,j_k,j_0).
\end{equation}
Recall that $\Lambda^2$ is the sub-lattice of $\Lambda^1$ generated by
the elements corresponding to the commutation relations.  These can be
written in various notations: $$r^{ij}_v=\chi_v^i
+\chi_{v-w_i}^j-\chi_v^j-\chi_{v-w_j}^i=
%\tilde\chi_{c^{ij}_v}=
\{v\}(i,j,-i,-j).$$

Putting these identities together gives the following lemma.
\begin{lemma}
  For any permutation $\sigma$ of $\{0,\dots,k\}$ one has 
  $$\{v\}(j_{\sigma(0)},\dots,j_{\sigma(k)}) =
  \{v\}(j_0,\dots,j_k)\mod \Lambda^2.$$
\label{lemma:permutation}
\end{lemma}
\vspace*{-\belowdisplayskip}
\vspace*{-\topsep}
\vspace*{-\partopsep}
\begin{proof}
It is enough to show that the equality holds for any 
transposition $\sigma=(q,q+1)$ of consecutive elements.
  Let $j_q=i$ and $j_{q+1}=l$. By the identities 
\eqref{eq:v-id2}--\eqref{eq:v-id3}
one has
$$
\{v\}(j_0,\dots,j_k) = 
\{v\}(j_0,\dots,j_{q-1})+\{v_q\}(i,l)+
\{v_{q+2}\}(j_{q+2},\dots,j_k),
$$ where  $v_a:=v-\sum_{i=0}^{a-1} w_{j_i}$ for any $a\in\{0,\dots, 
k\}$. Using 
$\{v_q\}(l,i,-l,-i)+\{v_q\}(i,l)=\{v_q\}(l,i)$, one sees that 
$$\{v\}(j_{\sigma(0)},\dots,j_{\sigma(k)})-\{v\}(j_0,\dots,j_k)
=\{v_q\}(l,i,-l,-i).$$ Now if $i,l\in\{1,\dots,n\} $ then the result 
follows because $\{v_q\}(l,i,-l,-i)=r^{il}_{v_q}$. If on the other hand
$-i,l\in\{1,\dots,n\} $, then the result is true because 
\begin{align*}
\{v_q\}(l,i,-l,-i) &=\{v_q +w_i\}(-i,l,i,-l),\qquad\text{ (by 
equation~\ref{eq:v-id4})}\\
&=r^{-i,l}_{v_q +w_i}.
\end{align*}
The other two possibilities ($(i,-l)$ and $(-i,-l)$) also follow in this way.
\end{proof}

Another useful result follows from these identities and
Corollary~\ref{lemma:0cycles}:

\begin{lemma}
  Any $\partial$-closed flow can be written as the basic flow
  associated to a single cycle in $\cQ$ (not necessarily contained in
  its support).
\label{lemma:null_single_cycle}
\end{lemma}
\begin{proof}
  Note that in the proof of Lemma~\ref{lemma:0cycles}, 
  $f\in\ker\partial$ was decomposed into a sum of basic flows for 
  cycles $c$.  The sum of two basic flows corresponding to cycles can 
  be written, using identities \eqref{eq:v-id2}--\eqref{eq:v-id4} as
\begin{multline}
\{v\}(j_1,\dots,j_k)+\{v'\}(j'_1,\dots,j'_{k'}) = \\
\{v\}(j_1,\dots,j_k,a_1,\dots,a_l,j'_1,\dots,j'_{k'},-a_l,\dots,-a_1),
\label{eq:v-sum}
\end{multline}
where $a_i\in\pm\{1,\dots,n\}$ are such that $v-\sum_i 
w_{a_i}=v'$, 
and one has $\sum_i w_{j_i}=\sum_i w_{j'_i}=0$ since the basic 
flows on the left-hand side correspond to cycles. This proves the 
lemma. 
\end{proof}

\subsection{Exactness Results}
\label{sec:2:exact:exact}
\aftersub 

The promised proof of the exactness of~\eqref{eq:exact1} can now be 
given. 
Recall the morphism of lattices $\pi\times\partial\colon\Lambda^1\to 
\Z^n\times\Lambda^{0,0}$ defined in 
section~\ref{sec:2:abel:toric}.
\begin{lemma}
 One has $\ker \pi\times\partial =\Lambda^2$, \ie 
$$\Lambda^2\to\Lambda^1\stackrel{\pi\times\partial}\to\Z^n\times\Lambda^{0,0}$$
is exact.
\label{lemma:kerc_lambda2}  
\end{lemma}
\begin{proof} Since $\pi\times\partial(r^{ij}_v)=0$ one only has to prove that 
  $\ker \pi\times\partial\subseteq\Lambda^2$. Let $f\in\ker \pi\times\partial$, and suppose that $c'$
  is a cycle such that $f=\tilde\chi_{c'}$, as
  described in Lemma~\ref{lemma:null_single_cycle}. Then $f$ is of the
  form
  $$f=\{v\}(j_0,\dots,j_k),$$ for some $v\in\cQ_0$ and
  $j_i\in\pm\{1,\dots,n\}$. Since $\pi f=0$, one has $\#\{i:
  j_i=j\}=\#\{i:j_i=-j\}$ for any $j\in\{1,\dots,n\}$.  Thus,
  up to permutations, $(j_0,\dots,j_k)$ can be rewritten as
  $(j_1,-j_1,j_2,-j_2,\dots,j_l,-j_l)$
  for some elements $j_k\in\{1,\dots,n\}$. Since
  $\{v\}(j_1,-j_1,j_2,-j_2,\dots,j_l,-j_l)=0$, the result
  follows by Lemma~\ref{lemma:permutation}.
\end{proof}

From this lemma it follows that~\eqref{eq:exact1} is exact, and 
induces an inclusion $\Lambda\hookrightarrow \Z^n\times\Lambda^{0,0}$.  
Recall that $\Pi=\ker\hat\rho$, where 
\map{\hat\rho}{\Z^n}{\gahat}{x}{\sum_i x_iw_i.}

\begin{lemma}
\label{lemma:Lambda_n_pi}

One has $\Pi=\pi(\ker\partial)$.
\end{lemma}
\begin{proof} 
 Suppose $x\in\Pi $.  One wants to find $f\in\ker\partial$ 
such that, for all $i\in\{1,\dots,n\}$, $f$ satisfies 
$$\sum_{a:\pi(a)=i} f(a) = x_i.$$
 This is easy: just take the basic flow given by 
\begin{equation} 
f=\{v\}(\underbrace{ 1,\dots, 1}_{x_1}, 
\underbrace{ 2,\dots, 2}_{x_2},\dots, 
\underbrace{ n,\dots, n}_{x_n}), 
\label{eq:cycle_lambda_n} 
\end{equation} for any $v\in\cQ_0$.  (If $x_i$ is negative under any 
brace, the notation is taken to mean $-x_i$ copies of $-i$.)  Now $\sum_{i=1}^n w_i{x_i} = 
0\text{ in }\gahat$ is equivalent to the fact that the basic flow 
corresponds to a cycle, and so $\partial f=0$.  The converse follows 
from Lemma~\ref{lemma:null_single_cycle} and the preceding sentence.
\end{proof}

\begin{lemma}
\label{lemma:imagec}
   There is an exact sequence of abelian groups
   $$0\to\Lambda \xrightarrow{\pi\times\partial}
   \Z^n\times\Lambda^{0,0} \xrightarrow{\hat\rho-\hat\nu} \Z^n/\Pi
   \cong\gahat \to 0,$$
   where
   \map{\hat\nu}{\Lambda^{0,0}}{\gahat}{\zeta}{\sum_{v\in\cQ_0}\zeta(v)v}
   is the morphism of lattices dual to the action of \/$\Gamma$ on
   $\End R$ by conjugation.
\end{lemma}
\begin{proof}
  One  needs to show that if $(x,\zeta)\in\Z^n\times\Lambda^{0,0}$ satisfies
  \begin{equation}
    \sum_{i=1}^n x_iw_i - \sum_v \zeta(v)v=0,
    \label{eq:zerosum}
  \end{equation}
then there exists a flow
  $f\in\Lambda^1$ such that $\pi(f)=x$ and $\partial f=\zeta$.  One
  can construct this flow in two steps.  For convenience, suppose that
  the weights of the action of $\ga$ on $Q$ have been normalised so
  that $w_1=1$.  First construct a flow $g$ such that $\partial
  g=\zeta$ with only arrows of type $1$: start at vertex $0$ and let
  $g(1\to 0)=\zeta(0)$.  Next, let $g(2\to 1)=\zeta(1)+\zeta(0)$, and
  continue in this way, setting
  $$g(k+1\to k) = \sum_{j=0}^k\zeta(j),$$ and $g=0$ on all the other
  arrows.  Then $\partial g =\zeta$ by construction, and
  $\pi(g)=X_1e_1$, where $X_1=\sum_j\zeta(j)j$.  Now add on any flow of the form
  $$g'=\{v\}(\underbrace{ 1,\dots, 1}_{x_1-X_1}, \underbrace{ 2,\dots,
    2}_{x_2},\dots, \underbrace{ n,\dots, n}_{x_n}),$$ using the same
  conventions as in the previous lemma for the case when the integers
  under the braces are negative.  Equation~\eqref{eq:zerosum} means
  that $g'$ is the basic flow associated to a cycle, and so $\partial
  g'=0$.  The flow $f=g+g'$ has the required properties.
\end{proof}

%\section{-}
\section{Singular Configurations}
\label{sec:2:sing}

In this  section  some comments are made regarding the singularities of $C_\zeta$.

Recall that the tangent cone to $C_\zeta =\pi F_\zeta $ at a point $x$
which has a configuration $S\in\cC^k$, is the cone $\pi  F_0(\SC)$.
This is singular if its dual $(\pi  F_0(\SC))^\vee$ is
generated by a part of a basis of ${\Pi}^*$.  In this case, 
$S$ is called a \emph{singular configuration}.\/ The set of
singular configurations is denoted $\cS$.

\begin{prop} 
The following statements are true.
  \begin{enumerate}
  \item $C_\zeta$ is singular at the faces corresponding to
    configurations $S\in \cS_\zeta$.
  \item $C_\zeta$ is non-singular in co-dimension $k$ if and only if
    $\cS_\zeta^{n-k}=\emptyset$.
  \item $C_\zeta$ is generically non-singular in co-dimension $k$ if
    and only if
    $\cS_{\spn}^{n-k}=\emptyset$.
  \item $C_\zeta$ is generically non-singular if and only if
    $\cS_{\spn}=\emptyset$.
  \end{enumerate}
\end{prop}
\begin{proof}
  Part (1) of the following proposition follows from
Theorem~\ref{thm:faces} and Lemma~\ref{lemma:generic_flow}. Part~(2) follows because $S\in\cS^{n-k} \iff \dim \sigma_S = k$. The last two
statements follow because $\cS_{\spn}=\cup_{\zeta\text{ generic }}\cS_\zeta$.
\end{proof}

\begin{question}
  What is the lowest $k$ for which $\cS^k_\zeta$ is empty (\ie in
  what co-dimension is $C_\zeta$ smooth)? The cone $C_0$ is smooth in
  co-dimension $1$ (it has an isolated singularity), so it seems likely
  that $\cS^k=\emptyset$ for $k\geq1$.  Of course, this is equivalent to
  the statement that $\cS^1=\emptyset$.
\end{question}
The statements about singular trees  are translated here  for the record.

\begin{conj}
  The polyhedra $C_\zeta$ are non-singular in co-dimension $n-1$,
  \ie $\cS^1=\emptyset$. If\/ $\ga\subset\SU(3)$, then $C_\zeta$ are
  non-singular for generic values of $\zeta$, \ie 
  $\cS^0_{\text{span}}=\emptyset$.  Furthermore, the Euler number of
  $C_\zeta$ for generic $\zeta$ is equal to the order of $\Gamma$,
  i.e\ $\card{\cT^0_\zeta/\!\!\sim}=\card{\ga}$.
\end{conj}

Let us look at the case of singular \emph{points.\/} Suppose that
$S\in\cC^0$ is a singular configuration. This means that the primitive
generators $\rho(S)=(\rho^1_S,\dots,\rho^k_S)$ of $\pi  F_0(S)$ do not
form a basis of $\Pi $. If $k=n$, $\pi  F_0(\SC)$ corresponds to a
finite abelian quotient singularity, whereas, if $k>n$, the
singularity is determined by the linear relations holding between the
$\rho_S^i$.  To determine whether a given set $S\in\cC^0$ is singular
or not, one must find all the cycles $c$ such that $c^-\subseteq S$,
calculate their type, and see what primitive vectors of $\Pi$
one obtains.  It is sufficient to restrict one's attention to the
cycles which are not decomposable into a union of cycles.  Let us look
at some examples.

\section{Examples and Computations}
\label{sec:2:examples}

\subsection{Commutators}
\label{sec:2:examples:comm}

In this section, the commutator of a configuration $S$ is defined;
this is specific to the McKay quiver and its purpose is purely
computational: it gives a necessary criterion for determining when
$S\in\cC^0$ which is extremely useful in practical calculations.

Let us begin by working out what the closure of a subset
$S\subseteq\cQ_1$ corresponds to in the case of the McKay quiver for
the action $\frac{1}{ r}(w_1,\dots,w_n)$. Let $\pi$ be the map
$\cQ_1\to\{1,\dots,n\}$ which assigns to each arrow $a_v^i = v\to
v-w_i$ its type $i$. This induces a map $\pi\colon\R^{\cQ_1}\to\R^n$
which assigns to the basis element $\chi_v^i=\chi_{a_v^i}$ the basis
element $e_i$ of $\R^n$.

\begin{figure}[htbp]
  \begin{center}
    \leavevmode
\epsfysize= 3cm \epsfbox{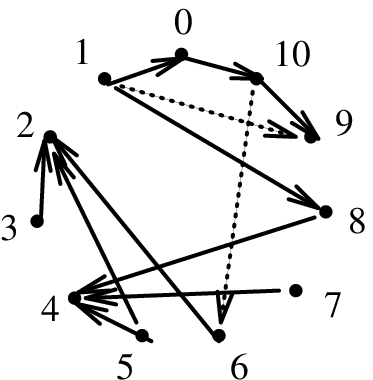}    
 \end{center}
  \caption{A configuration $S$ for which $S=S^{\bowtie}\neq \SC$. The dotted arrows indicate the arrows in $\SC\setminus S^{\bowtie}$}
  \label{fig:notcc}
\end{figure}

Recall that to calculate the closure of $S$, one must find the 
smallest over-set of $S$ which contains the positive part of cycles of 
type zero if and only if it contains the negative part.  The simplest 
cycles which have type zero are cycles with only four arrows: starting 
from vertex $v$, go forward along arrow $a_v^i$ to $v-w_i$, forward 
again along $a_{v-w_i}^j$ to $v-w_i-w_j$, back along $a_{v-w_j}^i$ to 
$v-w_j$ and back along $a_v^j$ to $v$.  This cycle is denoted by 
$c_v^{ij}$.  The basic flow corresponding to this $c^{ij}_v$ is 
$\{v\}(i,j,-i,-j)$ in the notation of Section~\ref{sec:2:exact:sequ}.  If 
$S$ contains $c^{ij+}_v=\{a_v^i,a_{v-w_i}^j\}$ for some $v,i,j$ the 
closure of $S$ must contain $c^{ij-}_v=\{a_v^j,a_{v-w_j}^i\}$.  The 
two pairs of arrows will be called \emph{complementary pairs}.\/ 
\begin{dfn} 
Denote by $p^{ij}_v$ the pair of arrows $\{a_v^i,a_{v-w_i}^j\}$.  
If $S$ is a subset of $\cQ_1$,  the 
smallest over-set $S^{\bowtie}\supseteq S$ satisfying the  
``commutation condition'' 
\begin{equation} p^{ij}_v \subseteq 
S^{\bowtie} \iff p^{ji}_v\subseteq S^{\bowtie}
\label{eq:S_commute}
\end{equation} is  called the \emph{commutator} of $S$.
\end{dfn} 
\begin{figure}[htbp]
 \begin{center} \leavevmode 
\epsfysize= 3cm \epsfbox{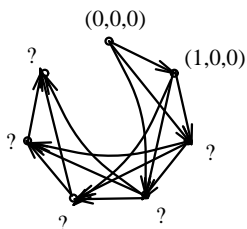} 
\end{center} 
\caption[The commutator of the set in 
Figure~\protect\ref{fig:uncycle7123}]{The commutator of the set in 
  Figure~\protect\ref{fig:uncycle7123}. Notice that this is not an
  invariant set: there is no way to assign elements of $\Z^n$ to the
  vertices in such a way that equation~\protect\ref{eq:S_inv} is
  satisfied.} \label{fig:commutator}
\end{figure}
 The fact  that $\SC\supseteq S^{\bowtie}$ is
extremely useful in practical computations, because it gives the
following easily checked necessary condition for $S$ to be an IC-set.
\begin{lemma} 
If  $S$ is a  spanning IC-set then the  morphism 
$\cW_S\colon \cQ_0 \to\Z^n$ defined by 
equations~\ref{eq:lambda_def} satisfies the following 
conditions for all $v\in\cQ_0$, and $i,j\in\{1,\dots,n\} $:
\begin{equation}
  p^{ij}_v \subseteq S \implies 
  \begin{cases}
  \cW_S({v-w_j})-\cW_S({v})&=e_j \\
 \cW_S({v-w_j-w_i})-\cW_S({v-w_j})&=e_i,
  \end{cases}
  \label{eq:lambda_commute}
\end{equation}
where $\{e_i\}$ denotes the standard basis of $\Z^n$.
\end{lemma}
It is rather surprising that this condition is in fact \emph{not\/} sufficient.  Figure~\ref{fig:notcc} gives an example of a set $S$ for which
$S=S^{\bowtie}\neq \SC$.

\subsection{Weightings and Stabilisers}
\label{sec:2:examples:weight}

There is   another way of understanding the condition $S\in \cC^k$ for a 
set $S\subseteq\cQ_1$ which relates directly to the toric geometry of 
$X_\zeta =\T^{\Pi,C_\zeta }$. The basic idea is that a $k$-dimensional 
face of $C_\zeta $ corresponds to elements of $X_\zeta $ which are 
fixed by a torus of codimension $k$. For instance, the extreme points 
correspond to fixed points of $\Tn$. Let us begin with this case for 
simplicity.

\subsubsection{Fixed Points and $n$-weightings}
\label{sec:2:examples:weight:fixed}
Recall that the lattice $\Z^{\cQ_0}$ is denoted by $\Lambda ^0$. The 
sub-lattice of co-rank $1$ defined by the equation $\sum_{v\in\cQ_0} 
\zeta(v) = 0$ is denoted $\Lambda ^{0,0}$. 

For any subset of arrows $S\subseteq\cQ_1$ one attempts to find 
a morphism $\cW_S\colon \Lambda^{0,0}\to\Z^n$ 
satisfying 
\begin{equation}
\cW_S(\partial \chi_a)=\pi(\chi_a),\quad\text{ for }a\in S.
  \label{eq:S_inv}
\end{equation}
If $\cW_S$ exists, it satisfies 
\begin{equation}
\cW_S(\partial\tilde{\chi }_p)=\pi\tilde{\chi }_p  
\label{eq:lambda_def}
\end{equation} for any path $p=(p^1,\dots,p^k)$ in $S$.  In 
particular, this equation implies that $\pi\tilde{\chi }_c=0$ must 
hold for all cycles $c\subseteq S$.  In other words, if $\cW_S$ 
exists, then $\rk \pi Z_0(S)=0$.  Conversely, if all cycles in $S$ 
have zero type then one can find a morphism $\cW_S$ 
satisfying~\eqref{eq:S_inv}.  Note that this is possible if and only 
if there is morphism $\cW'_S\colon \Lambda^0\to\Pi$ which extends 
it.  The latter corresponds to a labeling $$\cW'_S\colon\cQ_0\to\Z^n$$ of the vertices of the quiver by elements of 
$\Z^n$ in such a way that equation~\ref{eq:S_inv} holds.  This is 
called an \emph{$n$-weighting} of $S$, and 
$S$ is called~\emph{invariant} if it admits such a weighting.
\begin{figure}[htbp]
\begin{center} 
\leavevmode 
\epsfysize= 3cm \epsfbox{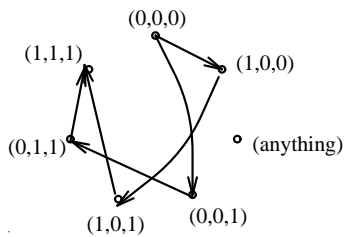} 
\end{center} 
\caption{Example of an $n$-weighting of a configuration of arrows.} 
\label{fig:lambdas7123} 
\end{figure}
In summary, $S$ is an extreme configuration if and only if $\SC$ is 
invariant, if and only if $\pi Z_0(S)=0$.  Such $S$ are called \emph{
invariant closure configurations\/} or \emph{IC-configuration}\emph{s\/} 
for short.  Particularly important is the set $\cT^0$ 
of \emph{IC-trees}; it is easy to 
determine whether $T\in \cT$ is an IC-tree: just check whether the 
arrows $a\in \TC$ all satisfy $\cW_T\partial\chi _a= \pi \chi _a$.

To see the relationship to the toric geometry of $X_\zeta $, recall 
 that $X_\zeta$ is the GIT quotient 
of a certain affine variety in $\C^{\cQ_1}$ by $T^{0,0}$. The extreme 
points of $C_\zeta $ therefore correspond to elements $\alphatw 
\in\C^{\cQ_1}$ which are mapped to the same $T^{0,0}$-orbit under the 
morphism of algebraic tori $$\widehat{\cW_S}\colon \CX^n\to T^1.$$ This 
is the case for any $\alphatw $ such that $\supp\alphatw \subseteq S^{\Box}$, where $S^{\Box}$ is the set of all arrows in $\cQ_1$ 
which satisfy~\eqref{eq:S_inv} (this includes of course $\clos S$). 

This discussion can be extended to higher dimensional faces, although 
it is not as practical for concrete calculations.

\subsubsection{Higher Dimensional Faces}
\label{sec:2:examples:weight:higher}
 The sub-lattice of integral flows on $\cQ$ which are zero outside $S$ 
will be denoted $\Z^S\subseteq\Z^{\cQ_1}$. Recall that $Z_0(S)$ is the 
set of $0$-flows which are zero outside $S$: $Z_0(S)= \Z^S\cap 
\ker\partial$. Define, for each $S\subseteq\cQ_1$, a morphism of 
lattices 
$$ \cW_S\colon \partial \Z^S\to{\pi \Z^S/\pi Z_0(S)}$$
defined by
\begin{equation}
\label{eq:ws}
{\partial\chi_a}:= \pi Z_0(S) + \pi\chi_a,\quad \rlap{$a\in S$.} 
\end{equation}
This fits into a commutative diagram
$$\begin{CD}
\label{CD:ws}
\partial\Z^S & @>{\cW_S}>> & {\pi\Z^S/\pi Z_0(S)}\\
@A{\partial}AA & & @AA{\pr}A\\
\Z^S & @>>{\pi_{|\Z^S}}> & \pi\Z^S.
\end{CD}
$$
 Taking the image of this under the functor $\ \widehat{\ } = \Hom(\ 
 \cdot\ ,\C^*)$, one obtains a corresponding diagram of algebraic tori:
 $$
\begin{CD}
\label{CD:tori}
T^{0,0}&\supseteq  T^{\partial\Z^S} 
           & @<{\widehat {\cW_S}}<< & T^{\pi\Z^S/\pi Z_0(S)}&\\
& @V{\widehat\partial}VV & & @VV{\widehat\pr}V &\\
T^1 & \supseteq T^{\Z^S} & @<<{\widehat {\pi_{|\Z^S}}}< & T^{\pi\Z^S}& 
\subseteq \CX^n .
\end{CD}
$$
This diagram shows that the action of the sub-torus $ T^{\pi \Z^S/\pi Z_0(S)}$ of 
$T^{\pi \Z^S}$   on $\alphatw \in \C^S$ via 
$\widehat{\pi_{|\Z^S}}\widehat{\pr}$ leaves $\alphatw $ in the same orbit of 
$T^{0,0}$. In fact, since the torus corresponding to $\Z^{S^c}$ 
(where $S^c$ denotes the complement of $S$ in $\cQ_1$) acts trivially 
on $\C^S$,  one sees that
$$ T^{\pi\Z^{\cQ_1}/\pi Z_0(S)}$$
acts on $\C^S$ fixing the $T^{0,0}$-orbits and is a 
sub-torus of $T^{\pi \Z^{\cQ_1}} = \CX^n$ of codimension $\rk\pi Z_0(S)$. 

\subsubsection{Existence and Uniqueness of $n$-weightings}
\label{sec:2:examples:weight:exist}

\begin{dfn} 
Two vertices in $v,v'\in \cQ$ are said to be \emph{connected by 
a subset\/} $S\subseteq\cQ_1$ if there is 
a  path $p$ in $S$ whose endpoints are $\{v,v'\}$. A subset $S\subseteq 
\cQ_1$ is called a \emph{spanning set} if it connects any 
two vertices in $\cQ_0$. 
\end{dfn} 

One sees easily that the condition that $S$ be a spanning subset is 
equivalent to the statement $\partial\Z^S = \Lambda ^{0,0}$. Thus if 
$S$ is a spanning subset, any $n$-weighting $\cW_S$ of $S$ is 
completely determined by~\eqref{eq:lambda_def}, and so is unique.
In fact, if $T\subseteq S$ is any spanning tree in $S$, then $T$ 
determines a unique morphism $\cW_T\colon \Lambda ^{0,0}\to\Z^n$, and 
one sees that $S$ is invariant if and only  
 \begin{equation}
 \label{eq:S_inv_wt}
 \cW_T\partial\chi _a = \pi \chi _a, \quad\forall a\in S\setminus T.
 \end{equation}

 The following lemma shows that, for generic values of $\zeta $, one
 can effectively use the condition above to check the invariance of
 configurations.

\begin{lemma}
\label{lemma:generic_flow}
 If $\zeta$ is generic in $\Lambda ^{0,0}_\R = \Lambda 
^{0,0}\otimes_\Z \R$ and $f$ is a $\zeta$-admissible flow then 
$\supp(f)$ is a spanning subset. 
\end{lemma}
\begin{proof}
  If $\supp(f)$ is not a spanning subset then there exists a partition
  of $\cQ_0$ into disjoint non-empty subsets $S_0$ and $S_1$ such that
  no element of $S_0$ is connected to any element of $S_1$. Consider
  the restriction of $f$ to the $S_i$; the previous statement implies
  that $\partial (f_{|S_i})=(\partial f)_{|S_i}$. Since, for any
  $\zeta$-flow $g$, ``the flow is conserved," \ie 
  $\sum_{v\in\cQ_0}(\partial g)_v=0$, one has
\begin{align*}
  \sum_{v\in S_i} \zeta_v &= \sum_{v\in S_i} (\partial f)_v \\ 
&=  \sum_{v\in\cQ_0} \left( \partial (f_{|S_i}) \right)_v \\ 
&=0,
\end{align*}
but this does not happen for a generic $\zeta$ in $\R^{\cQ_0}_0$.
\end{proof}

First,  an example which has some singular cones: the action of the
group of fifth roots of unity on $\C^3$ with weights $1$, $2$ and $3$.
In this case, there are whole cones of values of $\zeta$ for which
$X_\zeta$ is smooth, and others where $X_\zeta $ is singular.

\begin{rmk}[About the computations]
The computations were done using several computer programs.  A Pascal 
program was used to produce a list of all the IC-trees for any action 
of a cyclic group.  Then, for each value of $\zeta$, another Pascal 
program was used to determine which trees were admissible and to work 
out the corresponding flows and extreme points.  Then a Mathematica 
program was run to draw the pictures of the polyhedra and of the 
extreme flows.
\end{rmk}

\subsection{Example: the action $\qsing 1/5(1,2,3)$.}
\label{sec:2:ex:5123}

Recall the action $\frac{1}{ 5}(1,2,3)$ of the group $\mu_5$ of fifth
roots of unity on $\C^3$ with weights $1,2$ and $3$  considered in
example~\ref{ex:5123}.  The McKay quiver in this case
is the regular oriented graph with $5$ vertices and an arrow between
each vertex $v$ and $v-1$, $v-2$ and $v-3 \pmod 5$.

\begin{figure}[htbp]
  \begin{center}
    \leavevmode
\epsfysize= 9cm \epsfbox{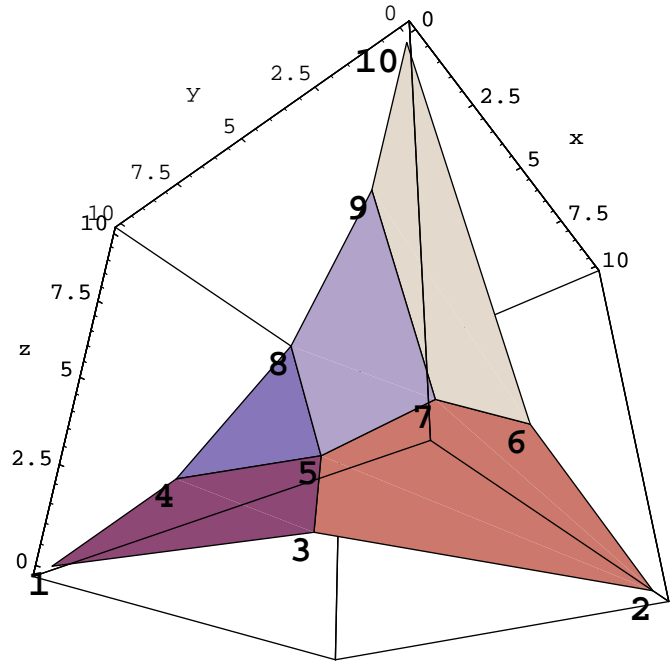}
  \end{center}
  \caption[$C_\zeta$ for $\qsing 1/5(1,2,3)$,
  $\zeta=(-1,-1,-1,-1,4)$.]{$C_\zeta$ for $\qsing 1/5(1,2,3)$,
    $\zeta=(-1,-1,-1,-1,4)$. (The view is from ``behind'', from the
    point with coordinates $(-1.3,-1,-1)$.  The vertices have been
    numbered (where possible) in order of increasing $z$ coordinate.)}
  \label{fig:poly5sing}
\end{figure}

\begin{figure}[htbp]
 \begin{center}
   \leavevmode
\epsfysize= 4cm  \epsfbox{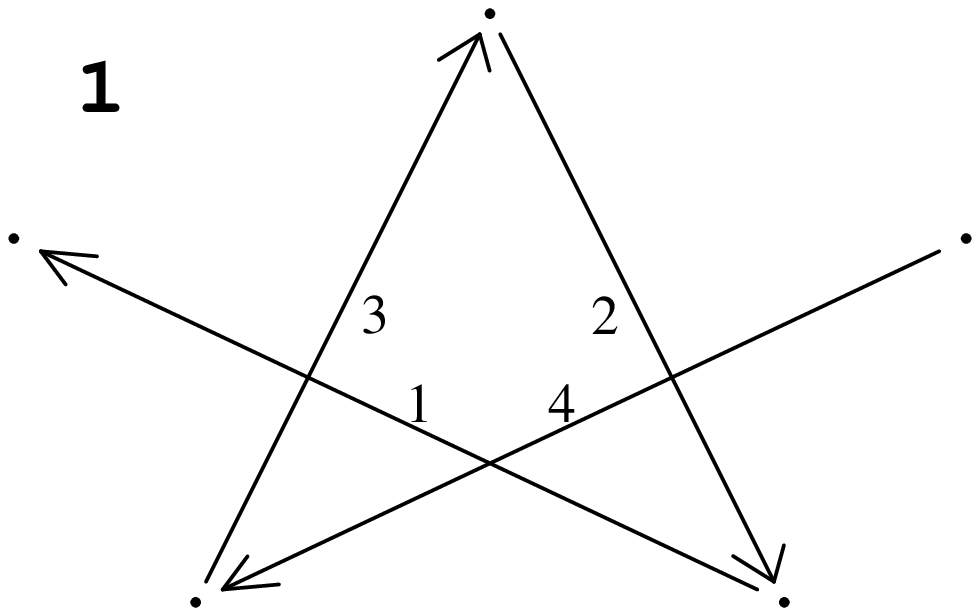}
\epsfysize= 4cm  \epsfbox{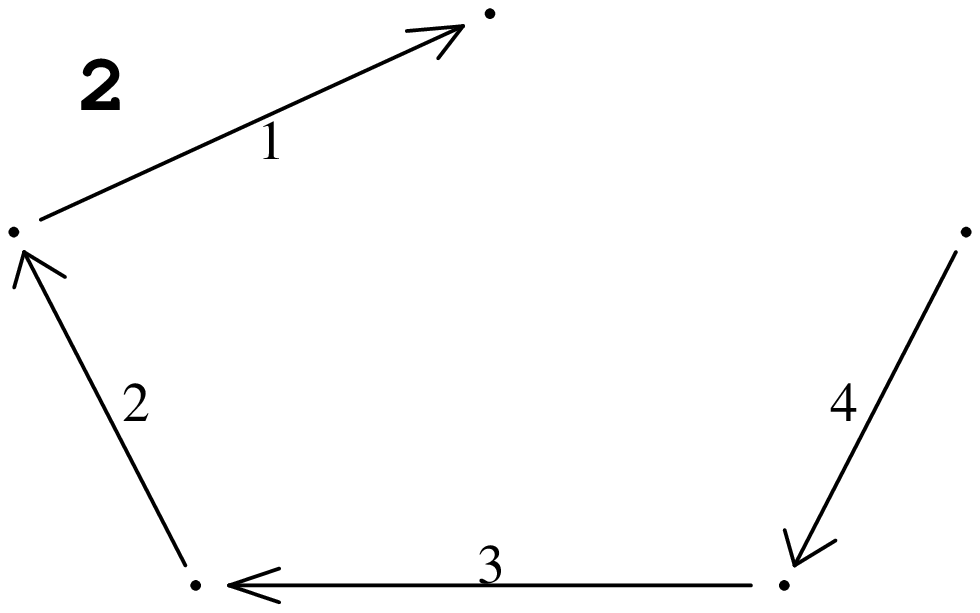}
\epsfysize= 4cm  \epsfbox{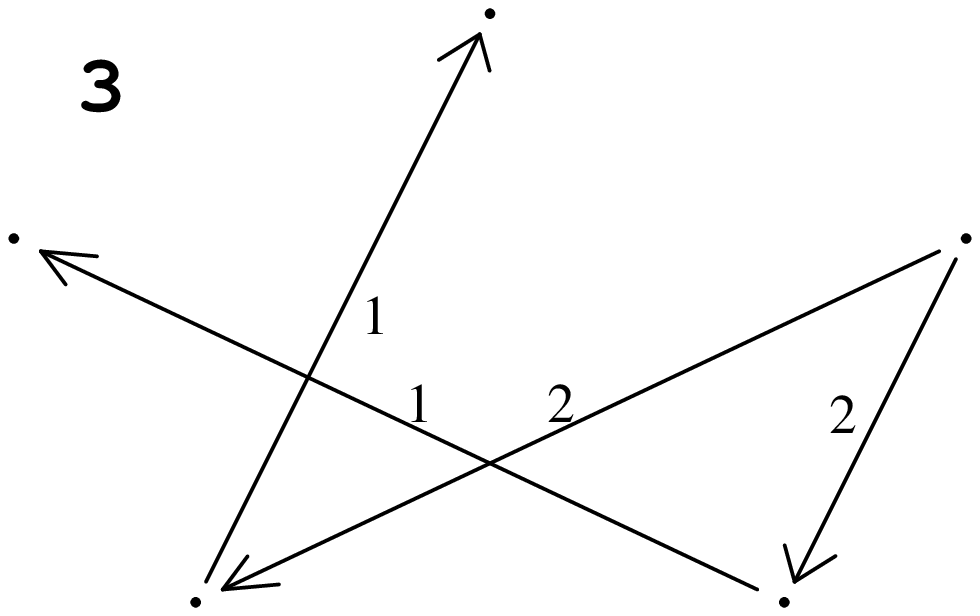}

\vspace*{-1.5cm}
   \leavevmode
\epsfysize= 4cm  \epsfbox{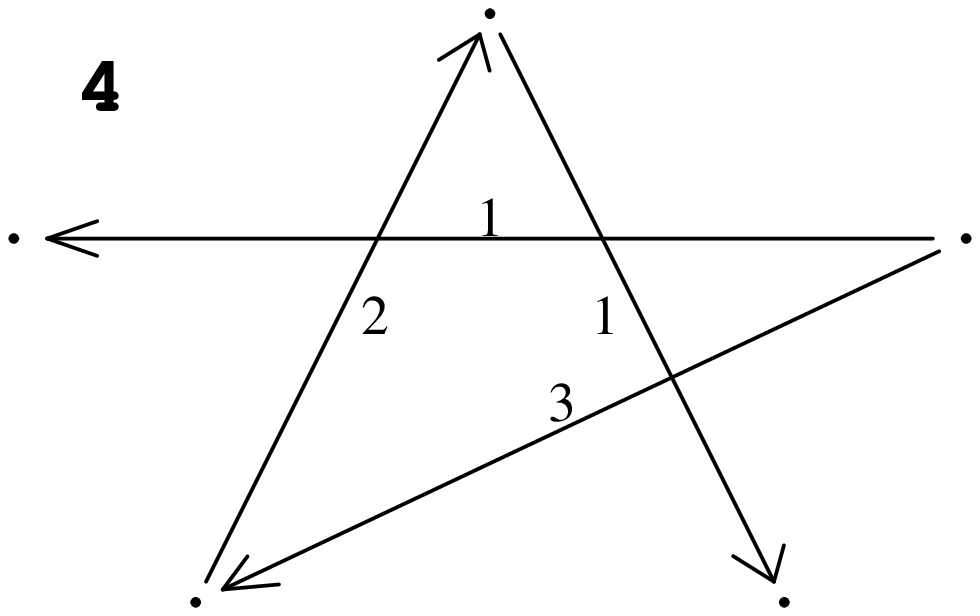}
\epsfysize= 4cm  \epsfbox{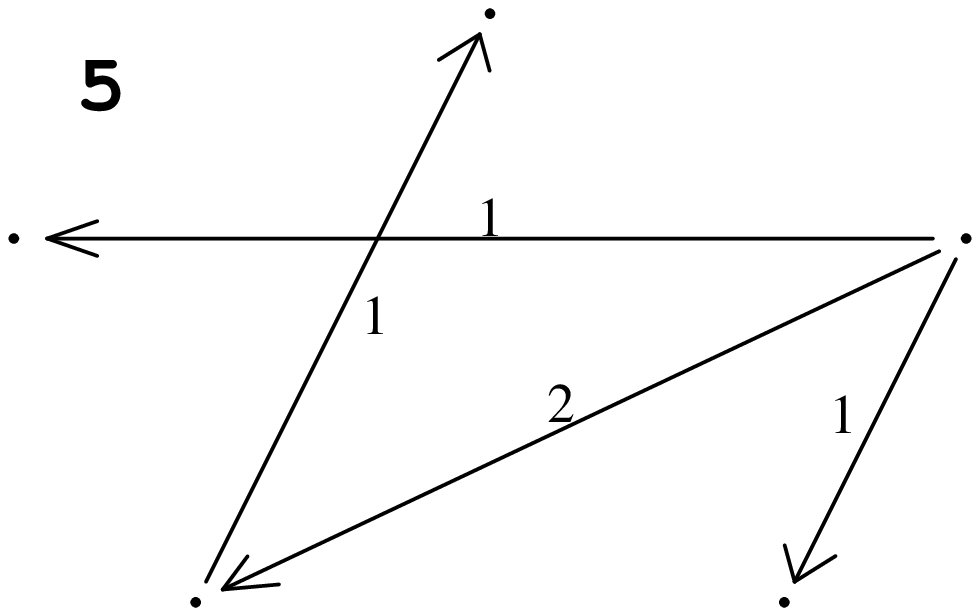}
\epsfysize= 4cm  \epsfbox{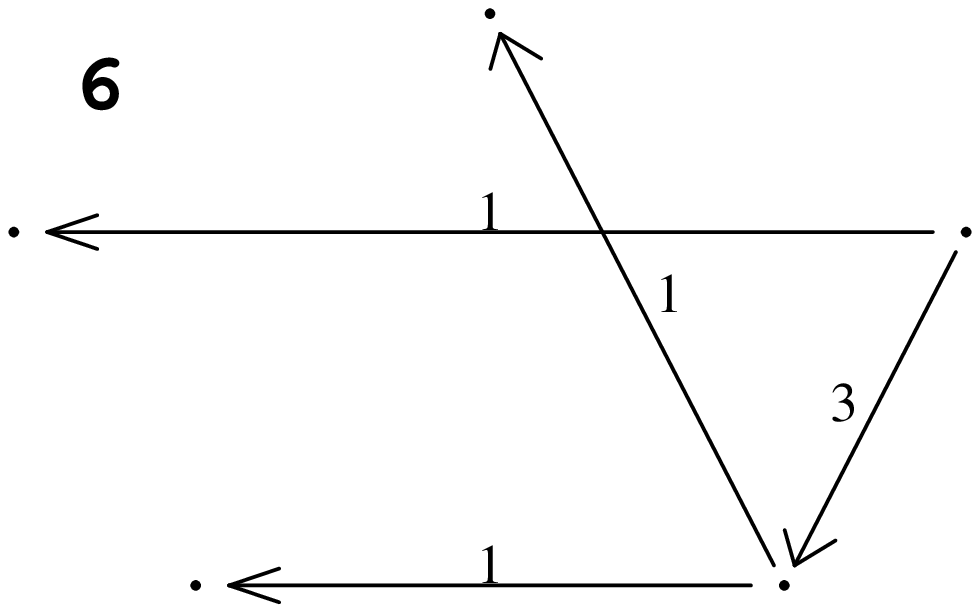}

\vspace*{-1.5cm}
   \leavevmode
\epsfysize= 4cm  \epsfbox{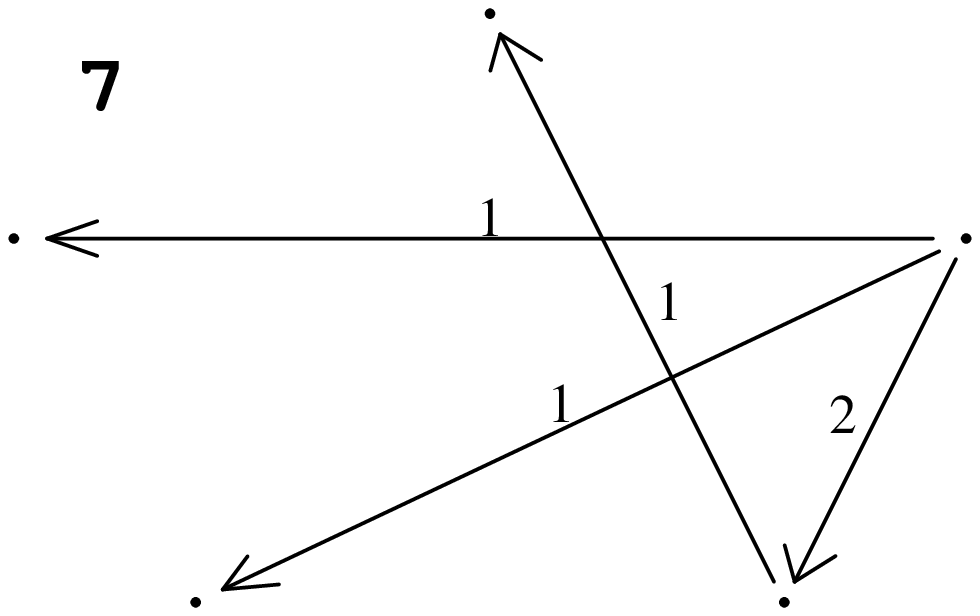}
\epsfysize= 4cm  \epsfbox{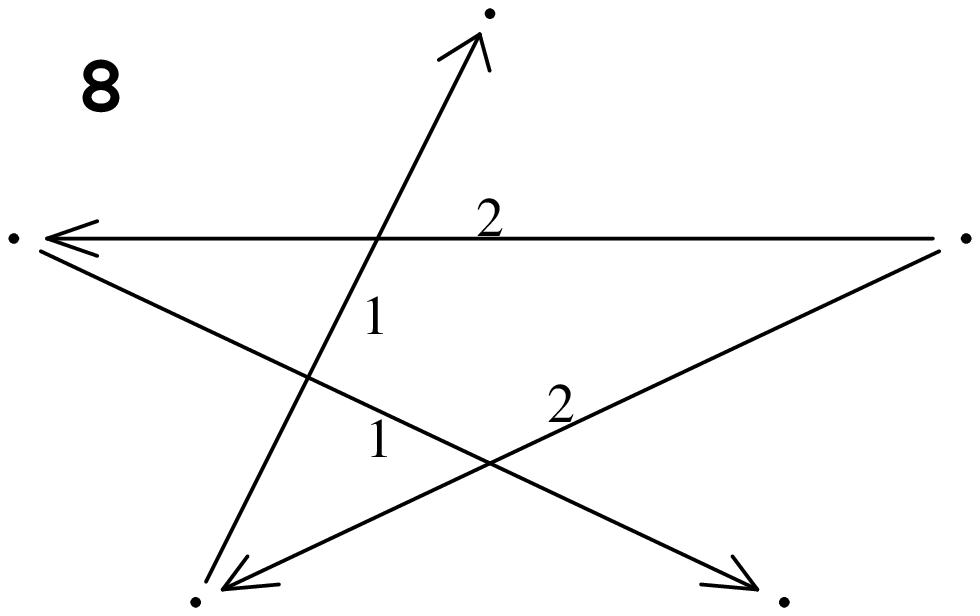}
\epsfysize= 4cm  \epsfbox{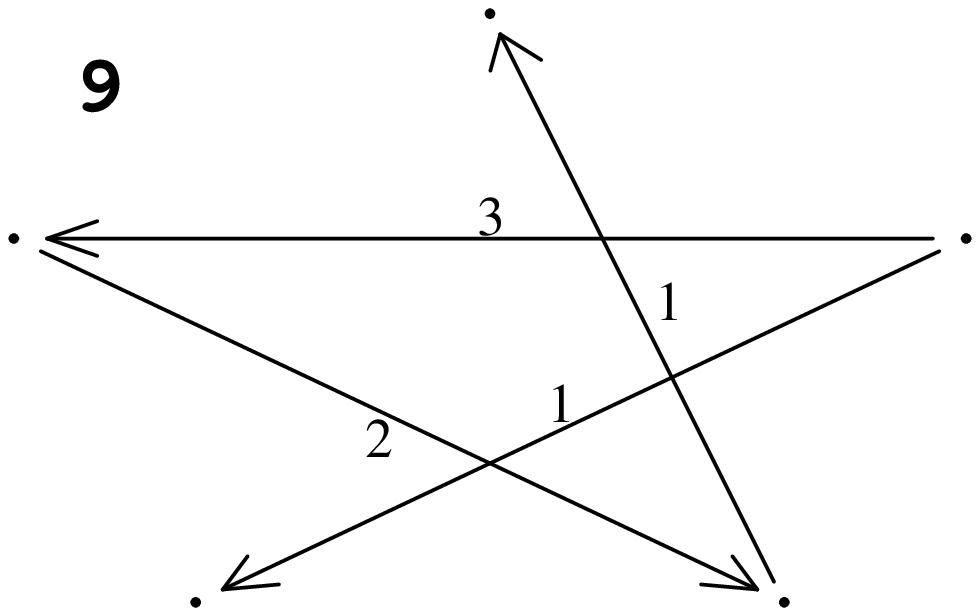}

\vspace*{-1.5cm}
   \leavevmode
\epsfysize= 4cm  \epsfbox{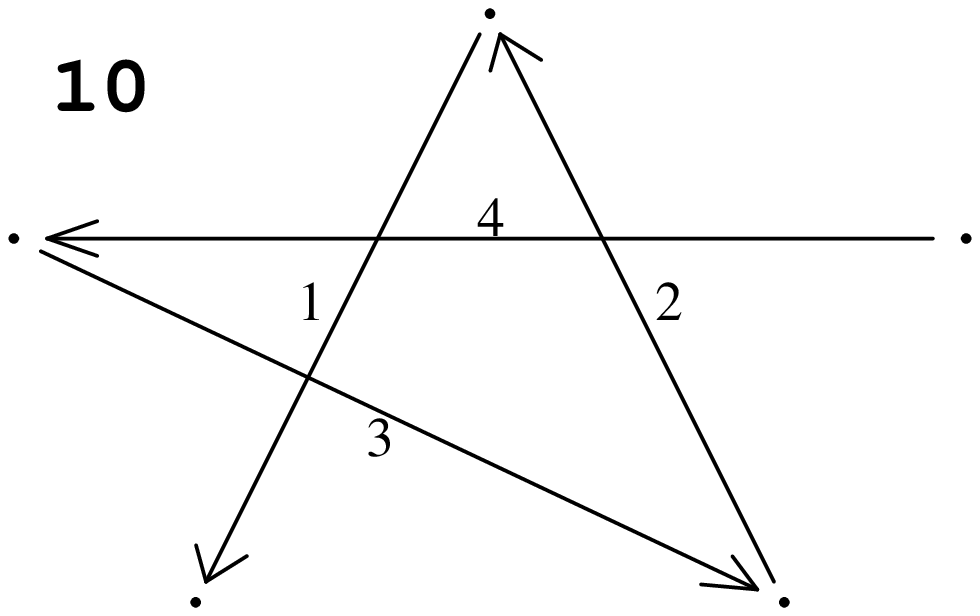}
\end{center}
\vspace*{-1cm}
  \caption[Extreme flows for $\qsing 1/5(1,2,3)$,
  $\zeta=(-1,-1,-1,-1,4)$.]{Extreme flows for $\qsing 1/5(1,2,3)$,
    $\zeta=(-1,-1,-1,-1,4)$.  (The values of the flows are indicated
    to the right of the arrows.  The numbers in the top left-hand
    corners correspond to the vertex numbers in
    Figure~\ref{fig:poly5sing}.)}
  \label{fig:tree5sing}
\end{figure}

The set $\cT^0$ of IC-trees contains a total of 55 trees, once one has
factored out by the symmetry which consists in permuting the vertices
cyclicly (the action of $\gahat$ on itself).  To draw $C_\zeta$,
choose values of $\zeta$, calculate the $\zeta$-flows on all the
IC-trees, discard those which are negative, and project the resulting
points to $\R^3$ via the map $\pi$.  The resulting convex polyhedra
$C_\zeta$ are the intersection of the positive orthant with a finite
number of half-spaces.  The polyhedron $C_\zeta$ for the value
$\zeta=(-1,-1,-1,-1,4)$ is shown in Figure~\ref{fig:poly5sing}.  The
trees and flows corresponding to the extreme points appear in
Figure~\ref{fig:tree5sing}.

One sees immediately from the figure that $X_\zeta$ has a singularity
at the point corresponding to
the extreme point $(1,3,1)$: the tangent cone there has four generators. The
other extreme points are the intersection of  three faces: in order to
determine whether they are in fact singular or not one must check
whether the primitive generators in $\Pi\subset \Z^3$ of the
tangent cone actually generate $\Pi$. In this case it turns out
that they do, so they are smooth points.

\begin{figure}[htbp]
  \begin{center}
\vspace*{-3cm}
    \leavevmode
\epsfbox{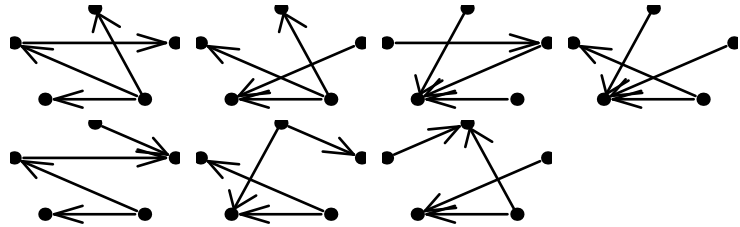}
  \end{center}
\vspace*{-3cm}
  \caption{Non-isomorphic singular IC-trees in $\cT^0$ for the action $\frac{1}{ 5}(1,2,3)$.}
  \label{fig:singtree}
\end{figure}

In fact, there are a total of 7 singular non-isomorphic IC-trees.  
These have been listed in Figure~\ref{fig:singtree}.  The first tree 
(rotated by $\frac{2\pi}{ 5}$) corresponds to the singular point 
$(1,3,1)$ above.  The corresponding tangent cone to $C_\zeta$ was 
described in Section~\ref{sec:2:flow:state} --- it has four 
generators $v_1,\dots,v_4$, any three of which generate $\Pi$, and 
satisfying a single relation of the form $v_1+v_3=v_2+v_4$.  The 
corresponding singularity is therefore of the type $xw=yz 
\subset\C^4$: a cone over a quadric surface.

\begin{figure}[htbp]
  \begin{center}
    \leavevmode
\epsfysize= 9cm \epsfbox{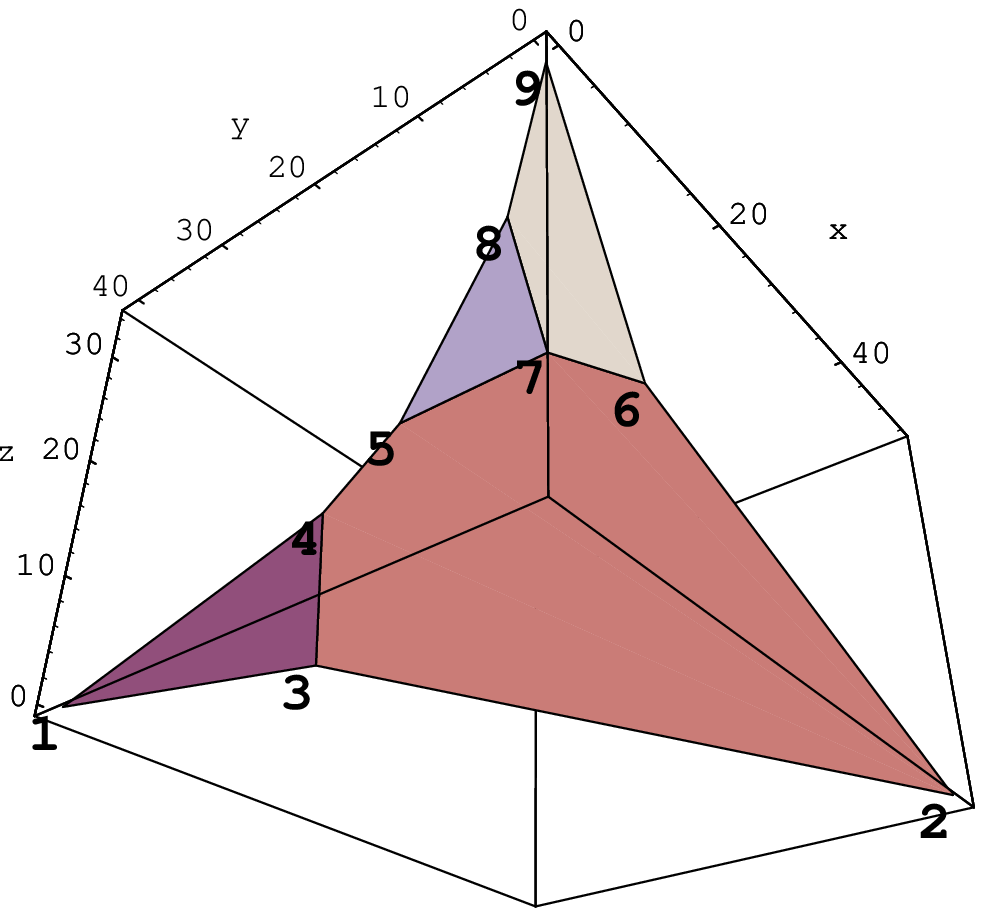}
  \end{center}
%\vspace*{-1cm}
  \caption[$C_\zeta$ for $\qsing 1/5(1,2,3)$,
  $\zeta=(9,8,-3,-2,-12)$.]{$C_\zeta$ for $\qsing 1/5(1,2,3)$,
    $\zeta=(9,8,-3,-2,-12)$.  (The view is from ``behind'', from the
    point with coordinates $(-1.3,-1,-1)$.  The vertices have been
    numbered (where possible) in order of increasing $z$ coordinate.)}
  \label{fig:poly5non}
\end{figure}

\begin{figure}[htbp]
 \begin{center}
   \leavevmode
\epsfysize= 4cm  \epsfbox{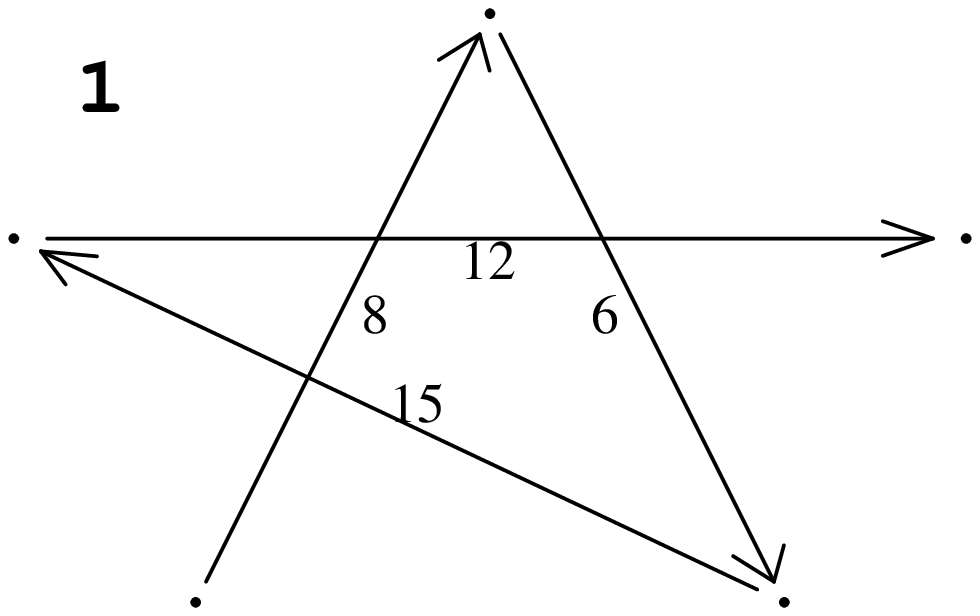}
\epsfysize= 4cm  \epsfbox{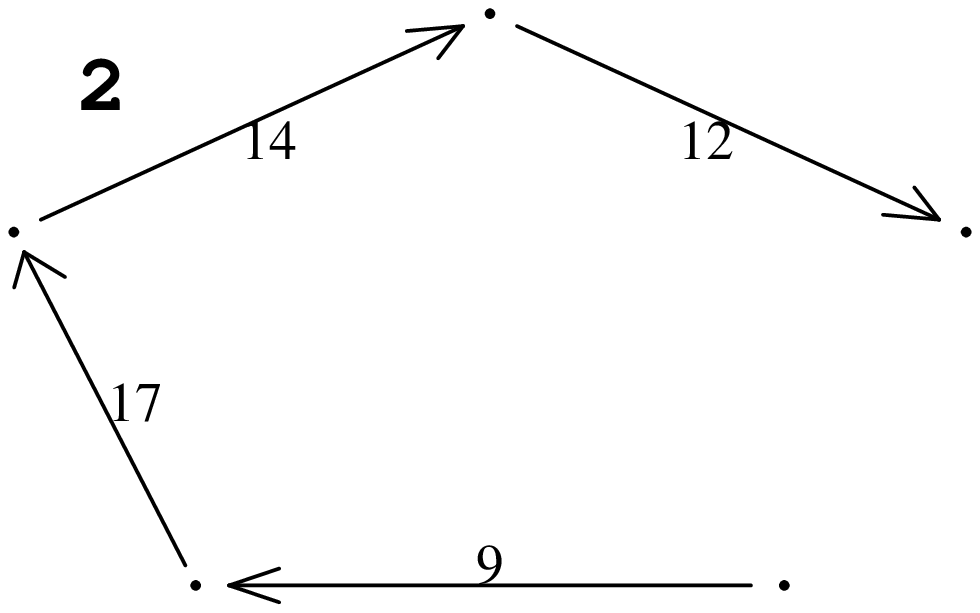}
\epsfysize= 4cm  \epsfbox{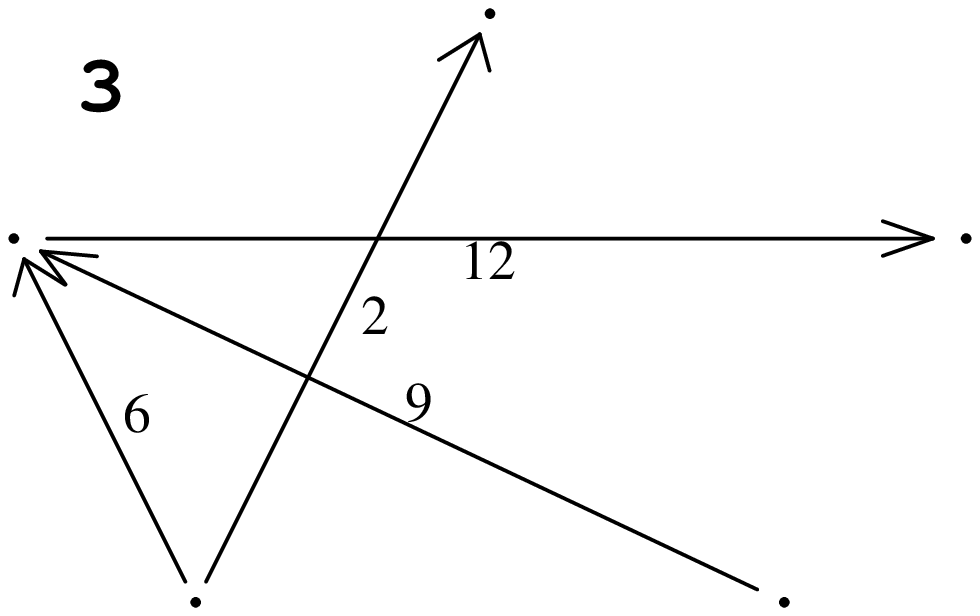}

\vspace*{-1.5cm}
   \leavevmode
\epsfysize= 4cm  \epsfbox{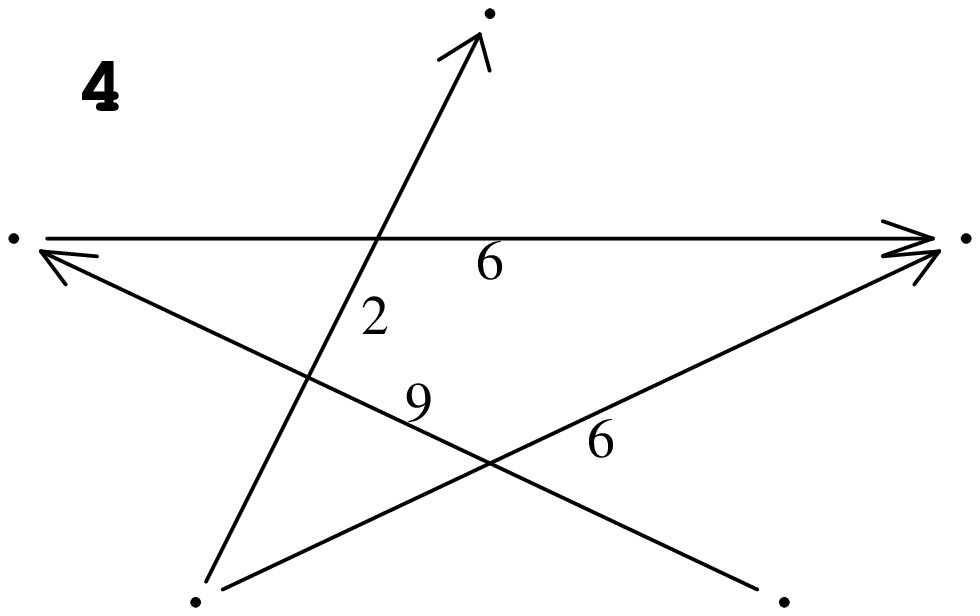}
\epsfysize= 4cm  \epsfbox{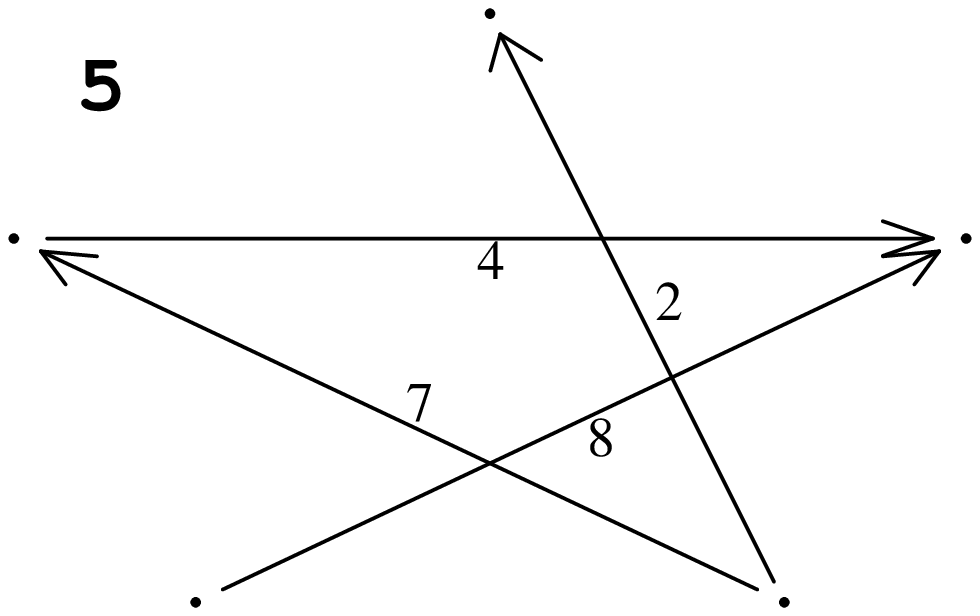}
\epsfysize= 4cm  \epsfbox{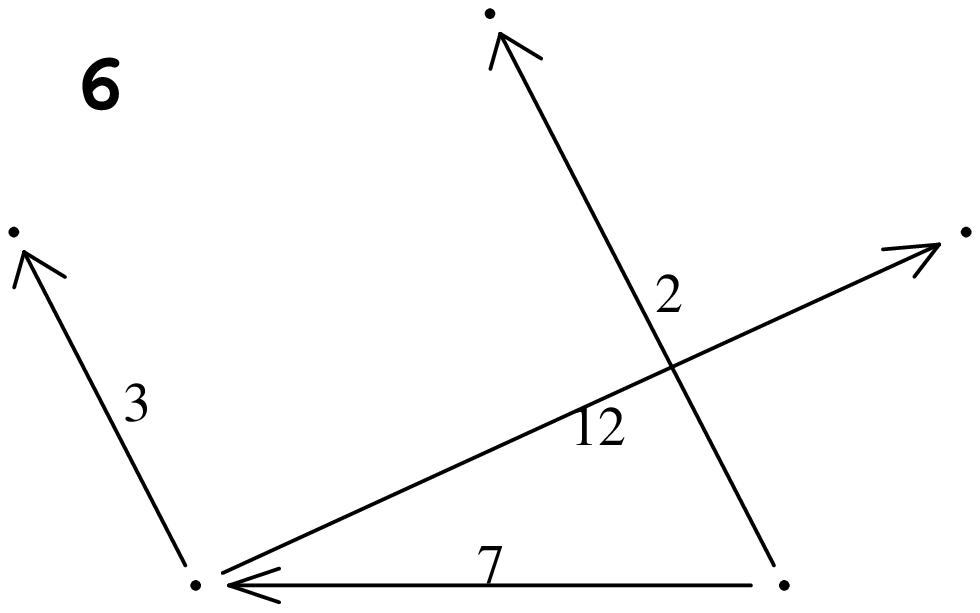}

\vspace*{-1.5cm}
   \leavevmode
\epsfysize= 4cm  \epsfbox{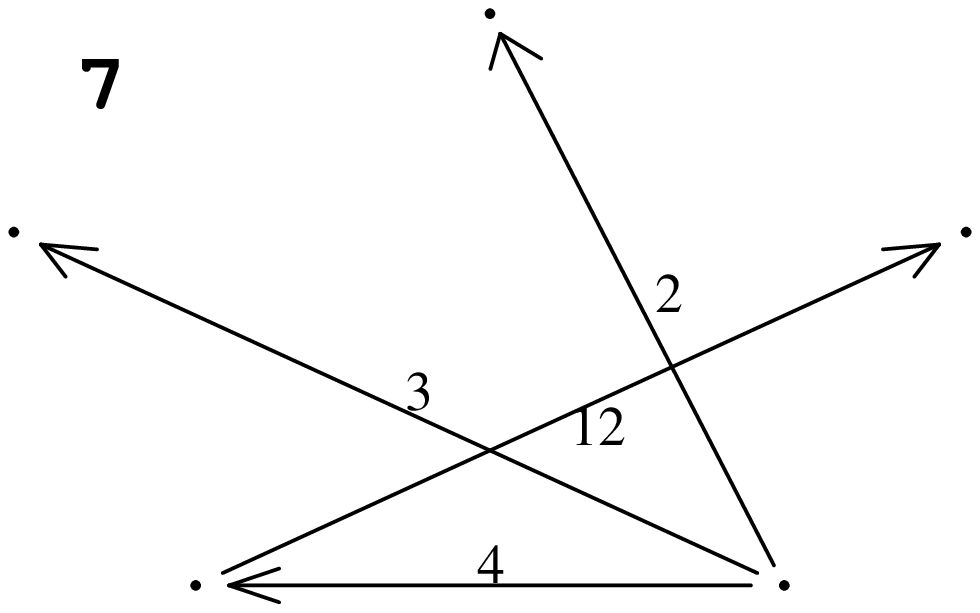}
\epsfysize= 4cm  \epsfbox{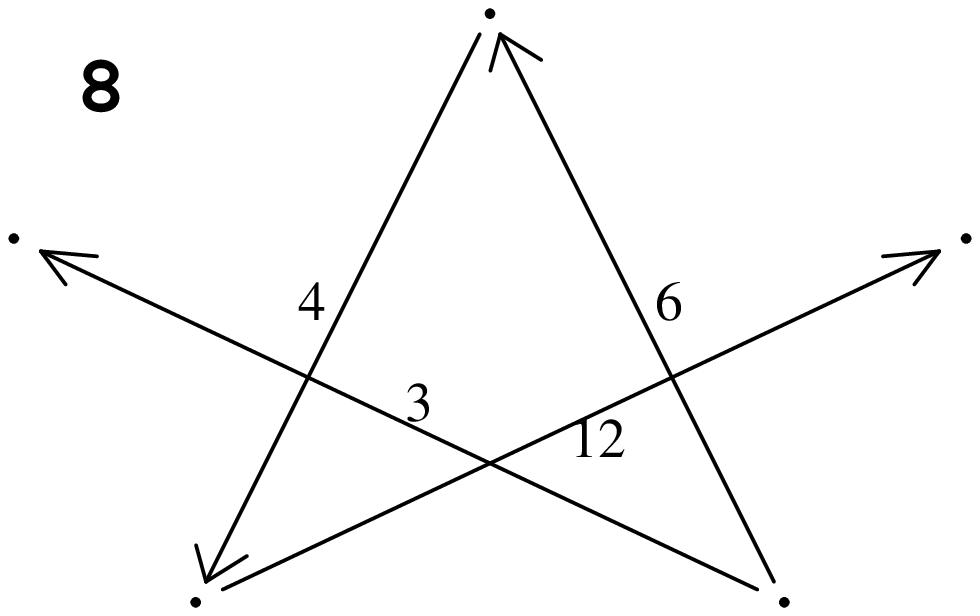}
\epsfysize= 4cm  \epsfbox{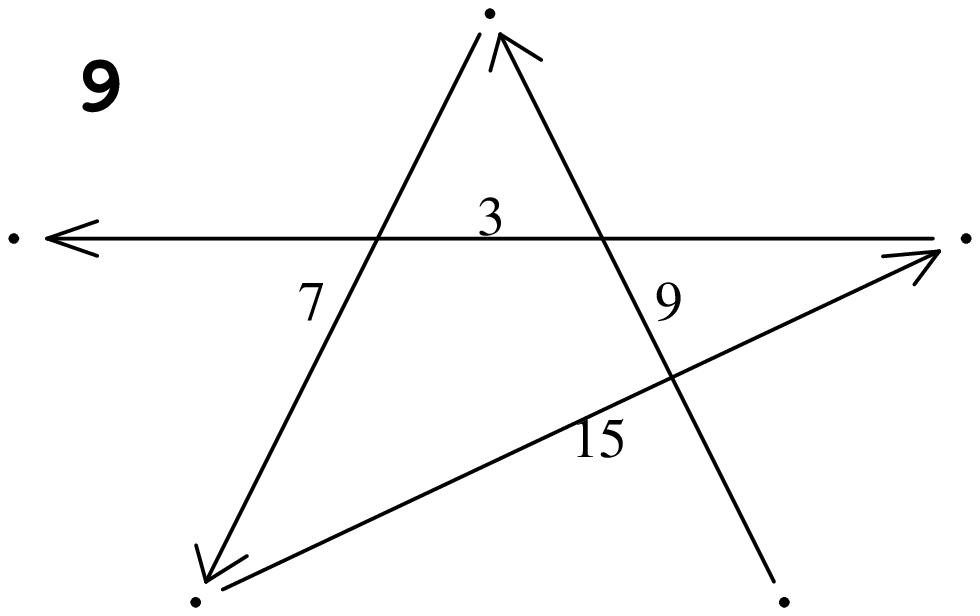}
\end{center}
\vspace*{-1cm}
  \caption[Extreme flows for $\qsing 1/5(1,2,3)$,
    $\zeta=(9,8,-3,-2,-12)$.]{Extreme flows for $\qsing 1/5(1,2,3)$,
    $\zeta=(9,8,-3,-2,-12)$.  (The values of the flows are indicated
    to the right of the arrows.  The numbers in the top left-hand
    corners correspond to the vertex numbers in Figure~\ref{fig:poly5non}.)}
  \label{fig:tree5non}
\end{figure}

Checking the cases listed in Figure~\ref{fig:singtree}, one sees that 
it is possible to find generic values of $\zeta$ for which none of 
these trees (nor their rotations by elements of $\hat\Gamma$) are admissible: for instance, $\zeta=(9,8,-3,-2,-12)$ is 
such a value; this gives a smooth resolution of $\C^3/\Z_5$ which has 
Euler number $9$.  The corresponding polyhedron and flows are shown in 
Figures~\ref{fig:poly5non} and~\ref{fig:tree5non}.

\begin{rmk}
  In fact, all the singular trees in Figure~\ref{fig:singtree}
correspond to the same type of singularity, namely a cone over a
quadric surface.\footnote{Note that this is consistent with the
  conjecture in~\cite{sacha:ale} regarding the
  quadratic nature of the singularities of $X_\zeta$.} All the
IC-trees whose cones are simplicial correspond to non-singular points.
Thus one can tell whether $X_\zeta$ is non-singular simply by checking
whether all extreme points have three edges emanating from them.  This
is also the case for other singularities: for instance, $\qsing
1/6(1,2,4)$, $\qsing 1/7(1,2,5)$, $\qsing 1/8(1,2,6)$, $\qsing 1/9(1,2,7)$; it
is not always the case however: for instance, for $\qsing 1/7(1,2,3)$ and
$\qsing 1/10(1,2,8)$ where there exists elements of $\cT^0$ which
correspond to $\Z_2$-quotient singularities.

\end{rmk}

%\section{-}
\subsection{Crepant Resolutions}
\label{sec:2:crep}

Let $\Sigma_\zeta $ denote the fan determined by the
polyhedron $C_\zeta$. By Corollary~\ref{cor:fan}, the one-skeleton of
$\Sigma_\zeta$ is given by
$$\Sigma_\zeta^{(1)} = \{ \sigma(S)^\vee: S\in\cC^{n-1}_\zeta\}, $$ so $X_\zeta$ has
trivial canonical bundle if the primitive generators $v_S$ of the
cones $\sigma(S)^\vee$ for $S\in\cC^{n-1}$ all belong to the hyper-plane
defined by the equation $\sum n_i = 1$ in $\Pi^* =
\Z^{3}+\frac{\Z}{ r}(w_1,w_2,w_3)$.

\begin{example}
\label{ex:3111}
  Consider the action $\qsing 1/3(1,1,1)$. This has only three
   IC-trees up to isomorphism. They consist of trees with two arrows of
  the same type. For generic values of $\zeta$, a little
  thought shows that the polyhedron $C_\zeta$ is the positive quadrant
  with a small equilateral triangle chopped off.  The dual fan is the
  barycentric subdivision of the positive quadrant by the ray passing
  through the point $v=\frac{1}{3}(1,1,1)\in \Pi^*$.  This is a
  non-singular fan, and since the point $v$ belongs to the hyper-plane
  $\sum n_i = 1$,  this has trivial canonical bundle. The
  variety $X_\zeta$ is the total space of the bundle $\cO(-3)$ over
  $\PP^2$.
\end{example}
To conclude this section,  a picture of a more complicated
example is draw.  Note that  $\ga\subset\SU(3)$ and that the variety 
$X_{\zeta}$ is a smooth crepant resolution with Euler number 11.
\begin{figure}[htbp]
  \begin{center}
    \leavevmode
\epsfysize= 9cm  \epsfbox{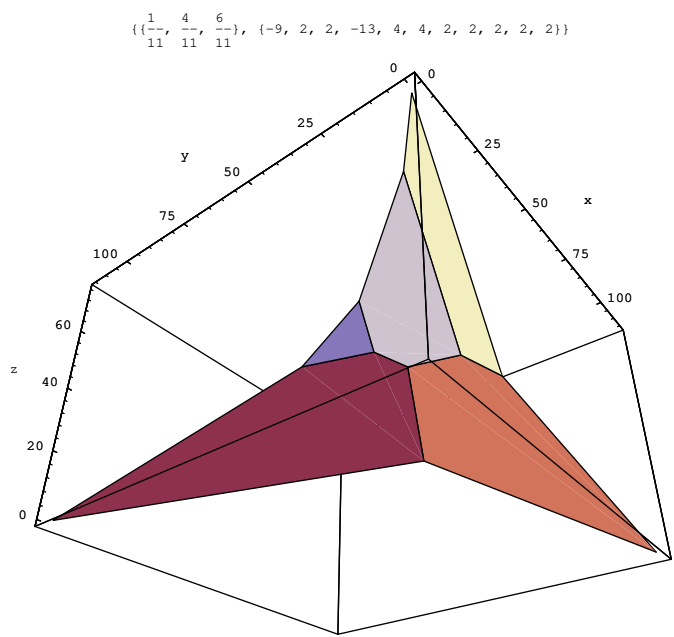}
  \end{center}
  \caption{An example of $C_\zeta$ for the action $\qsing 1/11(1,4,6)$.}
  \label{fig:11_146}
\end{figure}

\providecommand{\bysame}{\leavevmode\hbox to3em{\hrulefill}\thinspace}

\end{document}